\documentclass[11pt]{article}
\usepackage{exscale,relsize}
\usepackage{comment}
\usepackage{lscape}
\usepackage{amsmath}
\usepackage[sort,nocompress]{cite}%increasing order gor citations
\usepackage{amsfonts}
\usepackage[hidelinks]{hyperref}
\usepackage{amssymb}
\usepackage{calc}
\usepackage{theorem}
\usepackage{pifont}      %needed by dingautolist
\usepackage{array}
\usepackage{mathtools} % for ceiling function
\usepackage{cleveref}

\usepackage{color}

\usepackage{graphicx} %Enable figure
\usepackage{float} % Figure positioning
\usepackage{subcaption} % This package includes subfigure command

\newcommand{\barx}{{\bar x}}
\newcommand{\barv}{{\bar v}}

\newcommand{\Ball}{\mathbb{B}}
\newcommand{\sgn}{\ensuremath{\operatorname{sgn}}}
\newcommand{\N}{\mathbb{N}}   % natrual numbers

\newcommand{\norm}[1]{\left\lVert#1\right\rVert} %norm
\newcommand{\ip}[2]{\langle#1,#2\rangle} % inner product

\newcommand{\gph}{\ensuremath{\operatorname{gph}}}

\newcommand{\R}{\ensuremath{\mathbb R}}
\newcommand{\OR}{\ensuremath{\overline{\R}}}
\newcommand{\bR}{{\overline{\R}}}
\newcommand{\Rn}{\ensuremath{\mathbb R^n}}

\newcommand{\RX}{\ensuremath{\,\left]-\infty,+\infty\right]}}

\newcommand{\NN}{\ensuremath{\mathbb N}}

\newcommand{\menge}[2]{\big\{{#1} \mid {#2}\big\}}

\newcommand{\dom}{\ensuremath{\operatorname{dom}}}

\newcommand{\gra}{\ensuremath{\operatorname{gra}}}

\newcommand{\inte}{\ensuremath{\operatorname{int}}}

\newcommand{\card}{\ensuremath{\operatorname{card}}}

\newcommand{\ran}{\ensuremath{\operatorname{ran}}}

\newcommand{\conv}{\ensuremath{\operatorname{conv}}}

\newcommand{\clconv}{\ensuremath{\overline{\operatorname{conv}}\,}}

\newcommand{\Fix}{\ensuremath{\operatorname{Fix}}}

\newcommand{\Id}{\ensuremath{\operatorname{Id}}}

\newcommand{\prox}{\ensuremath{\operatorname{Prox}}}

\newcommand{\scal}[2]{\langle{{#1},{#2}}\rangle}

\newtheorem{theorem}{Theorem}[section]
\newtheorem{lemma}[theorem]{Lemma}
\newtheorem{fact}[theorem]{Fact}
\newtheorem{corollary}[theorem]{Corollary}
\newtheorem{proposition}[theorem]{Proposition}
\newtheorem{defn}[theorem]{Definition}

\theoremstyle{plain}{\theorembodyfont{\rmfamily}
	}
\theoremstyle{plain}{\theorembodyfont{\rmfamily}
	}
\theoremstyle{plain}{\theorembodyfont{\rmfamily}
	}
\theoremstyle{plain}{\theorembodyfont{\rmfamily}
	\newtheorem{example}[theorem]{Example}}
\theoremstyle{plain}{\theorembodyfont{\rmfamily}
	\newtheorem{remark}[theorem]{Remark}}
\theoremstyle{plain}{\theorembodyfont{\rmfamily}
	}
\def\proof{\noindent{\it Proof}. \ignorespaces}
\def\endproof{\ensuremath{\quad \hfill \blacksquare}}

\newcommand{\pluss}{{\hskip1pt \raise1pt\vbox{\hrule width6pt \vskip1pt
			\hrule width6pt}\kern-4pt{\lower1pt\hbox{\vrule height6pt \kern1pt\vrule
				height6pt}}\hskip5pt}}

\newcommand{\argmin}{\mathop{\rm argmin}\limits}

\newcommand{\im}{\mathrm{Im}}

\newcommand{\levp}{\partial_p^\lambda}
\newcommand{\olx}{\ensuremath{\bar{x}}}
\usepackage{tikz}
\usepackage{pgfplots}
\usepackage{tikz}
\usepackage{subcaption}

\usepackage{xcolor}

\begin{document}
	\title{Level proximal subdifferential, variational convexity,
and pointwise quadratic approximation}

\author{Honglin Luo\thanks{School of Mathematical Sciences, Chongqing Normal University, Chongqing, PRC. Email: \texttt{071025013@fudan.edu.cn}},
Xianfu Wang\thanks{Department of Mathematics, I.K. Barber Faculty of Science, University
of British Columbia Okanagan, Kelowna, British Columbia V1V 1V7,
Canada. E-mail: \texttt{shawn.wang@ubc.ca}},
Ziyuan Wang\thanks{Department of Mathematics, I.K. Barber Faculty of Science, University
of British Columbia Okanagan, Kelowna, British Columbia V1V 1V7,
Canada. E-mail: \texttt{ziyuan.wang@ubc.ca}},
and Xinmin Yang\thanks{School of Mathematical Sciences, Chongqing Normal University, Chongqing, PRC. Email: \texttt{xmyang@cqnu.edu.cn}}
}

\date{February 10, 2025; Revision December 12, 2025}

\maketitle

\vskip 8mm

\begin{abstract} \noindent
Level proximal subdifferential was introduced by Rockafellar recently for studying
proximal mappings of possibly nonconvex functions. In this paper a systematic study of level proximal subdifferential is given. We characterize variational convexity of a function by local firm nonexpansiveness of proximal mappings or local relative monotonicity of level proximal subdifferential, and use them to study local convergence of proximal gradient method and others for variationally convex functions. Variational sufficiency guarantees that proximal gradient method converges to local minimizers rather than just critical points.
We also
investigate the existence, single-valuedness and integration
of level proximal subdifferential, and quantify pointwise quadratic approximation (or Lipschitz smoothness) of a function.
As a powerful tool, level proximal subdifferential provides deep insights into variational analysis
and optimization.
\end{abstract}



\section{Introduction}
Proximal mappings are of central importance in both convex optimization
\cite{BC,beck2017,boyd2013, rockafellar1976monotone, moreau65, ryu2022large} and
variational analysis \cite{rockafellar_variational_1998, jourani2014, penot98, poliquin1996, thibault,cwp20,luo2024various}.
In convex optimization, proximal mappings are resolvents of subdifferentials of convex functions and firmly nonexpansive;
these underpin the convergence of various splitting algorithms.
In nonconvex optimization, a major barrier is that
{proximal mappings of nonconvex functions are set-valued in general and no longer resolvents of
their subdifferentials; thus no longer firmly nonexpansive.}
In fact, it is impossible to obtain firmly nonexpansive proximal mappings
 in the absence of convexity, given the following fact: for a proper, lsc and prox-bounded function \(f:\Rn\to\overline{\R}\) one has
\begin{equation}\label{key equivalence}
	\partial f \text{ is (maximally) monotone}
	\Leftrightarrow
	\text{\(f\) is convex}
	\Leftrightarrow
	(\exists\lambda>0)~
	P_\lambda f\text{ is firmly nonexpansive},
\end{equation}
where \(\partial f:\Rn\rightrightarrows\Rn\) is the limiting Mordukhovich subdifferential of \(f\) and
\(P_\lambda f:\Rn\rightrightarrows\Rn\) denotes the proximal operator of \(f\) with parameter \(\lambda>0\); see \cite[Theorem 12.17]{rockafellar_variational_1998} and \cite[Theorem 1]{rock2021}.
Bauschke, Moursi, and Wang investigated the conical nonexpansiveness \cite{bauschke2021generalized} which generalizes firm nonexpansiveness, partially extending the classical framework to nonconvex setting.
But this line of work is restricted to the class of hypoconvex functions (also called weakly convex functions),
see \cite[Proposition 6.4]{bauschke2021generalized}. All of these are in a global
setting.

In 2021 Rockafellar made a breakthrough by introducing the \emph{level proximal subdifferential}\footnote{We note that the level proximal subdifferential first appeared informally in \cite{rock2021},
but no terminology was given there. Hence we adopt the convention used in \cite{wang23level} for convenience.}, see \cite[Equation (2.13)]{rock2021}.
One amazing property of level
proximal subdifferential is that proximal mappings can be written as their resolvents
globally, no matter the functions are convex or nonconvex.
%It has the pleasant feature that every proximal operator is \textit{always} the resolvent of a level proximal subdifferential.
He also showed that \textit{variational convexity} of a prox-regular function at stationary points can be used to characterize firm nonexpansiveness of a localized proximal mapping; see \cite[Theorem 3]{rock2021}.
In contrast, for Mordukhovich limiting subdifferential, this only works for
hyperconvex functions \cite[Theorem 3.1]{wang23level}.
This shows that level proximal
subdifferential has significant advantages in studying proximal mappings of nonconvex functions in general.
To the best of our knowledge,
level proximal subdifferential has been much left unexplored in the literature. One major contribution of our work in this paper is that
in the absence of convexity for prox-regular functions we have
established a new local correspondence among various classical notions in convex optimization\label{convex-varia}:
\begin{figure}[H]
	\centering
	\begin{tikzpicture}[scale=1, every node/.style={scale=0.8}]
		% First tall rectangle with vertically aligned sentences
		\node[draw, minimum width=8cm, minimum height=6cm, align=center] (Rect1) at (0, 0) {
			\textbf{World I (global)}\\
			\\
			$f$ is convex \\
			$\Updownarrow$\\
			$\partial f$ is maximally monotone\\
			%		\parbox{5cm}{\centering $\partial f$ is maximally  \\monotone}\\
			$\Updownarrow$\\
			$e_\lambda f$ is convex \\
			$\Updownarrow$\\
			$\nabla e_\lambda f$ is $\lambda$-cocoercive\\
			$\Updownarrow$\\
			$P_\lambda f$ is nonexpansive\\
			$\Updownarrow$\\
			$P_\lambda f$ is averaged\\
			$\Updownarrow$\\
			%	\parbox{3cm}{\centering $P_\lambda f$ is  firmly \\ nonexpansive}
			$P_\lambda f$ is firmly nonexpansive
		};
		
		% Second tall rectangle with vertically aligned sentences, same dimensions
		\node[draw, minimum width=8cm, minimum height=6cm, align=center] (Rect2) at (8, 0) {
			\textbf{World II (local for prox-regular function)}\\
			\\
			$f$ is variationally convex at $\bar x$ for $\bar v\in\partial f(\bar x)$\\
			%\parbox{3cm}{\centering $f$ is variational \\ convex at $(\bar x,\bar v)$ }\\
			$\Updownarrow$\\
			$\partial_p^\lambda f$ is maximally monotone relative to $W_\lambda$\\
			%\parbox{5cm}{\centering $\partial_p^\lambda f$ is maximal \\monotone relative to $W_\lambda$}\\
			%	$\partial_p^\lambda f$ is maximally monotone\\
			$\Updownarrow$\\
			$e_\lambda f$ is convex on $U_\lambda$ \\
			$\Updownarrow$\\
			$\nabla e_\lambda f$ is $\lambda$-cocoercive on $U_\lambda$\\
			$\Updownarrow$\\
			$P_\lambda f$ is nonexpansive on $U_\lambda$\\
			$\Updownarrow$\\
			$P_\lambda f$ is averaged on $U_\lambda$\\
			$\Updownarrow$\\
			$P_\lambda f$ is firmly nonexpansive on $U_\lambda$
		};
		% Rightarrow between the two rectangles
		\draw[->, thick] (Rect1.east) -- (Rect2.west) node[midway, above] {};
	\end{tikzpicture}
	\caption{Correspondence between convex and variationally convex settings.}\label{fig:diagram}
\end{figure}
To demonstrate the practical significance of our findings, applications to localized proximal gradient method
and localized Krasnosel'ski\u{i}-Mann method for variationally convex optimization are given.

\emph{The goal of this paper is twofold. First, we provide a systematic study of level proximal subdifferential. This includes the existence of level proximal subdifferentials, level proximal subdifferentiable functions, local single-valuedness, and an integration result of level proximal subdifferential.
Second, we explore the connections of level proximal subdifferential to
variational convexity of functions and
their applications in variational convex optimization, and use level proximal subdifferential to
quantify pointwise quadratic approximation (or Lipschitz smoothness).
It turns out that level proximal subdifferential
is a convenient tool not only in variational analysis but also in optimization.
}

The rest of the paper is organized as follows. Section~\ref{s:prelim} contains
some preliminaries
on level proximal subdifferential used in the formulations and proofs of
main results later.
In section~\ref{s:charac} we characterize the existence and single valuedness of level proximal
subdifferential.
We establish that level proximal subdifferential exists at a point if and only if
the function is proximal and its proximal hull is subdifferentiable at the point;
and that
the level proximal subdifferential is a singleton
at a point if and only if the function is proximal and its
proximal hull is differentiable at the point.
%Various conditions for a function to be
%level proximal subdifferentiable are established. One interesting characterization is
%that the function is hypoconvex
%if its level proximal subdifferential is nonempty everywhere in
%the domain of the function and the domain is convex. We also focuses on consequences of single-valued
%level proximal
%subdifferential. One remarkable feature is that the level proximal sudifferential is a singleton
%at a point if and only if the function is proximal and its
%proximal hull is differentiable at the point.
Section~\ref{s:variat} presents new conditions for the variational convexity of a function
in terms of local relative monotonicity of level proximal subdifferential, and
in terms of local nonexpansiveness of proximal mapping. They improve recent important work
by Khan, Mordukhovich and Phat \cite{Khanh2023variational}
on characterizing variational convexity via local convexity of its Moreau
envelope, and by Rockafellar \cite{rock2021} on characterizing variational convexity via firm nonexpansiveness of
its localized proximal mapping. Applications of these results to local proximal gradient method and Krasnosel'ski\u{i}-Mann algorithm for variationally convex functions are given in section~\ref{s:algorithms}. Variational sufficiency guarantees
that these iterative methods converge to a local minimizer rather than just a critical point. Along the way,
we also develop two important calculus for variationally convex functions.
For comparison with section~\ref{s:variat}, we study
single-valued proximal operator on an open set in section~\ref{s:locally}. In essence the single-valuedness of proximal mapping of a function at a point
corresponds to the
strict differentiability of its Moreau envelope at the point.
Section~\ref{s:integ} is devoted to integration of level proximal subdifferentials.
Extending
Benoist and Daniilidis' integration technique for Fenchel subdifferentials \cite{benoist2002,penot2002}
%, we extend
%Rockafellar's integration technique for cyclically monotone operators
to level proximal subdifferentials, we demonstrate that the proximal hull of a function can be recovered.
In the concluding section~\ref{s:lips}, employing
level proximal subdifferential we are able to investigate
pointwise quadratic approximation (or Lipschitz smoothness) of functions. We show that if the pointwise quadractic
approximation property
holds on an open set, it does force the function to have its gradient to be Lipschitz on the open set.
None of the previous research has focused on
pointwise Lipschitz smoothness.

We remark that the paper deals with three topics: level proximal subdifferential, variational convexity, and pointwise quadratic approximation. Sections~\ref{s:prelim}, \ref{s:charac} are devoted to fundamental properties
of level proximal subdifferentials; Sections~\ref{s:variat}, \ref{s:algorithms} use level
proximal subdifferentials to
characterize variational convexity and establish proximal gradient algorithms for
variationally convex optimization problems; Sections~\ref{s:locally}, \ref{s:integ} and \ref{s:lips}
concerning single-valuedness of proximal mapping, integration of level proximal subdifferential, and
the equivalence between pointwise Lipschitz smoothness and nonemptiness of $\partial_{p}^{\lambda}(\pm f)$  are of independent interest, and they
are mostly complementary to results in Sections~\ref{s:prelim}, \ref{s:charac}.
As a whole package, they demonstrate the beauty and very promising applications of level proximal subdifferentials.

The notation that we employ is for the most part standard; however,
a partial list is provided for the readers' convenience.
The Euclidean space $\Rn$ is equipped with the usual scalar product
$\ip{x}{y}:=\sum_{i=1}^{n}\xi_{i}\eta_{i}$ for $x=(\xi_{1},\ldots, \xi_{n})$ and $y=(\eta_{1},\ldots, \eta_{n})$, and induced norm
$\norm{x}:=\sqrt{\ip{x}{x}}$.
The set of extended real numbers is $\overline{\R} :=(-\infty, +\infty]$.
By a proper lsc function $f:\Rn\rightarrow\bR$ we mean that $f$ is
lower semicontinuous on $\Rn$ and
$\dom f:=\{x\in\Rn| \ f(x)<+\infty\}\neq \varnothing$. $\conv f$ and $\clconv f$ are \emph{convex hull} and
\emph{lsc convex hull} of function $f$, respectively; see \cite[pages 56-57]{rockafellar_variational_1998}.
We let $j :=\|\cdot\|^2/2.$
For parameter value $\lambda>0$, the \emph{Moreau envelope} function $e_{\lambda}f$ and \emph{proximal mapping}
$P_{\lambda}f$ of $f$ are defined by
\begin{equation}
e_{\lambda}f(x) :=\inf_{y\in\Rn}\left\{f(y)+\frac{1}{2\lambda}\|y-x\|^2\right\},\text{ and }
P_{\lambda}f(x) :=\argmin_{y\in\Rn}\left\{f(y)+\frac{1}{2\lambda}\|y-x\|^2\right\},
\end{equation}
respectively.
%The $\lambda$-proximal hull of $f$ is $h_{\lambda}f=-e_{\lambda}[-e_{\lambda}f]$.
A function $f:\Rn\rightarrow\bR$ is \emph{prox-bounded} if there exists $\lambda>0$ such that
$e_{\lambda}f(x)>-\infty$ for some $x\in\Rn$. The supremum of the set of
all such $\lambda$ is the threshold $\lambda_{f}$ of prox-boundedness for
$f$. The positivity of $\lambda_{f}$ plays an important role in the paper.
If $\lambda_{f}$ is positive, then $(\forall \lambda\in (0,\lambda_{f}))$ the proximal mapping
$P_{\lambda}f$ is nonempty, compact valued, and outer semicontinuous, and
the Moreau envelope
$e_{\lambda}f$ is locally Lipschitz; see, e.g., \cite[Theorem 1.25, Example 10.32]{rockafellar_variational_1998}.
A proper lsc function $f:\Rn\rightarrow\overline{\R}$ is \emph{$1$-coercive} if
$\liminf_{\|x\|\rightarrow\infty}f(x)/\|x\|=+\infty.$ The \emph{Legendre-Fenchel conjugate} $f^*:\Rn\rightarrow\bR$ of $f$ is defined by
$y\mapsto \sup_{x\in\Rn}\{\ip{y}{x}-f(x)\}$ and $f^{**}=(f^*)^{*}$.
In term of Legendre-Fenchel conjugate, the \emph{$\lambda$-proximal hull} of $f$ defined
by $h_{\lambda}f :=-e_{\lambda}[-e_{\lambda}f]$ can be written as
$h_{\lambda}f=(f+\lambda^{-1}j)^{**}-\lambda^{-1}j$, see
\cite[Example 11.26(c)]{rockafellar_variational_1998}.

Generalized subdifferential constructions from variational analysis
will be widely used in the paper. Let $f:\Rn\to\overline{\R}$ be proper, $x\in\dom f$ and let $\lambda>0$. In \cite{rock2021} Rockafellar defined level proximal subdifferential.
A vector $v\in\Rn$ is a \emph{$\lambda$-level proximal subgradient} of $f$ at $x$, denoted by $v\in\levp f(x)$, if
	\begin{equation}\label{lep ineq}
		(\forall y\in\Rn)~f(y)\geq f(x)+\ip{v}{y-x}-\frac{1}{2\lambda}\norm{y-x}^2.
	\end{equation}
It is a special subclass of proximal subgradients. A vector $v\in\Rn$ is a \emph{proximal subgradient} of $f$ at $x$, denoted by $v\in\partial_{p}f(x)$, if
$(\exists \lambda>0)(\exists \delta>0)(\forall \|y-x\|<\delta) \ f(y)\geq f(x)+\ip{v}{y-x}-\frac{1}{2\lambda}\|y-x\|^2.$
Next is the classical Fenchel subdifferential.
A vector $v\in\Rn$ is a \emph{Fenchel subgradient} of $f$ at $x$, denoted by $v\in\partial_{F}f(x)$, if
$(\forall y\in\Rn)~f(y)\geq f(x)+\ip{v}{y-x}.
$
It is a special class of regular subgradients. A vector $v\in\Rn$ is a \emph{regular subgradient}
of $f$ at $x$, denoted by
$v\in\hat{\partial}f(x)$, if $(\forall y\in\Rn)\ f(y)\geq f(x)+\ip{v}{y-x}+o(\|y-x\|).$
It is a \emph{Mordukhovich limiting (or general) subgradient}, denoted by $v\in\partial f(x)$,
if $\exists v^{\nu}\in\hat{\partial}f(x^{\nu})$ such that $v^{\nu}\rightarrow v$,
$x^{\nu}\rightarrow x$, and $f(x^{\nu})\rightarrow f(x)$.
When $x\not\in\dom f$, all subdiffernetials are set to be empty.
%$\partial_{F}f(x)=\varnothing$.
See the landmark books \cite{mordukhovich2006, rockafellar_variational_1998} for
further details and references on subdifferentials.
We say that $f$ is \emph{prox-regular} at $\bar x$
for $\bar v\in\partial f(\bar x)$ if $(\exists \varepsilon>0, \rho\geq 0)(\forall y\in\Ball_{\varepsilon}(\bar x)) \
f(y)\geq f(x)+\ip{v}{y-x}-\dfrac{\rho}{2}\norm{y-x}^2$ for
$v\in\partial f(x), v\in\Ball_{\varepsilon}(\bar v), x\in\Ball_{\varepsilon}(\bar x),
f(x)<f(\bar x)+\varepsilon$, where $\Ball_{\varepsilon}(a)$ is the open ball centered at $a$ with
radius $\varepsilon$. If this holds for all $\bar v\in\partial f(\bar x)$,
$f$ is said to be prox-regular at
$\bar x$; see \cite{poliquin1996,rockafellar_variational_1998,thibault}.

Next we recall some notions for set-valued mappings.
Let $A:\Rn\rightrightarrows\Rn$ be a set-valued mapping and let $W$ be a nonempty subset of $\Rn\times\Rn$.
Then $A$ is \emph{monotone relative to $W$} if
\(
(\forall (x_1,v_1), (x_2,v_2)\in\gph A\cap W)~
\ip{x_1-x_2}{v_1-v_2}\geq0;
\)
and $A$ is \emph{$\lambda$-strongly monotone relative to $W$} if
\(
(\forall (x_1,v_1), (x_2,v_2)\in\gph A\cap W)~
\ip{x_1-x_2}{v_1-v_2}\geq\lambda\norm{x_1-x_2}^2,
\)
where $\lambda >0$. $A$ is \emph{maximally monotone relative to $W$} if there does not exist
 $A':\Rn\rightrightarrows\Rn$ such that $\gph A'\cap W\supseteq\gph A\cap W$ but
$\gph A'\cap W \neq \gph A\cap W$ and $A'$ is monotone
relative to $W$. See \cite{nghia2016,rock-vietnam,poliquin1996}
for more discussions. The \emph{resolvent}
 of $A$ is
$J_{A} :=(\Id+A)^{-1}$ with $\Id:\Rn\rightarrow\Rn: x\mapsto x$ being the identity mapping.
Let $U$ be a nonempty subset of $\Rn$ and let $T:U \to\Rn$. Then $T$  is \emph{firmly nonexpansive}
 if %$(\forall x,y\in U)$ $\norm{Tx-Ty}^2+\norm{\left(\Id-T\right)x-\left(\Id-T %\right)y}^2\leq\norm{x-y}^2$;
$(\forall x,y\in U)\ \norm{Tx-Ty}^2\leq \ip{x-y}{Tx-Ty}$;
$T$ is \emph{$\lambda$-cocoercive} with $\lambda>0$ if $\lambda T$ is firmly nonexpansive;
$T$ is \emph{nonexpansive} if $(\forall x,y\in U)\ \norm{Tx-Ty}\leq\norm{x-y}$; and $T$ is \emph{averaged} if
there exist nonexpansive $N:U\to\Rn$ and $0\leq \lambda<1$ such that $T:=(1-\lambda)\Id+\lambda N$.
We use $\Fix T$ for the set of fixed point of $T$.
See \cite[Chapter 4]{BC} for abundant results and applications of nonexpansive operators.
For a set $C\subseteq\Rn$, we use $\conv C$ and $\clconv C$ for its convex hull and closed convex hull, respectively.
The set-valued mapping
$\conv(A):\Rn\rightrightarrows\Rn$ is given by $\conv(A)(x):=\conv(A(x))$ for $x\in\Rn$.
The sets
$\dom A :=\{x\in\Rn|\ A(x)\neq \varnothing\},$ $\ran A :=\bigcup_{x\in\Rn}Ax=A(\Rn)$, and $\gph A
:=\{(x,v)\in\Rn\times \Rn|\ v\in Ax\}$ stand for
the domain, range, and graph of $A$, respectively.
%The identity mapping is $\Id:\Rn\rightarrow\Rn: x\mapsto x$.

\section{Preliminaries and auxiliary results}\label{s:prelim}
In this section, we provide
basic relationships of level proximal subdifferential with
well-known subdifferentials and collect
some key properties used in the sequel. See \cite{rock2021, wang23level} for more details.
%\subsection{Basics of level proximal subdifferential}

We begin with a few useful facts concerning level proximal subdifferentials.
\begin{fact}\emph{\cite[Proposition 3.6]{wang23level}}\label{f:levelsub:f}
Let $f:\Rn\to\overline{\R}$ be proper, lsc, and let $x\in\dom f$. Then the following hold:
\begin{enumerate}
\item\label{i:f1}
$\bigcap_{ \lambda>0}\partial_{p}^{\lambda}f(x)=\partial_{F}f(x)$;
\item\label{i:f2} If $f$ is prox-bounded, then
$\bigcup_{\lambda>0}\partial_{p}^{\lambda}f(x)=\partial_{p}f(x)$.
\end{enumerate}
\end{fact}

\begin{fact}\emph{\cite[Proposition 3.4]{wang23level}}\label{fact: level and Fenchel sub}
	Let $f:\Rn\to\overline{\R}$ be proper, lsc and prox-bounded with threshold $\lambda_{f}>0$.
Then
	$(\forall 0<\lambda<\lambda_{f})\
\partial_p^\lambda f=\partial_F\left(f+\lambda^{-1}j\right)-\Id/\lambda$.
\end{fact}
In particular,
Fact~\ref{fact: level and Fenchel sub} implies that $\partial_{p}^{\lambda}f$ is always
$1/\lambda$-hypomonotone.
\begin{fact}\label{fact: relations}
	\emph{\cite[Theorem 3.9]{wang23level}}
	Let $f:\Rn\to\overline{\R}$ be proper, lsc and prox-bounded with threshold
$\lambda_f>0$. Then the following hold:
	\begin{enumerate}
		\item\label{hypoconvex} Let $0<\lambda<\lambda_f$. Then $f$ is $1/\lambda$-hypoconvex $\Leftrightarrow\levp f=\partial f$ $\Leftrightarrow\levp f=\partial_p f\Leftrightarrow\partial_p^\lambda f=\hat\partial f$;
 %Consequently, $\levp f=\partial_p f=\partial f$, provided that $f$ is $1/\lambda$-hypoconvex.
		%If $(\exists \lambda>0)$ $\levp f=\partial f$, then $f$ is $1/\lambda$-hypoconvex.
		\item\label{convex} Suppose that $\lambda_f=\infty$. Then $f$ is convex $\Leftrightarrow (\forall \lambda>0)$ $\levp f=\partial f$ $\Leftrightarrow (\forall \lambda>0)$ $\levp f=\partial_p f$  $\Leftrightarrow (\forall \lambda>0)$ $\levp f=\hat\partial f$.
	\end{enumerate}
\end{fact}
%The following identity is immediate by the definition of level proximal subdifferential.
%\begin{lemma}\label{lem: scalar rule}
%	\(
%	(\forall \lambda>0)
%	\)
%	\(\partial_p^\lambda f=\lambda^{-1}\partial_p^1(\lambda f)\).
%\end{lemma}
One surprising property of $\partial_{p}^\lambda f$ is its connection to the proximal mapping
$P_{\lambda}f$,
no matter $f$ is convex or nonconvex.
\begin{fact}
%[resolvent representation of proximal mapping]
\label{fact: prox identity}
	\emph{\cite[Theorem 3.7]{wang23level}}
	Let $f:\Rn\to\overline{\R}$ be proper, lsc and prox-bounded with threshold
$\lambda_f>0$. Then
	\begin{equation}\label{prox and prox sub}
		(\forall 0<\lambda<\lambda_f)~P_\lambda  f =\big(\Id+\lambda\levp f\big)^{-1}.
	\end{equation}
Equivalently,
$ (\forall x\in\Rn)\ u\in P_{\lambda}f(x)\Leftrightarrow (u, \lambda^{-1}(x-u))\in\gph \partial_{p}^{\lambda}f$.
\begin{comment}
		Consequently, operator $\left(\Id+\lambda\levp f\right)^{-1}:\Rn\to2^{\Rn}$ has full domain, is compact-valued, and for every $x\in\Rn$
	\begin{align*}
		\left(\Id+\lambda\partial_p^\lambda f\right)^{-1}(x)=\big\{v\in\Rn:v_k\to v,\text{where }(\exists x_k\to x)(\exists \lambda_k\to\lambda)v_k\in\left(\Id+\lambda_k\partial_p^{\lambda_k} f\right)^{-1}(x_k)\big\}.
	\end{align*}
	Furthermore $\gph\levp f=T\left(\gph P_\lambda f\right)$, where
	$T=\begin{bmatrix}
		0&\Id\\\Id/\lambda&-\Id/\lambda
	\end{bmatrix}$ is an invertible linear transformation.
\end{comment}
\end{fact}
\begin{corollary}
Let $f:\Rn\to\overline{\R}$ be proper, lsc and prox-bounded with threshold
$\lambda_f>0$. Then
%for every $0<\lambda <\lambda_{f}$, it holds that
\begin{equation}\label{e:prox:e1}
(\forall 0<\lambda <\lambda_{f})\  P_{\lambda}f=[\partial_{F}(\lambda f+j)]^{-1}.
\end{equation}
Consequently, if $\exists 0<\lambda_{0}<\lambda_{f}$
such that $f+\lambda_0^{-1}j$ is convex, then $(\forall 0<\lambda\leq\lambda_{0})\ P_{\lambda}f=(\Id+\lambda \partial f)^{-1}$.
\end{corollary}
\proof
Combine Fact~\ref{fact: level and Fenchel sub} and Fact~\ref{fact: prox identity} to derive
\eqref{e:prox:e1}.
When $f+\lambda_0^{-1}j$ is convex, the function $f+\lambda^{-1}j$ is convex for $0<\lambda\leq\lambda_{0}$,
so $\partial_{p}^{\lambda}f=\partial f$ by Fact~\ref{fact: relations}\ref{hypoconvex}.
The result then follows from Fact~\ref{fact: prox identity}.
\endproof

Next we consider the domain of level proximal subdifferential.
\begin{proposition}\label{thm: always nonempty}
	Let $f:\Rn\to\overline{\R}$ be proper, lsc and prox-bounded with threshold $\lambda_f>0$.
	Then the following hold:
\begin{enumerate}
\item\label{i:prox-level1}
$(\forall 0<\lambda<\lambda_f)$ $\dom\partial_p^\lambda f\neq\varnothing$;
\item\label{i:prox-level2}
\begin{equation}\label{e:upperbound}
\lambda_{f}=\sup\menge{\lambda>0}{\dom\partial_{p}^{\lambda}f\neq\varnothing}.
\end{equation}
\end{enumerate}
\end{proposition}
\proof
\ref{i:prox-level1}:
Note that $(\forall 0<\lambda<\lambda_f)$
$
\ran \left(P_\lambda f\right)=\ran \left(\Id+\lambda\partial_p^\lambda f \right)^{-1}=\dom\left(\Id+\lambda\partial_p^\lambda f \right)=\dom\partial_p^\lambda f,
$
where the first equality owes to Fact~\ref{fact: prox identity}.
Appealing to \cite[Theorem 1.17]{rockafellar_variational_1998}  yields that $\dom P_\lambda f=\Rn$.
Hence $\dom\partial_p^\lambda f=\ran\left(P_\lambda f \right)\neq\emptyset$.

\ref{i:prox-level2}: Put $\gamma:=\sup\menge{\lambda>0}{\dom\partial_{p}^{\lambda}f\neq\varnothing}.$
By \ref{i:prox-level1}, $\lambda_{f}\leq \gamma$. If $\lambda_{f}=+\infty$, then done. If $\lambda_{f}<+\infty$ and
$\lambda_{f}<\gamma$, we take $\lambda\in ]\lambda_{f}, \gamma[$. The definition of
$\gamma$ implies $\dom\partial_{p}^{\lambda}f\neq\varnothing$,
say $x_{0}\in \dom\partial_{p}^{\lambda}f$. For $v_{0}\in \partial_{p}^{\lambda}f(x_{0})$, we have
\begin{equation}
		(\forall x\in\Rn)~f(x)\geq f(x_{0})+\ip{v_{0}}{x-x_{0}}-\frac{1}{2\lambda}\norm{x-x_{0}}^2.
	\end{equation}
It follows from \cite[Exercise 1.24]{rockafellar_variational_1998} that $\lambda_{f}\geq \lambda$, which is
a contradiction. Therefore, \eqref{e:upperbound} holds.
\endproof

It is crucial to consider level proximal subdifferential $\partial_p^\lambda f$ when $0<\lambda<\lambda_f$, as the following example illustrates. It also indicates that the supremum in \eqref{e:upperbound} might not be attained.
\begin{example}
	\label{example:prox threshold}
	Define $f:\R\rightarrow\R: x\mapsto e^{-x^2}-c^{-1}j(x)$ where $c>0$ is a constant.
	Then $f$ is prox-bounded with threshold $\lambda_f=c$ and
	$(\forall x\in\R)\ \partial_p^{\lambda_f} f(x)=\varnothing$.
\end{example}
\proof
%\old{
%	Define $r_f=\inf\{r>0: f+rj\text{ is bounded below}\}$.
%	Now we justify \(c^{-1}=r_f>0\), entailing that $\lambda_f=c$ by \cite[Exercise 1.24]{rockafellar_variational_1998}.
%	First, it is trivial to see $c^{-1}\geq r_f$.
%	Suppose now that there were $0<r<c^{-1}$ such that $f+rj$ is bounded below. Then $(\forall x\in\R)$ $(c^{-1}-r)j(x)\leq e^{-x^2}-\inf(f+rj)\leq 1-\inf(f+rj)$, which is absurd  and implies that $r_f\geq c^{-1}>0$.
%}
{
Since $\displaystyle\lim_{|x|\to\infty}f(x)/|x|^2=-1/2c$,
 the prox-threshold of $f$ is $\lambda_f=c$ in
 view of \cite[Exercise 1.24]{rockafellar_variational_1998}.
}
To find $\partial_p^{\lambda _f}f$,
define $g :=f+c^{-1}j$. Simple calculus gives that $
(\forall x\in\R)\ \partial_{F}g(x)=\varnothing$. Hence $(\forall x\in \R)\ \partial_p^{\lambda _f}f(x)=\varnothing$
by Fact~\ref{fact: level and Fenchel sub}.
%We claim that $\clconv g=0$.
%Since $\clconv g<g$,  Theorem~\ref{thm:domain} shows that $(\forall x\in \R)\ \partial_{p}^{\lambda_{f}}f(x)=\varnothing$.
%Indeed, it is trivial to see that $\clconv g\geq0$. Then
%\begin{align*}
%	(\forall x\in\R)~\left(\clconv g\right)(x)
%	&=
%	\inf_
%	{\substack {\lambda\in[0,1],\\ u,v\in\mathbb{R}}}
%	\Big\{
%	\lambda g(u)+(1-\lambda)g(v): x=\lambda u+(1-\lambda)v
%	\Big\}\\
%	&\leq
%	\inf_{s\geq0}
%	\left(
%	\frac{1}{2}g(x-s)+\frac{1}{2}g(x+s)
%	\right)
%	\leq
%	\lim_{s\to\infty}
%	\left(
%	\frac{1}{2}g(x-s)+\frac{1}{2}g(x+s)
%	\right)
%	=0,
%\end{align*}
%justifying the claim.
%\endproof

%\begin{remark}
%When $0<\lambda=\lambda_f<\infty$, $\partial_p^\lambda f$ may be empty everywhere; see Example \ref{example:prox threshold}.
%\end{remark}
%\subsection{Stationary conditions}
We end this section by showing that level proximal subdifferential can distinguish global and local
minimizers.
% a stationary point in the terms of level proximal
%subdifferential, i.e., $0\in\partial_p^\lambda f(x).$
%We start with optimality conditions.

\begin{proposition}[optimality condition]
	\label{thm:optimality}
	Let $f:\Rn\to\overline{\mathbb{R}}$ be proper, lsc, and let $x\in\dom f$.
	Then the following hold:
	\begin{enumerate}
\item\label{i:fixed}
 $(\forall \lambda>0)\ (\partial_{p}^{\lambda}f)^{-1}(0)=\Fix P_{\lambda}f$;
		\item\label{thm:optimality:global} $x$ is a global minimizer of $f$ if and only if $(\forall \lambda>0)~0\in\partial_p^\lambda f(x)$ if and only if $(\forall \lambda>0)$ $x\in P_\lambda f(x)$;
		\item\label{thm:optimality:local} Assume in addition that $f$ is prox-bounded with threshold $\lambda_f>0$.
		If $x$ is a local minimizer of $f$, then $(\exists 0<\overline\lambda<\lambda_f)~(\forall 0<\lambda\leq\overline\lambda)~0\in\partial_p^{\lambda}f(x).$
	\end{enumerate}
\end{proposition}
\proof
\ref{i:fixed}:
 Observe that
$
0\in\partial_p^\lambda f(x)\Leftrightarrow ~x\in \left(\Id+\lambda\partial_p^\lambda f\right)(x)\Leftrightarrow ~ x
\in P_\lambda f(x).
$

\ref{thm:optimality:global}:
$x$ is a global minimizer of $f$ if and only if
$
(\forall y\in\Rn)~f(y)\geq f(x)
\Leftrightarrow
0\in\partial_{F}f(x)=\bigcap_{\lambda>0}\partial_p^\lambda f (x)
$ by Fact~\ref{f:levelsub:f}\ref{i:f1},
from which the first equivalence readily follows. The second equivalence holds by \ref{i:fixed}.

\ref{thm:optimality:local}:
By assumption, $0\in\partial_p f(x)=\bigcup_{\lambda>0}\partial_p^\lambda f(x)$, implying that $(\exists \overline\lambda>0)$ $0\in\partial_p^{\overline\lambda}f(x)$. Noticing that $(\forall 0<\lambda<\overline\lambda)$ $\partial_p^{\overline\lambda}f(x)\subseteq\partial_p^{\lambda}f(x)$ completes the proof.
\endproof

%\begin{proposition}[optimality condition]
%\label{thm:optimality}
%Let $f:\Rn\to\overline{\mathbb{R}}$ be proper lsc, and let $x\in\dom f$.
%Then the following hold:
%\begin{enumerate}
%	\item\label{thm:optimality:global} $x$ is a global minimizer of $f$ if and only if $(\forall \lambda>0)~0\in\partial_p^\lambda f(x)$ if and only if $(\forall \lambda>0)$ $x\in P_\lambda f(x)$;
%	\item\label{thm:optimality:local} Assume in addition that $f$ is prox-bounded with threshold $\lambda_f>0$.
%	If $x$ is a local minimizer of $f$, then $(\exists 0<\overline\lambda<\lambda_f)~(\forall 0<\lambda\leq\overline\lambda)~0\in\partial_p^{\lambda}f(x).$
%\end{enumerate}
%\end{proposition}
%\proof
%\ref{thm:optimality:global}:
%$x$ is a global minimizer of $f$ if and only if
%$
%(\forall y\in\Rn)~f(y)\geq f(x)
%\Leftrightarrow
%0\in\partial_{F}f(x)=\bigcap_{\lambda>0}\partial_p^\lambda f (x),
%$
%from which the first equivalence readily follows. The second holds by observing that
%$
%(\forall \lambda>0)~0\in\partial_p^\lambda f(x)\Leftrightarrow(\forall \lambda>0)~x\in \left(\Id+\lambda\partial_p^\lambda f\right)(x)\Leftrightarrow(\forall \lambda>0)~ x
%\in P_\lambda f(x).
%$
%
%\ref{thm:optimality:local}:
%By assumption, $0\in\partial_p f(x)=\bigcup_{\lambda>0}\partial_p^\lambda f(x)$, implying that $(\exists \overline\lambda>0)$ $0\in\partial_p^{\overline\lambda}f(x)$. Noticing that $(\forall 0<\lambda<\overline\lambda)$ $\partial_p^{\overline\lambda}f(x)\subseteq\partial_p^{\lambda}f(x)$ completes the proof.
%\endproof

\begin{corollary}[local minimizer]
Let $f:\Rn\to\overline{\mathbb{R}}$ be proper, lsc and prox-bounded with threshold $\lambda_f>0$.
Let $x$ be a local but not global minimizer of $f$ and define
$S:=\{\lambda_{f}\geq \lambda>0|\ 0\in\partial_p^\lambda f(x)\}, \text{ and }\overline\lambda :=\sup S.$
Then $\overline\lambda<\infty$,
$
(\forall0<\lambda\leq\overline\lambda)~0\in\partial_p^{\lambda}f(x)
\text{ and }
(\forall \overline{\lambda}<\lambda\leq\lambda_{f})~0\notin\partial_p^\lambda f(x).
$
\end{corollary}
\proof
Suppose the contrary that $\overline\lambda=\infty$.
Then we can take $\lambda_k\in S$ such that
$\lambda_k\to\infty$ as $k\to\infty$.
% and $(\forall k\in\mathbb{N})~\lambda_k\in S$.
In turn $(\forall k\in\mathbb{N})$ $(\forall 0<\lambda<\lambda_k)$ $0\in\partial_p^{\lambda_k}f(x)$, implying that
$(\forall \lambda>0)~0\in\partial_p^\lambda f(x),$
which forces $x$ to be a global minimizer by Proposition~\ref{thm:optimality},
contradicting to our assumption.
Hence $\overline\lambda<\infty$. Now take a sequence  $(\lambda_k)_{k\in\mathbb N}\subseteq S$ such that $\lambda_k\to\overline\lambda$.
Then
$(\forall k\in\mathbb{N})~(\forall y\in\Rn)~
f(y)\geq f(x)-\frac{1}{2\lambda_k}\|y-x\|^2.
$
Passing $k$ to infinity yields $0\in\partial_p^{\overline\lambda} f(x)$.
\endproof

%In terms of Moreau envelopes, we have the following characterizations.
%\begin{theorem}
%%[connection to Moreau envelope]
%Let $f:\Rn\to\overline{\mathbb{R}}$ be proper, lsc and prox-bounded with threshold $\lambda_f>0$. Then $(\forall 0<\lambda<\lambda_f)$
%$
%0\in\partial_p^\lambda f(x)
%\Leftrightarrow
%e_\lambda f(x)=f(x)
%\Leftrightarrow
%x\in P_\lambda f(x).
%$
%\end{theorem}
%\proof
%Suppose first that $0\in\partial_p^\lambda f(x)$.
%Then
%$$
%(\forall y\in\Rn)~f(y)+\frac{1}{2\lambda}\|y-x\|^2\geq f(x)
%\Leftrightarrow
%e_\lambda f(x)\geq f(x),
%$$
%but it always holds that $e_\lambda f(x)\leq f(x)$ by \cite[Example 1.44]{rockafellar_variational_1998}, hence $e_\lambda f(x)=f(x)$.
%Conversely, assume $e_\lambda f(x)=f(x)$. Then
%$
%(\forall y\in\Rn)~
%f(y)+\frac{1}{2\lambda}\|y-x\|^2\geq f(x),
%$
%implying that $0\in\partial_p^\lambda f(x)$.
%The last equivalence follows from the definition of proximal mapping.
%\endproof
%
%\begin{corollary}
%Let $f:\Rn\to\overline{\mathbb{R}}$ be proper, lsc and prox-bounded with threshold $\lambda_f>0$.
%Let $x$ be a local minimizer but not global. Then $\exists 0<\overline\lambda<\lambda_f$ such that
%$$
%(\forall 0<\lambda\leq\overline\lambda)~e_\lambda f(x)=f(x)
%\text{ and }
%(\forall \lambda>\overline\lambda)~e_\lambda f(x)<f(x).
%$$
%\end{corollary}

\section{Level proximal subdifferential: existence and single-valuedness}\label{s:charac}
In this section, for a proper, lsc and prox-bounded function we investigate conditions under which the level proximal subdifferential is nonempty, and conditions under which the level proximal subdifferential is single-valued.
It turns out that the proximal hull plays an important role.
%\subsection{Level proximal subdifferentiable functions}
\subsection{Existence of level proximal subdifferential}

In general, even for locally Lipschitz functions,
$\dom\partial_{p}f$ can be a countable set only, so is $\dom \partial_{p}^{\lambda}f$; see, e.g., \cite{girgensohn98} or\cite[Theorem 6.1]{wolenski95}.
It is natural to ask under what conditions $\partial_{p}^{\lambda}f(x)\neq\varnothing$.
We start with an auxiliary result.
\begin{lemma}\label{lem: proximal and closed cvx hull}
	Let $f:\Rn\to\overline{\mathbb{R}}$ be proper, lsc and prox-bounded with threshold $\lambda_f>0$, and let $0<\lambda<\lambda_f$.
	Then the following hold:
	\begin{enumerate}
		\item\label{lem: closed cvx hull}
		$f+\lambda^{-1}j$ is 1-coercive.
		Consequently $\clconv (f+\lambda^{-1}j)=\conv (f+\lambda^{-1}j)$ is $1$-coercive, and $\dom[\clconv (f+\lambda^{-1}j)]=\conv(\dom f)$;
\item\label{i:convexity}
$h_\lambda f+\lambda^{-1}j=\conv (f+\lambda^{-1}j)$ is a proper, lsc and convex function;
		\item\label{lem:proximal hull}	
		$h_\lambda f(x)=f(x)\Leftrightarrow \conv (f+\lambda^{-1}j)(x)=(f+\lambda^{-1}j)(x).$
	\end{enumerate}
\end{lemma}
\proof
\ref{lem: closed cvx hull}:
Note that the function $(\forall \lambda<\lambda_0<\lambda_f)$ $f+\lambda_0^{-1}j$ is bounded below because $f$ is prox-bounded with threshold $\lambda_f$.
Then $f+\lambda^{-1}j=f+\lambda_0^{-1}j+(\lambda^{-1}-\lambda_0^{-1})j$ is 1-coercive.
By \cite[Corollary 3.47]{rockafellar_variational_1998} or \cite[Lemma 3.3]{benoist1996}, the function $\conv (f+\lambda^{-1}j)$ is proper, lsc, $1$-coercive, and
$\dom \conv (f+\lambda^{-1}j)=\conv \dom f$.

\ref{i:convexity}: Apply \ref{lem: closed cvx hull} and \cite[Example 11.26(c)]{rockafellar_variational_1998} to obtain
\[
h_\lambda f
=
\clconv\left(f+\lambda^{-1}j\right)-\lambda^{-1}j=\conv \left(f+\lambda^{-1}j\right)-\lambda^{-1}j.
\]
%Thus $h_\lambda f+\lambda^{-1}j$ is convex.% and consequently

\ref{lem:proximal hull}:
Again by \cite[Example 11.26(c)]{rockafellar_variational_1998}, we have $h_\lambda f(x)=(f+\lambda^{-1}j)^{**}(x)-\lambda^{-1}j(x).$
Then
$f(x)=h_\lambda f(x)=(f+\lambda^{-1}j)^{**}(x)-\lambda^{-1}j(x)
\Leftrightarrow
\conv(f+\lambda^{-1}j)(x)=\clconv(f+\lambda^{-1}j)(x)=(f+\lambda^{-1}j)^{**}(x)=(f+\lambda^{-1}j)(x)
$.
%\red{(TODO: Example why convex hull equals closed convex hull.)}
\endproof

%Strikingly, it turns out that the $\lambda$-level proximal subdifferentiability of a proper, lsc, and prox-bounded function $f$ is closely tied to properties of its proximal hull $h_\lambda f$ and the majorizing function $f+\lambda^{-1}j$ of $f$.
Below we say that $f$ is \emph{$\lambda$-level proximal subdifferentiable at $x$} if $\partial_p^\lambda f(x)\neq\varnothing$.
We now present the first main result of this section, which unveils an
intimate connection between the $\lambda$-level proximal subdifferentiability
of $f$ at a point and properties of its proximal hull $h_\lambda f$.

\begin{theorem}[level proximal subdifferentiability at a point]
\label{thm:domain}
Let $f:\Rn\to\overline{\mathbb{R}}$ be proper, lsc and prox-bounded with threshold $\lambda_f>0$.
Let $0<\lambda<\lambda_f$ and let $x\in\dom f$.
Then the following are equivalent:
\begin{enumerate}
	\item\label{thm:domain::domain}
%$\partial_p^\lambda f(x)\neq\emptyset$, i.e.,
$f$ is $\lambda$-level proximal subdifferentiable at $x$;
	\item\label{thm:domain::prox hull} $h_\lambda f(x)=f(x)$ and
$x\in\dom \partial h_{\lambda}f$;
	\item\label{thm:domain::cvx hull} $\conv (f+\lambda^{-1}j)(x)=(f+\lambda^{-1}j)(x)$
and $x\in\dom \partial [\conv (f+\lambda^{-1}j)]$.
\end{enumerate}
\end{theorem}
\proof
``\ref{thm:domain::domain}$\Rightarrow$\ref{thm:domain::prox hull}":
Let $u\in\partial_p^\lambda f(x)$.
Then $x\in P_\lambda f(x+\lambda u)$, which implies through \cite[Example 1.44]{rockafellar_variational_1998} yields that $f(x)=h_\lambda f(x)$.
Next we show that $x\in\dom \partial( h_\lambda f+\lambda^{-1}j)$.
	Define $v :=u+x/\lambda$. Then the assumption $u\in\partial_p^\lambda f(x)$ implies
	\begin{align*}
	(\forall y\in\Rn)~
	f(y)&\geq f(x)+\ip{u}{y-x}-\frac{1}{2\lambda}\norm{y-x}^2\\
	&=f(x)+\ip{u}{y-x}-\frac{1}{2\lambda}\norm{y}^2+\frac{1}{\lambda}\ip{x}{y-x}+
\frac{1}{2\lambda}\norm{x}^2\\
	&=
	\left(h_\lambda f+\lambda^{-1}j\right)(x)+\ip{v}{y-x}-\lambda^{-1}j(y)
	\end{align*}
	where the last equality owes to $f(x)=h_\lambda f(x)$.
	Equivalently $(\forall y\in\Rn)$ $f(y)+\lambda^{-1}j(y)\geq \left(h_\lambda f+\lambda^{-1}j\right)(x)+\ip{v}{y-x}$, entailing $y\mapsto \left(h_\lambda f+\lambda^{-1}j\right)(x)+\ip{v}{y-x}$ to be a convex minorant of $f+\lambda^{-1}j$.
	In turn
	$$
	(\forall y\in\Rn)~
	\left(h_\lambda f+\lambda^{-1}j\right)(y)=\clconv\left(f+\lambda^{-1}j\right)(y)\geq \left(h_\lambda f+\lambda^{-1}j\right)(x)+\ip{v}{y-x},
	$$
	thus $v\in \partial(h_{\lambda}f+\lambda^{-1}j)(x)$.
Since $\dom \partial(h_{\lambda}f+\lambda^{-1}j)=\dom \partial h_{\lambda}f$, we have
$x\in \dom \partial h_{\lambda}f$.

``\ref{thm:domain::prox hull}$\Rightarrow$\ref{thm:domain::cvx hull}":
Lemma~\ref{lem: proximal and closed cvx hull}\ref{lem:proximal hull} gives $\conv (f+\lambda^{-1}j)(x)=(f+\lambda^{-1}j)(x)$.
Moreover, $\dom \partial h_{\lambda}f=\dom \partial(h_\lambda f+\lambda^{-1}j)$ and the identity $h_\lambda f+\lambda^{-1}j=\conv(f+\lambda^{-1}j)$
from Lemma~\ref{lem: proximal and closed cvx hull}\ref{i:convexity} entail
that $x\in\dom \partial [\conv (f+\lambda^{-1}j)]$.

 ``\ref{thm:domain::cvx hull}$\Rightarrow$\ref{thm:domain::domain}":
 Let $u\in\partial[\conv(f+\lambda^{-1}j)](x)$. Applying the subgradient inequality yields that
\begin{align*}
(\forall y\in\Rn)~(f+\lambda^{-1}j)(y)
&\geq\conv (f+\lambda^{-1}j)(y)\geq\conv (f+\lambda^{-1}j)(x)+\langle u,y-x\rangle\\
&=(f+\lambda^{-1}j)(x)+\langle u,y-x\rangle,
\end{align*} implying that $u\in\partial_F\left(f+\lambda^{-1}j\right)(x)$. Put $v :=u-x/\lambda$.
Then $v\in\partial_p^\lambda f(x)$ by Fact~\ref{fact: level and Fenchel sub}.
\endproof

%\begin{remark}\label{rem:domain::prox-threshold}
%It is easy to see that Theorem~\ref{thm:domain} still holds with $\lambda=\lambda_f<\infty$.
%\end{remark}

The condition $x\in \dom\partial h_{\lambda}f$ in Theorem~\ref{thm:domain} is necessary,
as the example below demonstrates.
\begin{example} Let $\lambda>0$.
%[necessity of $x\in\dom \partial (\conv (f+\lambda^{-1}j))$]
Consider the function
\begin{align*}
f:\R\rightarrow\bR: x\mapsto
\begin{cases}
	x\ln(x), &\text{ if }x >0;\\
0, &\text{ if } x=0;\\
	+\infty,&\text{ otherwise}.
\end{cases}	
\end{align*}
Then $f$ is convex, implying that $h_\lambda f=f$ and
$\partial_p^\lambda f=\partial f$ by Fact~\ref{fact: relations}\ref{convex}.
However, one has $\partial_p^\lambda f(0)=\partial f(0)=\varnothing.$
\end{example}
One always has
\begin{equation}
(\forall \lambda>0)\ \dom \partial_{F}f\subset \dom \partial_{p}^{\lambda}f\subset \dom \partial_{p}f\subset \dom \hat{\partial}f\subset
\dom\partial f.
\end{equation}
For nonconvex functions, one should not expect $\dom\partial_{p}^{\lambda}f$ `large'.
In general, $\dom\partial_{p}^{\lambda}f$ is even a proper subset of $\dom\partial f$.
In fact,
if $\dom\partial_{p}^{\lambda}f=\dom f$ and $\dom f$ is convex, the function $f$ has to be hypoconvex,
as we proceed to show. To see this, we need:

%Having established pointwise equivalency of level proximal subdifferentiability in Theorem \ref{thm:domain}, we now work towards ``global'' characterizations.

\begin{lemma}\label{lem:prox of proximal hull}
	Let $f:\Rn\to\overline{\R}$ be proper, lsc and prox-bounded with threshold $\lambda_f>0$, and
	let $0<\lambda<\lambda_f$.
	Then the following hold:
	\begin{enumerate}
\item\label{i:rock:prox} $e_{\lambda}[h_{\lambda} f]=e_{\lambda}f$;
\item\label{identity}
		$P_\lambda\left(h_\lambda f\right)
		=
		\conv(P_\lambda f)$;
\item\label{mm}
		$\conv(P_\lambda f)$
		is maximally monotone.
	\end{enumerate}
\end{lemma}
\proof
%\red{TODO: Add proof}
Recall that $h_{\lambda}f+\lambda^{-1}j$ is a proper lsc convex function by
Lemma~\ref{lem: proximal and closed cvx hull}.

\ref{i:rock:prox}: See \cite[Example 1.44]{rockafellar_variational_1998}.

%\ref{i:rock:prox}: See \cite[Example 1.44]{rockafellar_variational_1998}.
\ref{identity}: We provide a proof similar to \cite[Lemma 2.9]{cwp20} for completeness.
By \ref{i:rock:prox} and \cite[Example 10.32]{rockafellar_variational_1998},
\[
(\forall x\in\Rn)~
\conv \left[P_\lambda( h_\lambda f)(x)\right]
=
\lambda\partial[-e_\lambda(h_\lambda f)](x)+x
=
\lambda\partial(-e_\lambda f)(x)+x=\conv \left[P_\lambda f(x)\right].
\]
Because $h_\lambda f+\lambda^{-1}j$ is convex, $P_\lambda (h_\lambda f)$ is convex-valued as the set of minimizers of a convex function, furnishing
\(
P_\lambda(h_\lambda f)(x)=\conv \left[P_\lambda( h_\lambda f)(x)\right]=\conv \left[P_\lambda f(x)\right]
\) for every  $x\in\Rn$.

\ref{mm}:
Using \ref{identity}, the desired result follows immediately from that $P_\lambda(h_\lambda f)$ is maximally monotone;
see \cite[Proposition 12.19]{rockafellar_variational_1998}.
\begin{comment}
	\ref{identity}
By assumption and~\cite[Theorem 1.25]{rockafellar_variational_1998}, $P_\lambda f$ is compact-valued and upper semicontinuous. Thus
\begin{equation}\label{conv f}
	(\forall x\in\Rn)~
	\left(\convP_\lambda f\right)(x)
	=
	\overline{\conv}
	\left(P_\lambda f(x)\right)
	=
	\conv\left(P_\lambda f(x)\right).
\end{equation}
Appealing to~\cite[Theorem 5.4]{cwp20} with $\alpha=1$ yields
\(
P_\lambda \varphi_\mu^1
=
\convP_\lambda f.
\)
Note that
\[
\varphi_\mu^1
=
-e_\mu(-e_\mu f)
=
h_\lambda f.
\]
%where the last equality holds due to~\cite[Fact 2.1(b)]{cwp20}.
Then combining the above identities with~(\ref{conv f}) gives
\[
P_\lambda\left(h_\lambda\right)
=
\prox_\mu(\varphi_\mu^1)
=
\convP_\lambda f
=
\conv(P_\lambda f).
\]
\ref{mm}
Invoking~\cite[Fact 2.1(b)]{cwp20} implies that
\[
h_\lambda f
=
\left(
f+\frac{1}{2\mu}\norm{\cdot}^2
\right)^{**}
-
\frac{1}{2\mu}
\norm{\cdot}^2,
\]
thus \(h_\lambda f\) is \(\mu\)-hypoconvex.
In turn~\cite[Proposition 12.19]{rockafellar_variational_1998} and \ref{identity} implies that $\conv(P_\lambda f)$ is maximally monotone.
\end{comment}
\endproof

Below we say that $f$
is \emph{$\lambda$-level proximal subdifferentiable} if
$(\forall x\in\dom f)$ $\partial_p^\lambda f(x)\neq\varnothing$.
%A full characterization of level proximal subdifferentiable functions comes as follows.
\begin{theorem}[level proximal subdifferentiable functions]
\label{thm: full domain}
Let $f:\Rn\to\overline{\mathbb{R}}$ be proper, lsc and prox-bounded with threshold $\lambda_f>0$, and let $0<\lambda<\lambda_f$.  Suppose that $\dom f$ is convex. Then the following are equivalent:
\begin{enumerate}
	\item\label{thm:full domain:: full domain} $f$ is $\lambda$-level proximal
subdifferentiable;
%i.e.,
%	$(\forall x\in\dom f)$ $\partial_p^\lambda f(x)\neq\emptyset$;% and $\dom \left[\partial \conv (f+\lambda^{-1}j)\right]=\dom f$;
		\item\label{thm:full domain:: cvx hull} $P_\lambda f=\conv( P_\lambda f)$ and $\dom \left[\partial \conv (f+\lambda^{-1}j)\right]=\dom f$;
		% and $\dom \left[\partial \conv (f+\lambda^{-1}j)\right]=\dom f$;
	\item\label{thm:full domain:: mm} $P_\lambda f$ is maximally monotone and $\dom \left[\partial \conv (f+\lambda^{-1}j)\right]=\dom f$;
	% and $\dom \left[\partial \conv (f+\lambda^{-1}j)\right]=\dom f$;
		\item\label{thm:full domain:: hypocvx} $f+\lambda^{-1}j$ is convex and $\dom \partial f=\dom f$;% and $\dom \left[\partial \conv (f+\lambda^{-1}j)\right]=\dom f$;
	\item\label{thm:full domain:: coincide sub} $\partial_p^\lambda f=\partial f$ and $\dom \partial f=\dom f$.
\end{enumerate}
%When one of the above holds, $\dom \left[\partial \operatorname{conv}(f+\lambda^{-1}j)\right]=\dom f$
\end{theorem}
\proof
%\red{TO DO: Add proof.}
%Apply~\cite[Proposition 12.19]{rockafellar_variational_1998}, Theorems \ref{fact: relations} and~\ref{thm:domain}.
Let $g :=f+\lambda^{-1}j$.
Then $g$ is 1-coercive with $\dom g=\dom f$ and $\dom (\conv g)=\conv (\dom g)=\conv(\dom f)=\dom f$, owing to Lemma \ref{lem: proximal and closed cvx hull} and that $\dom f$ is convex.

``\ref{thm:full domain:: full domain}$\Rightarrow$\ref{thm:full domain:: cvx hull}": Appealing to Theorem~\ref{thm:domain} yields $(\forall x\in\dom f)$
$
\partial_p^\lambda f(x)\neq\emptyset
\Leftrightarrow
(\forall x\in\dom f)~
h_\lambda f(x)=f(x)\text{ and }x\in\dom \left[\partial(\operatorname{conv}g )\right],
$
thus
\begin{equation}\label{dom inlcusion}
	\dom f\subseteq \dom[\partial (\conv g)].
\end{equation}
Hence the convexity of $\dom f$ and (\ref{dom inlcusion}) enforce that
\(
\dom f\subseteq \dom\partial( \conv g)\subseteq\dom(\conv g)=\dom f.
\)
To show that $P_\lambda f=\conv\left(P_\lambda f\right)$, we claim that $h_\lambda f=f$, which implies through Lemma \ref{lem:prox of proximal hull}\ref{identity} that
\(
P_\lambda f=P_\lambda\left(h_\lambda f \right)=\conv( P_\lambda f)
\).
Indeed, $\dom h_\lambda f=\conv\left(\dom f \right)=\dom f$ due to the convexity of $\dom f$ and \cite[Example 11.26(c)]{rockafellar_variational_1998}, justifying our claim.

``\ref{thm:full domain:: cvx hull}$\Rightarrow$\ref{thm:full domain:: mm}": Apply Lemma~\ref{lem:prox of proximal hull}\ref{mm}.

``\ref{thm:full domain:: mm}$\Rightarrow$ \ref{thm:full domain:: hypocvx}":
Apply \cite[Example 11.26(c)]{rockafellar_variational_1998} to see that $f+\lambda^{-1}j$ is convex.
In turn, the assumption $\dom \left[\partial \operatorname{conv}(f+\lambda^{-1}j)\right]=\dom f$ yields
\(
\dom \partial f=\dom \partial(f+\lambda^{-1}j)=\dom \left[\partial \operatorname{conv}(f+\lambda^{-1}j)\right]=\dom f.
\)

``\ref{thm:full domain:: hypocvx}$\Rightarrow$ \ref{thm:full domain:: coincide sub}$\Rightarrow$
\ref{thm:full domain:: full domain}":  The first implication is a direct application of Fact~\ref{fact: relations}\ref{hypoconvex}, while the second is trivial.
\endproof

\begin{remark}
\begin{enumerate}
\item One always has $\dom\partial_{p}^{\lambda}f=\ran P_{\lambda}f$ by Fact~\ref{fact: prox identity}.
\item
Since $\lambda$-level proximal subdifferentiable functions are hypoconvex,
they are prox-regular and subdifferentially
regular.
\item A $\lambda$-level proximal subdifferentiable function might have unbounded subdifferentials at
a boundary point of its domain. It suffices to consider the convex function $$f:\R\rightarrow\bR:x\mapsto\begin{cases}
x, & \text{ if $x\geq 0$;}\\
+\infty, & \text{ otherwise.}
\end{cases}
$$
\end{enumerate}
\end{remark}
%\begin{remark} Although in
%Theorem~\ref{thm: full domain}\ref{thm:full domain:: hypocvx}\&\ref{thm:full domain:: coincide sub},
%the limiting subdifferential $\partial f$ was used,
%it can be replaced by any subdifferential in Fact~\ref{fact: relations}\ref{hypoconvex}.
%\end{remark}
The following example shows that
equivalences in Theorem~\ref{thm: full domain} collapse, if $\dom f$ is not convex.
\begin{example}
Consider the proper, lsc and prox-bounded function $f
:=\delta_{\{-1,1\}}$. Then $\dom f=\{-1,1\}$ is not convex and
\[
(\forall \lambda>0)~
\partial_p^\lambda f(x)=
	\begin{cases}
	[-1/\lambda,\infty),&\text{if }x=1;\\
	(-\infty,1/\lambda], &\text{if }x=-1;\\
	\emptyset,&\text{otherwise, }
\end{cases}
\text{ and }
P_\lambda f(x)
=
\begin{cases}
	\{-1,1\},&\text{if }x=0;\\
	\sgn(x),&\text{if }x\neq0.
\end{cases}
\]
Clearly $(\forall \lambda>0)$ $f$ is $\lambda$-level proximal subdifferentiable,
but $P_\lambda f\subset \conv (P_\lambda f)$.
\end{example}
\proof See \cite[Example 4.1]{wang23level}.
\endproof

%\begin{comment}
%Equipped with Theorem \ref{thm: full domain}, we conclude the following striking corollary.
%\begin{corollary}
%	Let $f:\Rn\to\overline{\mathbb{R}}$ be proper, lsc and prox-bounded with threshold $\lambda_f>0$, and let $0<\lambda<\lambda_f$.
%	Then $P_\lambda f$ is single-valued if and only if $P_\lambda f$ is convex-valued.
%\end{corollary}
%\proof
%Obviously $P_\lambda f$ is convex-valued provided that it is single-valued.
%The converse implication holds due to Theorem \ref{thm: full domain} and \cite{wang2010Chebyshev}.
%\endproof
%\end{comment}

\begin{comment}
	\begin{corollary}[connection to Lipschitz smoothness]\label{cor: lip smooth}
Let $f:\Rn\to\R$ be lsc and prox-bounded with threshold $\lambda_f>0$, and let $0<\lambda<\lambda_f$.
Then the following are equivalent:
	\begin{enumerate}
		\item $\dom \partial_p^\lambda f=\dom \partial_p^\lambda(-f)=\Rn$;
		\item $f$ is $1/\lambda$-Lipschitz smooth.
	\end{enumerate}
\end{corollary}
\proof
%\red{TODO: Add proof.}
Note that $f$ is $1/\lambda$-Lipschitz smooth if and only if $(1/\lambda)j\pm f$ are convex; see for instance \cite[Lemma 2.5]{wang2022mirror}. The equivalence then follows from Theorem~\ref{thm: full domain}.
\endproof
\end{comment}

\subsection{Locally single-valued level proximal subdifferential}
%\label{s:single:level}
Another interesting question is under what conditions $\partial_{p}^{\lambda}f(x)$ is a singleton.
In this subsection, we will show that the level proximal subdifferential of $f$ is single-valued
at a point
if and only if its proximal hull $h_{\lambda}f$ is differentiable and $f$ is $\lambda$-proximal
at the point, i.e., $h_{\lambda}f$ and $f$ agree at the point.
The following result is well-known, and we omit its simple proof.

\begin{lemma}\label{lem: Fenchel of conv}
Let $f:\Rn\to\overline{\mathbb{R}}$ be proper and lsc, and let $x\in\dom f$.
\begin{enumerate}
\item\label{i:conv:sub} If $(\conv f )(x)= f (x),$
then
$\partial_F f (x)=\partial(\conv f )(x);$
\item\label{i:clconv:sub}
 If $(\clconv f )(x)= f (x),$
then
$\partial_F f (x)=\partial(\clconv f )(x).$
\end{enumerate}
\end{lemma}

Now we present characterizations of single-valued $\partial_p^\lambda f(x)$
via the differentiability of the proximal hull of $f$ or the convex hull of
$f+\lambda^{-1}j$.
Note that
in finite-dimensional spaces for a finite-valued convex function G\^ateaux and Fr\'echet differentiabilties are the same
\cite[Corollary 17.44]{BC}.
%the main result of this section.

\begin{theorem}[level proximal subdifferential single-valued at a point]
\label{thm: equi of singleval}
Let $f:\Rn\to\overline{\mathbb{R}}$ be proper, lsc and prox-bounded with threshold $\lambda_f>0$.
Let $0<\lambda<\lambda_f$ and let $x\in\dom\partial_p^\lambda f$.
Then the following are equivalent:
\begin{enumerate}
\item\label{thm: equi of singleval:single} $\partial_p^\lambda f(x)$ is a singleton;
\item\label{thm: equi of singleval:conv diff} $\conv\left(f+\lambda^{-1}j\right)$ is differentiable at $x$, and
$\conv\left(f+\lambda^{-1}j\right)(x)=\left(f+\lambda^{-1}j\right)(x)$;
\item\label{thm: equi of singleval:prox hull}
$h_\lambda f$ is differentiable at $x$, and $h_\lambda f(x)=f(x)$.
% i.e., $f$ is $\lambda$-proximal at $x$;
%\item\label{thm: equi of singleval:prox}
% There exist a unique $v\in\mathbb{R}^n$, a finite collection of points $\{x_i\}_{i\in I}\subseteq\dom f$ called by $x$ (for function $f+\lambda^{-1}j$)
% % and $\{\alpha_j\}_{j\in J}\subseteq[0,1]$
% such that
% \begin{equation}\label{e:allinprox}
% (\forall i\in I)\ x_i\in P_\lambda f(v), \text{ and } x\in P_\lambda f(v).
% \end{equation}
%  In particular, % $x=\sum_{j\in J}\alpha_jx_j$, $\sum_{j\in J}\alpha_j=1$, and
%$v=\lambda\nabla \conv(f+\lambda^{-1}j)(x)$.
\end{enumerate}
\end{theorem}
\proof
``\ref{thm: equi of singleval:single}$\Rightarrow$\ref{thm: equi of singleval:conv diff}":
The Fenchel subdifferential $\partial_F\left(f+\lambda^{-1}j\right)(x)$ is a singleton by assumption and Fact~\ref{fact: level and Fenchel sub}.
Theorem~\ref{thm:domain} yields that
$
\conv\left(f+\lambda^{-1}j\right)(x)=\left(f+\lambda^{-1}j\right)(x).
$
In view of Lemma \ref{lem: Fenchel of conv}, we conclude that $\partial\conv\left(f+\lambda^{-1}j\right)(x)$ is a singleton,
so $\conv\left(f+\lambda^{-1}j\right)$ is differentiable at $x$.

``\ref{thm: equi of singleval:conv diff}$\Rightarrow$\ref{thm: equi of singleval:single}":
The convex function $\conv\left(f+\lambda^{-1}j\right)$ being differentiable at $x$ implies that $\partial\conv\left(f+\lambda^{-1}j\right)(x)$ is single-valued,
so is $\partial_{F}(f+\lambda^{-1}j)(x)$ by Lemma~\ref{lem: Fenchel of conv} and
the assumption $\conv\left(f+\lambda^{-1}j\right)(x)=\left(f+\lambda^{-1}j\right)(x)$.
Therefore applying Fact~\ref{fact: level and Fenchel sub} again completes the proof.

``\ref{thm: equi of singleval:conv diff}$\Leftrightarrow$\ref{thm: equi of singleval:prox hull}":
Invoke the identity $h_\lambda f=\conv (f+\lambda^{-1}j)-\lambda^{-1}j$;
see Lemma~\ref{lem: proximal and closed cvx hull}\ref{i:convexity} or
\cite[Example 11.26(c)]{rockafellar_variational_1998}.
\endproof

A tantalizing consequence of Theorem~\ref{thm: equi of singleval} comes as follows.

\begin{corollary}\label{c:levelsub-sing}
	Let $f:\Rn\to\overline{\mathbb{R}}$ be proper, lsc and prox-bounded with threshold $\lambda_f>0$.
	Let $0<\lambda<\lambda _f$ and let $U\subseteq\Rn$ be a nonempty open convex set.
	Then $\partial_p^\lambda f$ is single-valued on $U$ if and only if $f$ is $\lambda$-proximal and differentiable on $U$.
Consequently,
$\partial_p^\lambda f$ is single-valued on $\Rn$ if and only if $f+\lambda^{-1}j$ is
convex and $f$ is differentiable on $\Rn$.
\end{corollary}

%A tantalizing consequence of Corollary~\ref{c:levelsub-sing} comes as follows.
%\begin{corollary}
%	Let $f:\Rn\to\overline{\mathbb{R}}$ be proper, lsc and prox-bounded with threshold $\lambda_f=+\infty$.
%If $\partial_p f$ is single-valued on $\Rn$, then $f$ is
%convex and differentiable on $\Rn$.
%\end{corollary}
%\begin{proof} By Fact~\ref{f:levelsub:f}\ref{i:f2}, for each $\lambda>0$, $P_{\lambda}f$ is single-valued.
%Corollary~\ref{c:levelsub-sing} shows that $f+\lambda^{-1}j$ is
%convex and differentiable on $\Rn$. The result follows by letting $\lambda\rightarrow\infty$.
%\end{proof}

\section{Variational convexity and variational strong convexity}\label{s:variat}
In this section, we establish equivalences among local relative monotonicity of level proximal
subdifferential, local nonexpaniveness of proximal mapping,
 and variational convexity of the function. Similar results are also established for variationally
 strong convex functions.
Our results extend those
by Khanh, Mordukhovich and Phat \cite[Theorem 3.2]{Khanh2023variational} and by Rockafellar
\cite[Theorem 2]{rock2021}. Our new characterizations contain
local firm nonexpansiveness or nonexpansiveness
of proximal mappings,
and they are crucial for algorithms in variationally convex optimization; see section~\ref{s:algorithms}.

Let us first recall Rockafellar's variational convexity.
%Recall that a lsc function $f:\Rn\rightarrow\bR$
%is called \emph{variationally convex} at $\bar x$ for $\bar v\in\partial f(\bar x)$
%if for some convex neighborhood $U\times V$ of
%$(\bar x, \bar v)$ there exist a lsc convex function $\varphi\leq f$ on $U$ and a number
%$\varepsilon>0$ such that
%$[U_{\varepsilon}\times V]\cap \gph\partial f=[U\times V]\cap \gph \partial \varphi
%\text { and } f(x)=\varphi(x) \text{ at the common element $(x,v)$,} $
%where $U_{\varepsilon}=\{x\in U|\ f(x)< f(\bar x)+\varepsilon\}$;
%see \cite[Definition 2]{rock-vietnam}.
\begin{defn}\label{def:vc}\emph{\cite[Definition 2]{rock-vietnam}}
A lsc function $f:\Rn\rightarrow\overline{\R}$ is called variationally convex at $\bar x$ for $\bar v\in\partial f(\bar x)$,
if for some convex neighborhood $U\times V$ of $(\bar x, \bar v)$ there exist a lsc convex function $\varphi\leq f$ on $U$ and a number $\varepsilon>0$ such that
$[U_{\varepsilon}\times V]\cap \gph\partial f=[U\times V]\cap \gph \partial \varphi\text { and } f(x)=\varphi(x) \text{ at the common element $(x,v)$,} $
where $U_{\varepsilon}:=\{x\in U|\ f(x)< f(\bar x)+\varepsilon\}$.
\end{defn}
Of course, if a function is convex on a neighborhood of a point, then it is variationally convex at the point.
In order for $f$ to be variationally convex at $\bar x$ for $\bar v\in\partial f(\bar x)$, it is necessary that
$\bar v\in\hat{\partial}f(\bar x)$; e.g., $x\mapsto -\|x\|$ is not variationally convex at $0$ for
every $\bar v\in\partial (-\|\cdot\|)(0)$ since $\hat{\partial}(-\|\cdot\|)(0)=\varnothing$.
We give two nonconvex examples to illustrate variational convexity.
\begin{example} Define the counting (or $\ell_{0}$) norm $\|\cdot\|_{0}:\Rn\rightarrow\R$ by
$x\mapsto \sum_{i=1}^{n}|x_{i}|_{0}$, where $|t|_{0}:=1$ if $t\neq 0$ and $0$ otherwise.
Then $\|\cdot\|_{0}$ is variationally convex everywhere on
$\Rn$ but it is not convex\footnote{In \cite[Example 2.5]{Khanh2023variational}, the authors
showed that $\|\cdot\|_{0}$ is variationally
convex at $0$. Analogous arguments can show that $\|\cdot\|_{0}$ is variationally
convex everywhere.}.
This is in stark contrast with the classical result: A finite-valued and locally convex function on an open convex set is a convex function; see \cite[Lemma 2.1]{penot88}.
\end{example}

\begin{example}
	Consider the function \(g:\Rn\to\R\) defined by
	\[x\mapsto \card\{1\leq i\leq n-1|\ x_i\neq x_{i+1}\},
	\]
where $\card(A)$ denotes the cardinality of set $A$.
	Then \(g\) is variationally convex at every \(\bar x\in\Rn\)
for every \(\bar v\in\partial g(\bar x)\).
\end{example}
\begin{proof}
	Define the difference operator \(D:\Rn\to\R^{n-1}\) by
	\[
	D=
	\begin{bmatrix}
		1 & -1 & 0 & \cdots & 0\\
		0 & 1 & -1 & \cdots & 0\\
		\vdots & \vdots & \vdots & \ddots & \vdots\\
		0 & 0 & 0 & \cdots & -1
	\end{bmatrix}.
	\]
	Then \(g(x)=\norm{Dx}_0\) for every \(x\in\Rn\) and \(D\) is surjective (so $D^*$ injective).
	Invoking the subdifferential chain rule \cite[Theorem 10.6]{rockafellar_variational_1998} gives
	\(
	(\forall x\in\Rn)~
	\partial g(x)
	=
	D^*\partial\norm{\cdot}_0(Dx).
	\)
	Because \(\norm{\cdot}_0\) is variationally convex at every point for every subgradient, it is so at \(D\bar x\) for \(\bar u\in\partial\norm{\cdot}_0(D\bar x)\) such that \(\bar v=D^*\bar u\).
	Suppose that there exists \(\varepsilon>0\) such that
	\begin{equation}\label{eq:jump function}
		\ip{u_1-u_2}{y_1-y_2}\geq0,
	\end{equation}
	whenever \(\norm{y_i-D\bar x}<\varepsilon\) with \(\norm{y_i}_0<\norm{D\bar x}_0+\varepsilon\) and \(u_i\in\partial\norm{\cdot}_0(y_i)\) with \(\norm{u_i-\bar u}<\varepsilon\).
	Now pick \((x_i,v_i)\in\gra\partial g\) such that \(\norm{x_i-\bar x}<\norm{D}^{-1}\varepsilon\) and \(\norm{v_i-\bar v}<\sigma(D^*)\varepsilon\), where \(\sigma(D^*)>0\) is the minimal singular value of \(D^*\), and \(g(x_i)<g(\bar x)+\varepsilon\). Here $\sigma(D^*)>0$ because $D^*$ is injective.
	Then to each \(v_i\) corresponds \(u_i\in\partial\norm{\cdot}_0(Dx_i)\) such that \(v_i=D^*u_i\) and
	\(
		\norm{u_i-\bar u}
	\leq
		\norm{D^{*}}^{-1}\norm{v_i-\bar v}<\varepsilon
	\).
	Moreover, it holds that
	\(
		\norm{Dx_i-D\bar x}
	\leq
		\norm{D}\norm{x_i-\bar x}
	<
		\varepsilon
	\).
	Hence, the pair \((Dx_i,u_i)\in\gra\partial\norm{\cdot}_0\) satisfies all condition for \eqref{eq:jump function}, and it follows that
	\[
		\ip{v_1-v_2}{x_1-x_2}
	=
		\ip{D^*u_1-D^*u_2}{x_1-x_2}
	=
		\ip{u_1-u_2}{Dx_1-Dx_2}
	\geq 0.
	\]
	Note that \(\partial g(\bar x)=D^*\partial\norm{\cdot}_0(D\bar x)=D^*\hat\partial\norm{\cdot}_0(D\bar x)=\hat\partial g(\bar x)\).
	Then \(\bar v\in\hat\partial g(\bar x)\), which means that all conditions in \cite[Theorem 1]{Rockafellar2024} are satisfied, completing the proof.
\end{proof}

Coming back to our characterization of variational convexity, we need three auxiliary results.
The first one is a localized Minty's parametrization.
\begin{lemma}\label{l:minty} Let $U$ be a nonempty set in $\Rn$, let $\lambda>0$, and let
$W_{\lambda} :=\{(u,v)\in\Rn\times\Rn| \ u+\lambda v\in U\}.$
\begin{enumerate}
\item \label{i:operator}
Suppose that $A:\Rn\rightrightarrows\Rn$ and that $J_{\lambda A}$ is single-valued on $U$.
Then
\begin{equation}\label{e:mintylocal}
\gph A \cap W_{\lambda}=\{(J_{\lambda A}(x), \lambda^{-1}[x-J_{\lambda A}(x)])|\ x\in U\}.
\end{equation}
\item\label{i:function}
 Suppose that $f:\Rn\rightarrow\bR$ is proper, lsc and prox-bounded with threshold $\lambda_{f}$, and $0<\lambda<\lambda_{f}$. If $P_{\lambda}f$ is single-valued on $U$, then
 \begin{align}
 \gph \partial_{p}^{\lambda}f\cap W_{\lambda} &=\{(P_{\lambda}f(x), \lambda^{-1}[x-P_{\lambda}f(x)])|\ x\in U\}\label{e:proxf}\\
 & =\{(x-\lambda \nabla e_{\lambda}f(x),
    \nabla e_{\lambda}f(x))|\ x\in U\}.\label{e:gradf}
 \end{align}
\end{enumerate}
\end{lemma}
\proof
\ref{i:operator}:
Denote the right-hand side of \eqref{e:mintylocal} by $E_{\lambda}$. We have
\begin{align*}
(u,v)\in E_{\lambda} & \Leftrightarrow (\exists x\in U)\ u=J_{\lambda A}(x), v=\lambda^{-1}[x-J_{\lambda A}(x)]\\
& \Leftrightarrow (\exists x\in U)\ x\in u+\lambda Au, x=u+\lambda v\\
& \Leftrightarrow (\exists x\in U)\ v\in Au, x=x+\lambda v
\Leftrightarrow (u,v)\in \gph A \cap W_{\lambda}.
\end{align*}
\ref{i:function}:  Apply \ref{i:operator} with $A=\partial_{p}^{\lambda}f$ to obtain
\eqref{e:proxf}.
When $0<\lambda<\lambda_{f}$ and $P_{\lambda}f$ is single-valued on $U$,
\cite[Example 10.32]{rockafellar_variational_1998} implies
$(\forall x\in U)\ \nabla e_{\lambda}f(x)=\lambda^{-1}[x-P_{\lambda}f(x)]$.
Hence \eqref{e:gradf} holds,
\endproof

The second one connects local cocoercivity of resolvents and
relative monotonicity through Minty parametrization.
\begin{lemma}\label{lem:local coco Minty}
	Let $A:\Rn\rightrightarrows\Rn$ be a set-valued mapping, let $\lambda>0$, and
	let $\varnothing\neq U\subseteq\dom J_{\lambda A}$.
	Define $W_\lambda \coloneqq \{(u,v)\in\Rn\times\Rn|\ u+\lambda v\in U\}$.
	Then for \(\gamma>0\) the following are equivalent:
\begin{enumerate}
\item\label{i:resolvent:firm}
	$J_{\lambda A}$ is \(\gamma\)-cocoercive on $U$;
\item\label{i:monotone:hypo}
$A$ is \(\tfrac{1}{\lambda}\left(\gamma-1\right)\)-monotone relative to $W_\lambda$.
\end{enumerate}
%	\red{\textbf{(This should be comono because \(\gamma-1\) may be negative.)}}
%	\begin{enumerate}
%		\item\label{thm:local Minty::coco}
%		$J_{\lambda A}$ is \(\gamma\)-cocoercive on $U$;
%		\item\label{thm:local Minty::str mono}
%		$A$ is \(\tfrac{1}{\lambda}\left(\gamma-1\right)\)-comonotone relative to $W_\lambda$.
%	\end{enumerate}
\end{lemma}
\proof
``\ref{i:resolvent:firm}$\Rightarrow$\ref{i:monotone:hypo}'':
	Take $(u_{i},v_{i})\in \gph A
	\cap W_{\lambda}$ and set $x_{i} \coloneqq u_{i}+\lambda v_{i}$ for $i=1, 2$. Then
	$x_{i}\in U$ and
	$(u_{i},v_{i})=(J_{\lambda A}(x_{i}), \lambda^{-1}[x_{i}-J_{\lambda A}(x_{i})])$
	by Lemma~\ref{l:minty}\ref{i:operator}.
	Then the assumption implies that
	\begin{align*}
		\ip{u_1-u_2}{v_1-v_2}
		&=
		\tfrac{1}{\lambda}\ip{J_{\lambda A}(x_1)-J_{\lambda A}(x_2)}{(x_1-x_2)-(J_{\lambda A}(x_1)-J_{\lambda A}(x_2))}\\
		&\geq
		\tfrac{1}{\lambda}\left(\gamma-1\right)\norm{J_{\lambda A}(x_1)-J_{\lambda A}(x_2)}^2
		=
		\tfrac{1}{\lambda}\left(\gamma-1\right)\norm{u_1-u_2}^2,
	\end{align*}
	as claimed.

``\ref{i:monotone:hypo}$\Rightarrow$\ref{i:resolvent:firm}'':
	Suppose that $x\in U\subseteq \dom J_{\lambda A}$.
	Then $x=u+\lambda v$ for some $(u,v)\in\gph A$, consequently $(u,\lambda^{-1}(x-u))\in\gph A\cap W_\lambda$ and $u\in J_{\lambda A}(x)$.
	We claim that $J_{\lambda A}$ is single-valued.
	Indeed, let $u'\in J_{\lambda A}(x)$.
	Then the same process as above gives $(u',\lambda^{-1}(x-u'))\in\gph A\cap W_\lambda$, and the
\(\tfrac{1}{\lambda}\left(\gamma-1\right)\)-monotonicity of
$\gph A\cap W_\lambda$ implies
	\[
	\ip{u-u'}{\tfrac{1}{\lambda}(x-u)-\tfrac{1}{\lambda}(x-u')}
	\geq
	\tfrac{1}{\lambda}(\gamma-1)\norm{u-u'}^2
	\Leftrightarrow
	\gamma\norm{u-u'}^2\leq 0 \Leftrightarrow u=u'.
	\]
	Next, let $u_i \coloneqq J_{\lambda A}(x_i)$ with $x_i\in U$ for $i=1,2$.
	Then there exists \(v_i\in A(u_i)\) such that \(x_i=u_i+\lambda v_i\) and
	the \(\tfrac{1}{\lambda}\left(\gamma-1\right)\)-monotonicity of $\gph A\cap W_\lambda$ gives
	$
	\ip{u_1-u_2}{v_1-v_2}\geq\tfrac{1}{\lambda}\left(\gamma-1\right)\norm{u_1-u_2}^2
	$, i.e.,
$\ip{u_1-u_2}{(\lambda v_1+u_{1})-(\lambda v_2+u_{2})}\geq \gamma\norm{u_1-u_2}^2$. Thus,
$\ip{J_{\lambda A}(x_{1})-J_{\lambda A}(x_{2})}{x_{1}-x_{2}}\geq \gamma\norm{J_{\lambda A}(x_{1})-J_{\lambda A}(x_{2})}^2$
 and this is precisely \ref{i:resolvent:firm}.
\endproof

The third one relates local firm nonexpansiveness of resolvents to relative maximal monotonicity; and it
extends \cite[Proposition 1]{rock2021}.

\begin{lemma}\label{thm:local Minty}
Let $A:\Rn\rightrightarrows\Rn$ be a set-valued mapping, let $\lambda>0$, and
let $\varnothing\neq U\subseteq\dom J_{\lambda A}$.
Define $W_\lambda :=\{(u,v)\in\Rn\times\Rn|\ u+\lambda v\in U\}$.
Then the following are equivalent:
\begin{enumerate}
	\item\label{thm:local Minty::fnexp}
	 $J_{\lambda A}$ is firmly nonexpansive on $U$;
	\item\label{thm:local Minty::mono}
	$A$ is monotone relative to $W_\lambda$;
	\item\label{thm:local Minty::max mono}
	 $A$ is maximally monotone relative to $W_\lambda$.
\end{enumerate}
\end{lemma}
\proof
Observe that
$
	[(x,u)\in\gph J_{\lambda A} \text{ and } x \in U]
	\Leftrightarrow
	\left(u,\lambda^{-1}(x-u)\right)\in\gph A\cap W_\lambda.
$
%Indeed
%$
%	(z,x)\in\gph J_{\lambda A}, z\in U
%	\Leftrightarrow
%	x\in J_{\lambda A}(z),z\in U
%	\Leftrightarrow
%	(z-x)/\lambda\in Ax,z\in U
%$.

``\ref{thm:local Minty::fnexp}$\Leftrightarrow$
\ref{thm:local Minty::mono}'': Apply Lemma~\ref{lem:local coco Minty} with $\gamma=1$.
 %Take $(u_{i},v_{i})\in \gph A
%\cap W_{\lambda}$ and set $x_{i} :=u_{i}+\lambda v_{i}$ for $i=1, 2$. Then
%$x_{i}\in U$ and
%$(u_{i},v_{i})=(J_{\lambda A}(x_{i}), \lambda^{-1}[x_{i}-J_{\lambda A}(x_{i})])$
%by Lemma~\ref{l:minty}\ref{i:operator}.
%Using \ref{thm:local Minty::fnexp}, we derive
%\begin{align*}
%\ip{u_{1}-u_{2}}{v_{1}-v_{2}} &= \ip{J_{\lambda A}(x_{1})-J_{\lambda A}(x_{2})}{\lambda^{-1}[x_{1}-J_{\lambda A}(x_{1})]-\lambda^{-1}[x_{2}-J_{\lambda A}(x_{2})]}\\
%&=\lambda^{-1}\bigg(\ip{J_{\lambda A}(x_{1})-J_{\lambda A}(x_{2})}{x_{1}-x_{2}}-\norm{J_{\lambda A}(x_{1})-J_{\lambda A}(x_{2})}^2\bigg)\geq 0.
%\end{align*}
%%Denote by $x_i=J_{\lambda A}(z_i)$ for $z_i\in U$ and $i=1,2$.
%%We have
%%\[
%%\ip{x_1-x_2}{z_1-z_2}\geq\|x_1-x_2\|^2.
%%\]
%%Dividing by $\lambda>0$ gives
%%\[
%%\ip{x_1-x_2}
%%{
%%	\frac{z_1-x_1}{\lambda}-\frac{z_2-x_2}{\lambda}
%%}
%%\geq0.
%%\]
%Hence $A$ is monotone relative to $W_\lambda$.

``\ref{thm:local Minty::mono}$\Rightarrow$\ref{thm:local Minty::max mono}":
Take $(u,v)\in W_\lambda$ that is monotonically related to $\gph A\cap W_\lambda$.
We claim that $(u,v)\in\gph A$, which entails the desired maximal monotonicity.
Indeed, $(u,v)\in W_\lambda$ implies $u+\lambda v\in U\subseteq \dom J_{\lambda A}$.
Let $y\in J_{\lambda A}(u+\lambda v)$.
Then
$
\left(y,\lambda^{-1}(u+\lambda v-y)\right)\in\gph A\cap W_\lambda.
$
By the assumption that $(u,v)$ is monotonically related to $\gph A\cap W_\lambda$, we obtain
\begin{align*}
\ip{u-y}{\lambda^{-1}(y-u)}
=
\ip{u-y}{v-\lambda^{-1}(u+\lambda v-y)}\geq 0
\Rightarrow
u=y,
\end{align*}
enforcing that $u\in J_{\lambda A}(u+\lambda v)$.
In turn, $u+\lambda v\in u+\lambda Au \Leftrightarrow v\in Au$.

``\ref{thm:local Minty::max mono}$\Rightarrow$\ref{thm:local Minty::mono}":
Clear.
%Suppose that $x\in U\subseteq \dom J_{\lambda A}$.
%Then $x=u+\lambda v$ for some $(u,v)\in\gph A$, consequently $(u,\lambda^{-1}(x-u))\in\gph A\cap W_\lambda$ and $u\in J_{\lambda A}(x)$.
%We claim that $J_{\lambda A}$ is single-valued.
%Indeed, let $u'\in J_{\lambda A}(x)$.
%Then the same process as above implies $(u',\lambda^{-1}(x-u'))\in\gph A\cap W_\lambda$.
%Appealing to the monotonicity of $\gph A\cap W_\lambda$ yields
%\[
%\ip
%{
%u-u'
%}
%{
%\lambda^{-1}(x-u)-\lambda^{-1}(x-u')
%}
%\geq 0
%\Leftrightarrow
%\ip{u-u'}{-\lambda^{-1}(u-u')}\geq 0 \Leftrightarrow u=u'.
%\]
%Next, let $u_i :=J_{\lambda A}(x_i)$ with $x_i\in U$ for $i=1,2$.
%Then the monotonicity of $\gph A\cap W_\lambda$ again gives
%$$
%\ip{u_1-u_2}{\lambda^{-1}(x_1-u_1)-\lambda^{-1}(x_2-u_2)}\geq0
%%\Leftrightarrow\ip{x_1-x_2}{(z_1-z_2)-(x_1-x_2)}\geq0
%\Leftrightarrow
%\ip{u_1-u_2}{x_1-x_2}\geq\|u_1-u_2\|^2,
%$$
%justifying \ref{thm:local Minty::fnexp}.
\endproof

Now we are ready to derive the first main result of this section showing the remarkable
interconnections among variational convexity
of the function, local relative monotonicity of level proximal subdifferential,
and local nonexpaniveness of proximal mapping.

\begin{theorem}[characterizations of variational convexity]\label{thm:coco env grad}
Suppose that $f:\mathbb{R}^n\to\overline{\mathbb{R}}$ is proper, lsc and prox-bounded
with threshold $\lambda_{f}>0$,
and that $f$ is prox-regular at $\bar x$ for $\bar v\in\partial f(\bar x)$.
%{\color{gray}
%Then for all $\lambda>0$ sufficiently small there exists an open convex neighborhood $U_{\lambda}$
%of $\bar{x}+\lambda \bar{v}$ on which the following are equivalent:
%}%
Then there exists \(0<\bar\lambda\leq \lambda_{f}\) such that $(\forall 0<\lambda<\bar{\lambda})\ P_{\lambda}f(\barx+\lambda\barv)=\barx$, and $\exists U_{\lambda}$
an open convex neighborhood of \(\bar x+\lambda\bar v\) for which the following equivalent properties hold:
	\begin{enumerate}
		\item\label{thm:coco env grad:: v cvx}
		$f$ is variationally convex at $\bar x$ for $\bar v$;
		\item\label{thm:coco env grad:: mono}
		\((\forall 0<\lambda<\bar\lambda)\)
		$\partial_p^\lambda f$ is monotone relative to $W_\lambda$;
% and $\bar v\in\partial_{p}^{\lambda}f(\bar x)$;
		\item\label{thmn:coco env grad:: max mono}
		\((\forall 0<\lambda<\bar\lambda)\)
		$\partial_p^\lambda f$ is maximally monotone relative to $W_\lambda$;
% and $\bar v\in\partial_{p}^{\lambda}f(\bar x)$;
		\item\label{thm:coco env grad:: firmly nonexap prox}
		\((\forall 0<\lambda<\bar\lambda)\)
		$P_\lambda f$ is firmly nonexpansive on $U_\lambda$;
		\item\label{thm:coco env grad:: averaged}
		\((\forall 0<\lambda<\bar\lambda)\) \((\exists0<\alpha<1)\)
		$P_\lambda f$ is $\alpha$-averaged on $U_\lambda$;
% and $P_{\lambda}f(\barx+\lambda\barv)=\barx$;
		\item\label{thm:coco env grad:: ne}
		\((\forall 0<\lambda<\bar\lambda)\)
		$P_\lambda f$ is nonexpansive on $U_\lambda$;
% and $P_{\lambda}f(\barx+\lambda\barv)=\barx$;
	%	\item\label{thm:coco env grad:: comono}
	%	\((\forall 0<\lambda<\bar\lambda)\) \((\exists0<\alpha<1)\)
	%	$\partial_p^\lambda f$ is $\rho$-monotone relative to \(W_\lambda\) with \(\rho=((2\alpha)^{-1}-1)\lambda^{-1}\);
		\item\label{envel:convex}
		\((\forall 0<\lambda<\bar\lambda)\)
		$e_\lambda f$ is convex on $U_\lambda$;
		\item\label{envel:convex spec}
		\((\exists0<\lambda<\bar \lambda)\)
		\(e_\lambda f\) is convex on $U_\lambda$;
		%($\partial_F\left(\lambda^{-1}j+f\right)$ is $1/\lambda$-strongly monotone
		%relative to $\Rn\times\left[ \lambda^{-1}U_\lambda\right]$.)
		\item\label{thm:coco env grad:: coco}
		\((\forall 0<\lambda<\bar\lambda)\)
		$\nabla e_\lambda f$ is $\lambda$-cocoercive on $U_\lambda$;
	\end{enumerate}
	\vspace*{-2mm}
		where $W_\lambda :=\{(x,v)\in\Rn\times\Rn|\ x+\lambda v\in U_\lambda\}$.
\end{theorem}

\proof
By
\cite[Proposition 5.3]{thibault},
there exists \(0<\bar\lambda\leq \lambda_{f}\) such that to every \(0<\lambda<\bar\lambda\) corresponds an open convex neighborhood $U_\lambda$ of $\bar x+\lambda\bar v$ on which:
\begin{equation}\label{e:grad1}
\text{$P_{\lambda}f$ is monotone, single-valued and Lipschitz continuous; $P_{\lambda}f(\bar x+\lambda{\bar v})=\bar x$;}
\end{equation}
and $e_\lambda f$ is differentiable with
\begin{equation}\label{e:grad2}
\nabla e_{\lambda}f=\lambda^{-1}[\Id-P_\lambda f].
\end{equation}
See also \cite[Proposition 13.37]{rockafellar_variational_1998} (for $\bar v=0$).
In particular, \eqref{e:grad1} implies $\bar{v}\in\partial_{p}^{\lambda}f(\bar{x})\subseteq\hat{\partial}f(\bar x)$;
see Fact \ref{fact: prox identity}.
%on (see \cite[Proposition 2.2]{Khanh2023variational}).

\begin{comment}
	``\ref{thm:coco env grad:: v cvx}$\Leftrightarrow$\ref{envel:convex}
$\Leftrightarrow$\ref{thm:coco env grad:: coco}":
Although only
$\ref{thm:coco env grad:: v cvx}\Leftrightarrow\ref{envel:convex}$ stated
in \cite[Theorem 3.2]{Khanh2023variational},
the proof there
 actually gives $\ref{thm:coco env grad:: v cvx}
\Rightarrow \ref{thm:coco env grad:: coco}\Rightarrow\ref{envel:convex}
\Rightarrow \ref{thm:coco env grad:: v cvx}$.
\end{comment}

``\ref{thm:coco env grad:: v cvx}$\Rightarrow$\ref{thm:coco env grad:: mono}'':
Since $f$ is variationally convex at $\bar x$ for $\bar v\in\hat{\partial}f(\bar{x})$,
there exist a convex neighborhood \(U\times V\) of \((\bar x,\bar v)\) and \(\varepsilon>0\) such that
\begin{equation}\label{e:subf:monotone}
	\gph\partial f\cap\left[U_\varepsilon\times V\right]
	\text{ is monotone},
\end{equation}
where \(U_\varepsilon :=\{x\in U|\ f(x)<f(\bar x)+\varepsilon\}\); see \cite[Theorem 1]{rock-vietnam}.
For every \(0<\lambda<\bar\lambda\), \eqref{e:grad1} and \eqref{e:grad2} entail that
$(\forall u\in U_{\lambda})\ u-\lambda\nabla e_\lambda f(u)=P_{\lambda}f(u)$,
$P_{\lambda}f$ is Lipschitz, and
\(\nabla e_\lambda f(\bar x+\lambda \bar v)=\bar v\).
Hence,
\begin{equation}\label{eq:limits}
\lim_{u\to\bar x+\lambda\bar v}\nabla e_\lambda f(u)=\bar v,
\lim_{u\to\bar x+\lambda\bar v}(\Id-\lambda\nabla e_\lambda f)(u)=\bar x,
\text{ and }\lim_{u\to\bar x+\lambda\bar v}f(P_{\lambda}f(u))=f(\bar x),
\end{equation}
where the last limit holds because
\begin{align}
\lim_{u\to\bar x+\lambda\bar v}	f(P_{\lambda}f(u)) & =\lim_{u\to\bar x+\lambda\bar v }\Big(e_\lambda f(u)-\tfrac{1}{2\lambda}\norm{u-P_{\lambda}f(u)}^2\Big)\\
&=e_\lambda f(\bar x+\lambda\bar v)-\tfrac{1}{2\lambda}\norm{\bar x+\lambda\bar v-\bar x}^2=f(\bar x).
\end{align}
Owing to \eqref{eq:limits}, by using a smaller $U_{\lambda}$ if necessary, we can assume
that for \(u\in U_\lambda\)
\begin{equation}\label{e:f:attentative}
	u-\lambda\nabla e_\lambda f(u)\in U, \nabla e_\lambda f(u)\in V,
\text{ and } f(u-\lambda \nabla e_\lambda f(u))<f(\bar x)+\varepsilon.
\end{equation}
Moreover,
\(u-\lambda\nabla e_\lambda f(u)=P_\lambda f(u)\) implies, by Fact~\ref{fact: prox identity}, that
\begin{equation}\label{e:in:graph}
(u-\lambda\nabla e_\lambda f(u),\nabla e_\lambda f(u))\in
\gra \partial_{p}^{\lambda}f\subseteq \gra\partial f.
\end{equation}
%\(\nabla e_\lambda f(u)\in\partial f(u-\lambda\nabla e_\lambda f(u))\).
Combining \eqref{e:f:attentative}, \eqref{e:in:graph} and \eqref{e:subf:monotone} we obtain
\[
	\{(u-\lambda\nabla e_\lambda f(u), \nabla e_\lambda f(u))|\ u\in U_\lambda\}\subseteq\gph\partial f\cap\left[U_\varepsilon\times V\right]
	\text{is monotone},
\]
and so is \(\gph\partial_p^\lambda f\cap W_\lambda\) by Lemma~\ref{l:minty}.

``\ref{thm:coco env grad:: mono}$\Leftrightarrow$\ref{thmn:coco env grad:: max mono}$\Leftrightarrow$\ref{thm:coco env grad:: firmly nonexap prox}'':
Note that \(P_\lambda f=(\Id+\lambda\partial_p^\lambda f)^{-1}\).
Then applying Lemma \ref{thm:local Minty} completes the proof.

``\ref{thm:coco env grad:: firmly nonexap prox}$\Rightarrow$\ref{thm:coco env grad:: averaged}$\Rightarrow$\ref{thm:coco env grad:: ne}": Clear.

``\ref{thm:coco env grad:: ne}$\Rightarrow$\ref{envel:convex}": The assumption implies that $\nabla e_{\lambda}f$ is monotone on $U_{\lambda}$, so $e_{\lambda}f$ is convex on $U_{\lambda}$
{in view of \cite[B Theorem 4.1.4]{hiriart2004fundamentals}}.
Indeed, by the nonexpansiveness of $P_\lambda f$ on $U_\lambda$ and the Cauchy-Schwarz inequality,
\begin{equation}\label{e:nonexp:grad}
(\forall x_{1}, x_{2}\in U_\lambda)~
\ip{P_\lambda f(x_{1})-P_\lambda f(x_{2})}{x_{1}-x_{2}}\leq\norm{x_{1}-x_{2}}^2.
\end{equation}
Since
$\nabla e_{\lambda}f=\lambda^{-1}[\Id-P_{\lambda}f] \text{ on $U_{\lambda}$},$
it follows from \eqref{e:nonexp:grad} that
$$(\forall x_{1},x_{2}\in U_\lambda)~\ip{x_{1}-x_{2}}{\nabla e_\lambda f(x_{1})-\nabla e_\lambda f(x_{2})}\geq 0,$$
hence $\nabla e_{\lambda}f$ is monotone on $U_{\lambda}$.

``\ref{envel:convex}$\Rightarrow$\ref{envel:convex spec}$\Rightarrow$\ref{thm:coco env grad:: v cvx}'':
The first implication is trivial.
We prove the second following the same line of proof in \cite[Theorem 3.2]{Khanh2023variational} for the sake of completeness.
By the prox-regularity of $f$ at $\bar{x}$ for $\bar{v}\in\partial f(\bar{x})$,
there exist \(\gamma>0\) and \(T:\Rn\rightrightarrows\Rn\) defined by
\begin{align*}
	T(x):=
	\begin{cases}
		\{v\in\partial f(x)|\norm{v-\bar v}<\gamma\},&\text{if }\norm{x-\bar x}<\gamma, f(x)<f(\bar x)+\gamma;\\
		\varnothing, &\text{otherwise,}
	\end{cases}
\end{align*}
such that \(\nabla e_\lambda f=\lambda^{-1}[\Id-P_{\lambda}f]=
\left(\lambda\Id+T^{-1}\right)^{-1}\) on an open convex neighborhood \(U_\lambda\) of $\bar{x}+\lambda\bar{v}$; see
\cite[Proposition 13.37]{rockafellar_variational_1998} or \cite[Proposition 4.4]{thibault}.
By taking a smaller $U_{\lambda}$ if necessary, assumption \ref{envel:convex spec} implies that
\begin{equation}\label{eq:cvx env}
	(\forall x,y\in U_\lambda)~
	e_\lambda f(y)\geq e_\lambda f(x)+\ip{\nabla e_\lambda f(x)}{y-x}.
\end{equation}
Take \(\varepsilon>0\) sufficiently small so that \(\varepsilon<\gamma\) and
\begin{equation}\label{e:uv}
(\forall x\in U)(\forall v\in V)\ x+\lambda v\in U_\lambda,
\end{equation}
where \(U:=\Ball_{\varepsilon}(\bar x)\) and \(V:=\Ball_{\varepsilon}(\bar v)\).
Then \(v\in T(x)\) for \((x,v)\in\gph\partial f\cap[U_\varepsilon\times V]\) with
\(U_\varepsilon:=\{x\in U|\ f(x)<f(\bar x)+\varepsilon\}\), implying that
\begin{align}
	\nabla e_\lambda f(x+\lambda v) &=\left(\lambda\Id+T^{-1}\right)^{-1}(x+\lambda v)=v, \text{ and }
	\label{eq:env and prox}\\
	P_\lambda f(x+\lambda v) &=(x+\lambda v)-\lambda\nabla e_{\lambda}f(x+\lambda v)=x.
\label{eq:env and prox1}
\end{align}
Therefore, for every \(y\in U\) and \((x,v)\in\gph\partial f\cap[U_\varepsilon\times V]\), we have \(y+\lambda v, x+\lambda v\in U_\lambda\) by \eqref{e:uv} so that
\begin{align*}
	f(y)
&\geq
	e_\lambda f(y+\lambda v)-\tfrac{1}{2\lambda}\norm{y+\lambda v-y}^2\\
&\geq
	e_\lambda f(x+\lambda v)+\ip{\nabla e_\lambda f(x+\lambda v)}{y-x}-\tfrac{1}{2\lambda}\norm{x+\lambda v-x}^2\\
&=e_\lambda f(x+\lambda v)-\tfrac{1}{2\lambda}\norm{x+\lambda v-x}^2+\ip{v}{y-x}\\
&=
	f(x)+\ip{v}{y-x},
\end{align*}
where the second inequality owes to \eqref{eq:cvx env}, the first equality by \eqref{eq:env and prox}, and the second
equality by \eqref{eq:env and prox1}.
Appealing to \cite[Theorem 1(b)$\Leftrightarrow$(c)]{rock-vietnam} completes the proof.

``\ref{thm:coco env grad:: firmly nonexap prox}$\Leftrightarrow$\ref{thm:coco env grad:: coco}'':
Observe that \ref{thm:coco env grad:: coco} amounts to \(\Id-P_\lambda f=\lambda\nabla e_\lambda f\) is firmly nonexpansive on \(U_\lambda\), which is equivalent to \ref{thm:coco env grad:: firmly nonexap prox}.
\endproof

Observe that for a fixed $\lambda$ the same $U_{\lambda}$ works for Theorem~\ref{thm:coco env grad}\ref{thm:coco env grad:: mono}--\ref{thm:coco env grad:: coco},
but $U_{\lambda}$ relies on $\lambda$.
To illustrate Theorem~\ref{thm:coco env grad}, consider the following simple example.
\begin{example} Define $f:\R\rightarrow\R: x\mapsto |x|_{0}$. Then $f$ is prox-bounded with $\lambda_{f}=+\infty$,
and $f$ is prox-regular at $0$ for $0\in\partial f(0)$.  The following hold for every $\lambda>0$:
\begin{enumerate}
\item $f$ is variationally convex at $0$ for $0\in\partial f(0)$.
\item Since
$$\partial_{p}^\lambda f(x)=
\begin{cases}
\varnothing, & \text{ if $0<|x|<\sqrt{2\lambda}$;}\\
[\sqrt{-2/\lambda},\sqrt{2/\lambda}], & \text{ if $x=0$;}\\
0, & \text{ if $|x|\geq\sqrt{2\lambda}$,}
\end{cases}
$$
for $U_{\lambda}=(-\sqrt{2\lambda},\sqrt{2\lambda})$ we have
$\gra\partial_{p}^\lambda f\cap W_{\lambda}=\{0\}\times (-\sqrt{2/\lambda},\sqrt{2/\lambda})$,
so monotone.

\item Since
$$P_{\lambda}f(x)=\begin{cases}
x, & \text{ if $|x|>\sqrt{2\lambda}$;}\\
\{0, \sqrt{2\lambda}\}, & \text{ if $x=\sqrt{2\lambda}$;}\\
0, & \text{ if $|x|<\sqrt{2\lambda}$,}\\
\{0, -\sqrt{2\lambda}\}, & \text{ if $x=-\sqrt{2\lambda}$,}
\end{cases}
$$
we have $P_{\lambda}f(x)=0$ on $U_{\lambda}=(-\sqrt{2\lambda},\sqrt{2\lambda})$, so firmly nonexpansive.
\item The Moreau envelope
$$e_{\lambda}f(x)=\begin{cases}
1, & \text{ if $|x|\geq \sqrt{2\lambda}$;}\\
\frac{1}{2\lambda}x^2, & \text{ if $|x|\leq\sqrt{2\lambda}$}
\end{cases}
$$
is convex on $U_{\lambda}=(-\sqrt{2\lambda},\sqrt{2\lambda}).$
\item $\nabla e_{\lambda}f(x)=x/\lambda$
on $U_{\lambda}=(-\sqrt{2\lambda},\sqrt{2\lambda})$, and
is $\lambda$ cocoercive.
\end{enumerate}
See \cite[Example 4.2]{wang23level}, and Examples~\ref{example:coincidience} and \ref{example:fail coin} in Section~\ref{s:integ} for details on $\partial_{p}^{\lambda}f, e_{\lambda}f$ and $P_{\lambda}f$.
\end{example}

\begin{remark}
Some comments are in order.
\vspace*{-2mm}
\begin{enumerate}
	\item
	Theorem \ref{thm:coco env grad} establishes a new correspondence between various classic notions in convex optimization in the absence of convexity.
% extending Fact \ref{fact:convex equivalence} significantly; see Figure \ref{fig:diagram} for the big picture.\
%	We primarily highlight the relation \ref{thm:coco env grad:: v cvx}$\Leftrightarrow$\ref{thmn:coco env grad:: max mono}$\Leftrightarrow$\ref{thm:coco env grad:: firmly nonexap prox}$\Leftrightarrow$\ref{thm:coco env grad:: averaged}$\Leftrightarrow$\ref{thm:coco env grad:: ne}, which reveals that various notions of nonexpansiveness of proximal operators reemerge in nonconvex setting locally \textit{if and only if} the underlying function is variationally convex \textit{if and only if the level proximal subdifferential is relatively monotone}.
	Hence variationally convex functions and local relative monotonicity of level proximal subdifferentials are natural counterparts of convex functions and monotonicity of Mordukhovich limiting subdifferentials in the absence of convexity;
see Figure~\ref{fig:diagram}.

\item One distinguished feature of Theorem~\ref{thm:coco env grad}\ref{envel:convex spec}$\Rightarrow$\ref{thm:coco env grad:: v cvx} is that \emph{only
    one $\lambda$} is required instead of all $0<\lambda<\bar{\lambda}$; see \cite[Theorem 3.2]{Khanh2023variational}.
    This is exactly what we need
    in Corollary~\ref{c:proxreg:neigh} and Theorem~\ref{cor:application} for
    variational sufficiency.
\item
 Recently, in \cite{khanh-phat} among many other results
 Kanhn, Khoa, Mordukhovich and Phat have extended the equivalence of
Theorem~\ref{thm:coco env grad}\ref{thm:coco env grad:: v cvx}$\Leftrightarrow$\ref{envel:convex}
to infinite-dimensional spaces.

	%Akin to the observation in \cite[Remark 3.3]{Khanh2023variational}, the reemergence of classic notions in convex optimization investigated in Theorem \ref{thm:coco env grad} opens the door for a promising new venue to analyze first-order methods in the absence of convexity---
	%\item
%%	Despite Theorem~\ref{thm:coco env grad} is established with a fixed parameter \(0<\lambda<\bar\lambda\),
%%	we note that Theorem~\ref{thm:coco env grad}\ref{envel:convex} can be strengthened to ``to \textit{every} \(0<\lambda<\bar\lambda\) corresponds to an open convex neighborhood \(U_\lambda\) of \(\bar x+\lambda\bar v\) on which \(e_\lambda f\) is convex''.
%%	The same can be said about every condition in Theorem~\ref{thm:coco env grad}.	
%	Although conditions in Theorem~\ref{thm:coco env grad} emphasize \textit{every} \(0<\lambda<\bar \lambda\), we note that in fact
%	\(
%	(\exists0<\lambda<\bar\lambda)~e_\lambda f\text{ is convex}\Leftrightarrow\ref{envel:convex}
%	\).
%	Clearly \ref{envel:convex} implies the former condition, from which an identical argument in \cite{Khanh2023variational} but with a fixed \(\lambda\) proves \ref{thm:coco env grad:: v cvx}.
\item Variational convexity implies prox-regularity, but not vice versa. For example, the function
$f:\Rn\rightarrow\Rn: x\mapsto -\|x\|^2/2$ is prox-regular but not variationally convex at every $x\in\Rn$.
\end{enumerate}
\end{remark}

Next we turn to characterizations of variational strong convexity, again coined by Rockafellar.

\begin{defn}\label{def:svc}\emph{\cite[Definition 2]{rock-vietnam}}
	A lsc function $f:\Rn\rightarrow\overline{\R}$ is called variationally strongly convex at $\bar x$
for $\bar v\in\partial f(\bar x)$ with modulus \(\sigma>0\),
	if for some convex neighborhood $U\times V$ of $(\bar x, \bar v)$ there exist a lsc \(\sigma\)-strongly convex function $\varphi\leq f$ on $U$ and a number $\varepsilon>0$ such that
	$[U_{\varepsilon}\times V]\cap \gph\partial f=[U\times V]\cap \gph \partial \varphi\text { and } f(x)=\varphi(x) \text{ at the common element $(x,v)$,} $
	where $U_{\varepsilon}:=\{x\in U|\ f(x)< f(\bar x)+\varepsilon\}$.
\end{defn}	

Below is the second main result of this section that establishes the equivalence
among variational strong convexity, local relative strong monotonicity of level proximal subdifferential, and
local cocoercivity of
proximal mapping with precise moduli.
The new features of our result are that we use
local properties of proximal mappings and level proximal subdifferentials.
The result complements \cite[Theorem 4.4]{Khanh2023variational}
which quantifies variational strong convexity via Moreau envelopes only; and
significantly improve \cite[Lemma 4.13]{planidenw} which quantifies strong convexity
via proximal mappings in a global setting.

\begin{theorem}[characterizations of variational strong convexity]\label{thm:variational str cvx}
	Suppose that $f:\mathbb{R}^n\to\overline{\mathbb{R}}$ is proper, lsc and prox-bounded
	with threshold $\lambda_{f}>0$,
	and that $f$ is prox-regular at $\bar x$ for $\bar v\in\partial f(\bar x)$.
	%{\color{gray}
		%Then for all $\lambda>0$ sufficiently small there exists an open convex neighborhood $U_{\lambda}$
		%of $\bar{x}+\lambda \bar{v}$ on which the following are equivalent:
		%}%
	Then there exists \(0<\bar\lambda\leq\lambda_{f}\) such that $(\forall 0<\lambda<\bar{\lambda})\ P_{\lambda}f(\barx+\lambda\barv)=\barx$, and $\exists U_{\lambda}$
an open convex neighborhood of \(\bar x+\lambda\bar v\) for which the following equivalent properties hold for \(\sigma>0\):
	\begin{enumerate}
		\item\label{thm:variational str cvx:str cvx}
		\(f\) is variationally strongly convex at \(\bar x\) for \(\bar v\in\partial f(\bar x)\) with modulus
\(\sigma\);
		\item\label{thm:variational str cvx:str mono}
		\((\forall 0<\lambda<\bar\lambda)\)
		\(\partial_p^\lambda f\) is \(\sigma\)-strongly monotone relative to \(W_\lambda\);
		\item\label{thm:variational str cvx:coco}
		\((\forall 0<\lambda<\bar\lambda)\)
		\(P_\lambda f\) is \((1+\lambda\sigma)\)-cocoercive on \(U_\lambda\);
		\item\label{thm:variational str cvx:lip}
		\((\forall 0<\lambda<\bar\lambda)\)
		\(P_\lambda f\) is \(\tfrac{1}{1+\lambda\sigma}\)-Lipschitz on \(U_\lambda\);
% and $P_{\lambda}f(\barx+\lambda\barv)=\barx$;
		\item\label{thm:variational str cvx:env}
		\((\forall0<\lambda<\bar\lambda)\)
		\(e_\lambda f\) is \(\tfrac{\sigma}{1+\sigma\lambda}\)-strongly convex on \(U_\lambda\);
	\end{enumerate}
	where  $W_\lambda :=\{(x,v)\in\Rn\times\Rn|\ x+\lambda v\in U_\lambda\}$.
\end{theorem}
\proof
By
\cite[Proposition 5.3]{thibault}, there exists \(0<\bar\lambda\leq \lambda_{f}\) such that to every \(0<\lambda<\bar\lambda\) corresponds an open convex neighborhood $U_\lambda$ of $\bar x+\lambda\bar v$ on which:
\begin{equation}\label{e:proxreg1}
	\text{$P_{\lambda}f$ is monotone, single-valued and Lipschitz continuous; $P_{\lambda}f(\bar x+\lambda{\bar v})=\bar x$; and}
\end{equation}
%$e_\lambda f$ is differentiable with
\begin{equation}\label{e:proxreg2}
\text{$e_\lambda f$ is differentiable with }	\nabla e_{\lambda}f=\lambda^{-1}[\Id-P_\lambda f].
\end{equation}
See also \cite[Proposition 13.37]{rockafellar_variational_1998} for $\bar v=0$.
In particular, $\bar{v}\in\partial_{p}^{\lambda}f(\bar{x})\subseteq\hat{\partial}f(\bar{x})$
by Fact~\ref{fact: prox identity}.
%\begin{itemize}[left=0pt, labelsep=0em, itemindent=0pt, label={}]

``\ref{thm:variational str cvx:str cvx}$\Rightarrow$\ref{thm:variational str cvx:str mono}":
%	Our proof is inspired by \cite[Theorem 3.2]{Khanh2023variational}.
Since $f$ is variationally strongly convex at $\bar x$ for $\bar v\in\hat{\partial}f(\bar x)$
with modulus $\sigma$,
by \cite[Theorem 2]{rock-vietnam}, there exist a neighborhood $U$ of $\bar x$, a neighborhood $V$ of $\bar v$,
and $\varepsilon>0$ such that
\begin{equation}\label{e:strong:var}
\gph \partial f\cap (U_{\varepsilon}\times V) \text{ is $\sigma$-strongly monotone,}
\end{equation}
where $U_{\varepsilon}:=\{x\in U|\ f(x)<f(\bar{x})+\varepsilon\}.$

Now \eqref{e:proxreg1} and \eqref{e:proxreg2} show that when
$x\in U_{\lambda}$, we have
\begin{align}
\lim_{x\rightarrow \bar{x}+\lambda\bar{v}}(x-\lambda \nabla e_{\lambda}f(x))& =\lim_{x\rightarrow \bar{x}+\lambda\bar{v}}P_{\lambda}f(x)= P_{\lambda}(\bar{x}+\lambda\bar{v})=\bar{x},\\
\lim_{x\rightarrow \bar{x}+\lambda\bar{v}}\nabla e_{\lambda}f(x) &= \lim_{x\rightarrow \bar{x}+\lambda\bar{v}}\frac{x-P_{\lambda}f(x)}{\lambda}=
 \frac{\bar{x}+\lambda\bar{v}-\bar{x}}{\lambda}=\bar{v},
\end{align}
and
\begin{align}
\lim_{x\rightarrow \bar{x}+\lambda\bar{v}}f(x-\lambda \nabla e_{\lambda}f(x)) &=\im_{x\rightarrow \bar{x}+\lambda\bar{v}} f(P_{\lambda}f(x))=\lim_{x\rightarrow \bar{x}+\lambda\bar{v}}\Big(e_{\lambda}f(x)-\tfrac{1}{2}\|x-P_{\lambda}f(x)\|^2\Big)\\
& =e_{\lambda}f(\bar{x}+\lambda\bar{v})-\tfrac{1}{2}\|\bar{x}+\lambda\bar{v}-\bar{x}\|^2
=f(\bar{x}).
\end{align}
Therefore, we can choose $U_{\lambda}$ smaller, if necessary, such that
\begin{equation}\label{e:gradenv1}
(\forall x\in U_{\lambda})\ x-\lambda \nabla e_{\lambda}f(x)\in U, \nabla e_{\lambda}f(x) \in V, \text{ and }
f(x-\lambda \nabla e_{\lambda}f(x))<f(\bar x)+\varepsilon;
\end{equation}
moreover,
$x-\lambda \nabla e_{\lambda}f(x)=P_{\lambda}f(x)$ implies by Fact~\ref{fact: prox identity} that
\begin{equation}\label{e:gradenv2}
(x-\lambda \nabla e_{\lambda}f(x),\nabla e_{\lambda}f(x))\in\gph \partial_{p}^{\lambda}f\subset\gph \partial f.
\end{equation}
Combining \eqref{e:gradenv1} and \eqref{e:gradenv2} we arrive at
\begin{equation}
\{(x-\lambda \nabla e_{\lambda}f(x),\nabla e_{\lambda}f(x))|\ x\in U_{\lambda}\}\subseteq \gph \partial f\cap
(U_{\varepsilon}\times V),
\end{equation}
which implies by \eqref{e:strong:var} that
$\{(x-\lambda \nabla e_{\lambda}f(x),\nabla e_{\lambda}f(x))|\ x\in U_{\lambda}\}$ is $\sigma$-strongly
monotone. The result then follows from Lemma~\ref{l:minty}\ref{i:function}.

%Prox-regularity implies that there exists \(\gamma>0\) such that set-valued operator \(T:\Rn\rightrightarrows\Rn\) defined by
%	\begin{align*}
%	T(x) :=
%	\begin{cases}
%		\{v\in\partial f(x): \norm{v-\bar v}<\gamma\}, &\text{ if }\norm{x-\bar x}<\gamma\text{ and }f(x)< f(\bar x)+\gamma;\\
%		\varnothing, &\text{ otherwise},
%	\end{cases}
%	\end{align*}
%	satisfies \(\nabla e_\lambda f(x)=\left(\lambda\Id+T^{-1}\right)^{-1}(x)\) for \(x\in U_\lambda\), and therefore for \(x\in U_\lambda\),
%	\(
%		\left(\nabla e_\lambda f(x), x-\lambda\nabla e_\lambda f(x)\right)\in\gph T.
%	\)
%	Invoking Lemma \ref{l:minty} yields
%	\begin{equation}\label{gph inclusion}
%		\gph\partial_p^\lambda f\cap W_\lambda\subseteq\gph T.
%	\end{equation}
%	Shrinking \(\gamma\) if necessary, condition \ref{thm:variational str cvx:str cvx} together with \cite[Theorem 2]{rock-vietnam} enforces \(\gra T\) to be \(\sigma\)-strongly monotone, and hence in particular \(\gph\partial_p^\lambda f\cap W_\lambda\).
	
``\ref{thm:variational str cvx:str mono}$\Leftrightarrow$\ref{thm:variational str cvx:coco}'':
	Apply Lemma \ref{lem:local coco Minty}.
	
``\ref{thm:variational str cvx:coco}$\Rightarrow$\ref{thm:variational str cvx:lip}'':
	For every \(x_1,x_2\in U_\lambda\),
	\(
		\norm{P_\lambda f(x_1)-P_\lambda f(x_2)}^2
	\leq
		\tfrac{1}{1+\lambda\sigma}\ip{x_1-x_2}{P_\lambda f(x_1)-P_\lambda f(x_2)}
	\leq
		\tfrac{1}{1+\lambda\sigma}\norm{P_\lambda f(x_1)-P_\lambda f(x_2)}\norm{x_1-x_2}
	\),
	completing the proof.
	
``\ref{thm:variational str cvx:lip}$\Rightarrow$\ref{thm:variational str cvx:env}'':
	In view of \eqref{e:proxreg2}, for every \(x_1,x_2\in U_\lambda\),
	\(
		\ip{\nabla e_\lambda f(x_1)-\nabla e_\lambda f(x_2)}{x_1-x_2}
	=
		\tfrac{1}{\lambda}\norm{x_1-x_2}^2
	-
		\tfrac{1}{\lambda}\ip{x_1-x_2}{P_\lambda f(x_1)-P_\lambda f(x_2)}
	\geq
		\tfrac{1}{\lambda}\left(1-\tfrac{1}{1+\sigma\lambda}\right)\norm{x_1-x_2}^2
=\tfrac{\sigma}{1+\sigma\lambda}\norm{x_1-x_2}^2
	\), where the inequality follows from the Cauchy-Schwartz inequality and \ref{thm:variational str cvx:lip}.
This justifies that $e_{\lambda}f$ is $\sigma/(1+\sigma\lambda)$-strongly convex on $U_{\lambda}$;
see, e.g., \cite[Excersie 12.59]{rockafellar_variational_1998} or \cite[Theorem 5.24]{beck2017}.
	
``\ref{thm:variational str cvx:env}$\Rightarrow$\ref{thm:variational str cvx:str cvx}'':
	Apply \cite[Theorem 4.4]{Khanh2023variational}.
\endproof

Observe that for a fixed $\lambda$ the same $U_{\lambda}$ works for Theorem~\ref{thm:variational str cvx}\ref{thm:variational str cvx:str mono}--\ref{thm:variational str cvx:env},
but $U_{\lambda}$ relies on $\lambda$.

For $\lambda>0$ and a subset $B\subseteq \Rn$, define the \emph{localized proximal mapping} by
$$(\forall x\in B)\ P_{\lambda}^{B}f(x) \coloneqq \argmin_{y\in B}\{f(y)+\frac{1}{2\lambda}\|y-x\|^{2}\}.$$
Significant progress has been made by Rockafellar in \cite[Theorem 3]{rock2021} in which he
showed:``Let $f$ be prox-regular for $(\bar{x},0)\in\gph \partial f$. Then the variational convexity of $f$ for $(\bar{x},0)$ is not just sufficient but also necessary for the existence of an open ball
$B$ centered at $\bar{x}$ on which $P^{B}_{\lambda}f$ is firmly nonexpansive."
Our Theorem~\ref{thm:coco env grad} utilizes $P_{\lambda}f$ instead of $P_{\lambda}^{B}f$ and provide more new characterizations. In addition, the following result shows that
\cite[Theorem 3]{rock2021} can be recovered.
\begin{theorem}
Let \(f:\Rn\to\overline{\R}\) be proper, lsc and prox-bounded with prox-bound $\lambda_{f}>0$.
Suppose that $f:\Rn\rightarrow\overline{\R}$ is variationally convex at $(\bar{x}, 0)\in\gph \partial f$.
% and that
%	$f$ is prox-bounded.
Then for all $\lambda>0$ sufficiently small there exists an open ball
	$B$ centered at $\olx$ on which
	$P_{\lambda}f=P_{\lambda}^B f$ and $P_{\lambda}f$ is firmly nonexpansive on $B$.
\end{theorem}
\proof Since $f$ is variationally convex at $\olx$ for $0\in\partial f(\olx)$, it is prox-regular
$\olx$ for $0\in\partial f(\olx)$.
Then there exists a ball $B_{1}$
	centered at $\olx$
	on which $P_{\lambda}f$ is single-valued, Lipschitz and
	$P_{\lambda}f(\olx)=\olx$; see, e.g.,
\cite[Proposition 13.37]{rockafellar_variational_1998}.
	By Theorem~\ref{thm:coco env grad},
	we can find an open ball $B_{2}$ centered at $\olx$ and $B_{2}\subseteq B_{1}$ on which
	$P_{\lambda}f$ is firmly nonexpansive. Because
	$$(\forall x\in B_{2})\ \|P_{\lambda}f(x)-\olx\|=\|P_{\lambda}f(x)-P_{\lambda}f(\olx)\|\leq \|x-\olx\|,$$
	this implies that $P_{\lambda}f(x)\in B_{2}$ if $x\in B_{2}$. Since $P_{\lambda}f(x)$ is a global unique minimizer
	of
	$$y\mapsto f(y)+\frac{1}{2\lambda}\|y-x\|^2,$$
	it must be a unique minimizer on $B_{2}$. With $B \coloneqq B_{2}$, we have
	$(\forall x\in B)\ P_{\lambda}f(x)=P_{\lambda}^{B}f(x)$.
% for every $x\in B$.
\endproof

%For $\lambda>0$ and a subset $B\subseteq \Rn$, define the \emph{localized proximal mapping} by
%$$(\forall x\in B)\ P_{\lambda}^{B}f(x) :=\argmin_{y\in B}\{f(y)+\frac{1}{2\lambda}\|y-x\|^{2}\}.$$
%
%\begin{theorem} Suppose that $f:\Rn\rightarrow\overline{\R}$ is variationally convex at $(\bar{x}, 0)\in\gph \partial f$.
%% and that
%%	$f$ is prox-bounded.
%Then for all $\lambda>0$ sufficiently small there exists an open ball
%	$B$ centered at $\olx$ on which
%	$P_{\lambda}f=P_{\lambda}^B f$ and $P_{\lambda}f$ is firmly nonexpansive on $B$.
%\end{theorem}
%\proof Since $f$ is variationally convex at $\olx$ for $0\in\partial f(\olx)$, it is prox-regular
%$\olx$ for $0\in\partial f(\olx)$, in particular, prox-bounded. Then
%there exists a ball $B_{1}$
%	centered at $\olx$
%	on which $P_{\lambda}f$ is single-valued, Lipschitz and
%	$P_{\lambda}f(\olx)=\olx$; see, e.g.,
%\cite[Proposition 13.37]{rockafellar_variational_1998}.
%	By Theorem~\ref{thm:coco env grad},
%	we can find an open ball $B_{2}$ centered at $\olx$ and $B_{2}\subseteq B_{1}$ on which
%	$P_{\lambda}f$ is firmly nonexpansive. Because
%	$$(\forall x\in B_{2})\ \|P_{\lambda}f(x)-\olx\|=\|P_{\lambda}f(x)-P_{\lambda}f(\olx)\|\leq \|x-\olx\|,$$
%	this implies that $P_{\lambda}f(x)\in B_{2}$ if $x\in B_{2}$. Since $P_{\lambda}f(x)$ is a global unique minimizer
%	of
%	$$y\mapsto f(y)+\frac{1}{2\lambda}\|y-x\|^2,$$
%	it must be a unique minimizer on $B_{2}$. With $B :=B_{2}$, we have
%	$(\forall x\in B)\ P_{\lambda}f(x)=P_{\lambda}^{B}f(x)$.
%% for every $x\in B$.
%\endproof

Acute readers might be wondering whether one can lift Theorem~\ref{thm:coco env grad} to a global
setting. To this end, Lemma~\ref{thm:local Minty} is essential.

%\begin{corollary}
%Let $f:\Rn\to\overline{\R}$ be proper, lsc, and prox-bounded with threshold $\lambda_f>0$.
%Let $U\subseteq\Rn$ be a nonempty subset and define $(\forall 0<\lambda<\lambda_f)$ $W_\lambda=\{(x,v): x+\lambda v\in U\}$.
%The following are equivalent:
%\begin{enumerate}
%	\item $P_\lambda f$ is firmly nonexpansive on $U$;
%	\item $\partial_p^\lambda f$ is monotone relative to $W_\lambda$;
%	\item  $\partial_p^\lambda f$ is maximally monotone relative to $W_\lambda$.
%\end{enumerate}
%\end{corollary}
%\proof
%Appealing to Theorem \ref{thm:local Minty} with $A=\partial_p^\lambda f$ completes the proof.
%Note that $\dom J_{\lambda\partial_p^\lambda f}=\dom P_\lambda f=\Rn$, therefore the condition $U\subseteq\dom J_{\lambda A}$ is automatically satisfied.
%\endproof
\begin{theorem}\label{t:global:stuff}
Let $f:\Rn\to\overline{\R}$ be proper, lsc and prox-bounded with threshold $\lambda_f>0$, and let
$0<\lambda<\lambda_f$.
Then the following are equivalent:
\begin{enumerate}
	\item\label{i:sun1} $P_\lambda f$ is firmly nonexpansive on $\Rn$;
	\item\label{i:sun2} $\partial_p^\lambda f$ is monotone;
	\item\label{i:sun3} $\partial_p^\lambda f$ is maximally monotone;
\item\label{i:sun4} $e_{\lambda} f$ is convex on $\Rn$;
	\item\label{i:sun5} $f$ is convex on $\Rn$.
\end{enumerate}
\end{theorem}
\proof
``\ref{i:sun1}$\Leftrightarrow$\ref{i:sun2}$\Leftrightarrow$\ref{i:sun3}":
Apply Lemma~\ref{thm:local Minty} to $A :=\partial_{p}^{\lambda}f$ and $U :=\Rn$.

``\ref{i:sun1}$\Rightarrow$\ref{i:sun4}": Invoke Lemma~\ref{l:minty}\ref{i:function}
to obtain that $\nabla e_{\lambda}f$ is cocoercive on $\Rn$,
so monotone on $\Rn$.
This implies that $e_{\lambda}f$ is convex by \cite[Theorem 2.14(a)]{rockafellar_variational_1998}
or \cite[Proposition 17.7(iii)]{BC}.

``\ref{i:sun4}$\Rightarrow$\ref{i:sun5}": Apply \cite[Theorem 3.17]{wang2010Chebyshev}.
``\ref{i:sun5}$\Rightarrow$\ref{i:sun1}": This is well-known,
see, i.e., \cite[Proposition 12.28]{BC}.
% and \cite[Theorem 3.2]{luo2024various}.
\endproof

\begin{remark}
Theorems~\ref{thm:coco env grad} and \ref{t:global:stuff} together establish the beautiful and
new correspondence between convex functions and variationally convex functions
given in Figure~\ref{fig:diagram} on page~\pageref{convex-varia},
which is not known in the literature.
\end{remark}

%\begin{theorem} Suppose that $f:\Rn\rightarrow\overline{\R}$ is prox-regular at $(\bar{x}, 0)\in\gr \partial f$,
%	and that $f$ is prox-bounded.
%	Let $B$ be a open ball centered at $\olx$. Then for all $\lambda>0$ sufficiently small there exists an open ball
%	$U\subseteq B$ centered at $\olx$ on which
%	$P_{\lambda}f=P_{\lambda}^B f$ on $U$.
%\end{theorem}
%
%\proof Since $f$ is prox-regular at $\olx$ for $0\in\partial f(\olx)$, there exists a ball $U\subseteq B$
%	centered at $\olx$
%	on which $P_{\lambda}f$ is single-valued, Lipschitz with modulus $\ell>0$ and
%	$P_{\lambda}f(\olx)=\olx$; see, e.g., \cite[Proposition 13.37]{rockafellar_variational_1998}. Since
%	$$(\forall x\in U)\ \|P_{\lambda}f(x)-\olx\|=\|P_{\lambda}f(x)-P_{\lambda}f(\olx)\|\leq \ell \|x-\olx\|,$$
%	there exists $\delta>0$ such that $P_{\lambda}f(x)\in U\subseteq B$ if $\|x-\olx\|<\delta$.
%	Since $P_{\lambda}f(x)$ is a global unique minimizer
%	of $$y\mapsto f(y)+\frac{1}{2\lambda}\|y-x\|^2,$$
%	it must be a unique minimizer on $B$. Therefore,
%	$P_{\lambda}f(x)=P_{\lambda}^{B}f(x)$ whenever $\|x-\olx\|<\delta$.
%\endproof

\section{Algorithms for variationally convex functions}\label{s:algorithms}

The theory developed in Theorems~\ref{thm:coco env grad} and \ref{thm:variational str cvx} opens the door, in algorithm design, to deploying classical convex machinery locally, in the absence of convexity. In this section, we will propose local
proximal gradient algorithm and
Krasnosel'ski\u{i}-Mann algorithm for variationally convex functions, and establish their convergence respectively.

In a series of articles, Rockafellar \cite{rock-vietnam,rockafellar2023augmented,rockafellar2023convergence} explored the applications of variational convexity on proximal point and augmented Lagrangian methods. These work use the very
powerful notion of \emph{variational sufficiency}, which guarantees that {stationary points of variationally convex functions are local minimizers}; see \cite[Section 2]{rockafellar2023augmented}, \cite[Section 6]{Khanh2023variational}, and \cite{Wang2023}. We shall employ variational sufficiency below.

\begin{defn} Let $f:\Rn\rightarrow\OR$ be proper and lsc with $\inf_{\R^n}f$ finite. Consider the unconstrained
minimization problem
\begin{equation}\label{e:varia}
\min\{f(x)|\ x\in\Rn\}.
\end{equation}
We say that \eqref{e:varia}
satisfies the variational sufficiency condition at $\olx$ if $f$ is variationally convex at $\olx$
for $0\in\partial f(\olx)$.
\end{defn}
The beauty of variational sufficiency is:
\begin{fact}\label{f:localmin:v}
 Suppose that \eqref{e:varia} satisfies the variational sufficiency condition at $\olx$. Then
$\olx$ is a local minimizer of $f$.
\end{fact}
\begin{proof} See \cite[p. 550]{rock-vietnam} or apply \cite[Theorem 1(c)]{rock-vietnam}.
\end{proof}

Consider the structured minimization problem
$$
                \min\left\lbrace
                   f(x)+g(x)|\ x\in\R^n
                \right\rbrace
            $$
            where $f:\R^n\rightarrow\OR$ is a proper, lsc and prox-bounded function; and $g:\R^n\rightarrow\R$
            is a $C^1$ function whose gradient $\nabla g$ is Lipschitz continuous.
For $\lambda>0$, the proximal gradient operator is
\begin{equation}\label{e:prox:oper}
T:=P_{\lambda}f(\Id-\lambda\nabla g).
\end{equation}
 Level-proximal subdifferential provides
a new explanation of fixed points of $T$.
\begin{lemma}\label{l:prox:fixset}
 With $T$ given by \eqref{e:prox:oper}, we have
\begin{equation}\label{e:prox:fix}
\Fix T=\menge{x\in\R^n}{0\in \nabla g(x)+\partial_{p}^{\lambda}f(x)}.
\end{equation}
\end{lemma}
\begin{proof} Indeed,
\begin{align}
x\in Tx & \Leftrightarrow x\in P_{\lambda}f(x-\lambda\nabla g(x))
 \Leftrightarrow x-\lambda \nabla g(x) \in (\Id+\lambda\partial_{p}^{\lambda}f)(x) \label{e:cone}\\
& \Leftrightarrow 0\in \lambda\nabla g(x)+\lambda\partial_{p}^{\lambda}f(x)
\Leftrightarrow 0\in \nabla g(x)+\partial_{p}^{\lambda}f(x). \label{e:ctwo}
\end{align}
\end{proof}

Lemma~\ref{l:prox:fixset} shows that the proximal gradient method finds precisely $x$ such that
$0\in \nabla g(x)+\partial_{p}^{\lambda}f(x)$. We call such an $x$ as a \emph{$\lambda$-level $p$-critical point}
of $f+g$.

\begin{remark}
A critical point $x$ of
$f+g$ is usually defined by $0\in \nabla g(x)+\partial f(x)$, which might not be
a $\lambda$-level $p$-critical point,
since $\partial_{p}^{\lambda}f(x)\subseteq\partial f(x)$. When $f$ is convex, they coincide.
\end{remark}

We start with a sum rule on variationally convex functions.
\begin{lemma}\label{l:c1sum}
 Suppose that $f:\R^n\rightarrow\OR$ is proper, lsc and $\barx\in\dom f$,
and that $g:\R^n\rightarrow\RX$ is convex and $C^1$ on a neighborhood of $\barx$. If
$f$ is variationally convex at $\barx$ for $\barv\in
\partial f(\barx)$, then
$f+g$ is variationally convex at $\barx$ for $\barv+\nabla g(\barx)\in\partial(f+g)(\barx)$.
\end{lemma}
\begin{proof}
Since $f$ is variationally convex at $\barx$ for $\barv$, by \cite[Theorem 1]{rock-vietnam},
there exist convex neighborhoods
$U$ of $\barx$ and $V$ of $\barv$, and $\varepsilon>0$ such that
$V-\barv=-(V-\barv)$,
\begin{equation}\label{e:rock:char}
(\forall x\in U)(\forall (u,v)\in\gra \partial f\cap (U_{\varepsilon}\times V))\ f(x)\geq f(u)+\scal{v}{x-u},
\end{equation}
where $U_{\varepsilon}:=\menge{u\in U}{f(u)<f(\barx)+2\varepsilon}$.
Choose a smaller neighborhood $U$, if necessary, to make sure that $g$ is convex and $C^1$ on $U$.

Next, because $g$ is $C^1$, we may choose a smaller
neighborhood $U_{1}$ of $\barx$ such that $U_{1}\subseteq U$,
$$(\forall u\in U_{1})\ |g(u)-g(\barx)|< \varepsilon, \text{ and }
 \nabla g(u)\in \nabla g(\barx)+\frac{V-\barv}{2}.$$
For
\begin{align}
v &\in\partial(f+g)(u) = \partial f(u)+\nabla g(u),\label{e:fplusg1}\\
(u,v) & \in U_{1,\varepsilon}\times \bigg(\barv+\nabla g(\barx)+\frac{V-\barv}{2}\bigg), \text{ where }\\
U_{1,\varepsilon} & :=\menge{u\in U_{1}}{(f+g)(u)<(f+g)(\barx)+\varepsilon},\label{e:fplusg2}
\end{align}
we have: $v-\nabla g(u)\in\partial f(u)$,
\begin{equation}
v-\nabla g(u)= v-(\barv+\nabla g(\barx))-(\nabla g(u)-\nabla g(\barx))+\barv
\in\frac{V-\barv}{2}+\frac{V-\barv}{2}+\barv=V-\barv+\barv=V,
\end{equation}
and
\begin{equation}
f(u)-f(\barx)=[(f+g)(u)-(f+g)(\barx)]-(g(u)-g(\barx))<\varepsilon+\varepsilon=2\varepsilon.
\end{equation}
Thus, applying \eqref{e:rock:char} to $(u,v-\nabla g(u))\in\gra \partial f$ we obtain
\begin{equation}\label{e:function:f}
(\forall x\in U_{1})\ f(x)\geq f(u)+\scal{v-\nabla g(u)}{x-u}.
\end{equation}
Since $g$ is convex on $U_{1}$,
\begin{equation}\label{e:function:g}
(\forall x\in U_{1})\ g(x)\geq g(u)+\scal{\nabla g(u)}{x-u}.
\end{equation}
Adding \eqref{e:function:f} and \eqref{e:function:g} yields
\begin{equation}
(\forall x\in U_{1})\ f(x)+g(x)\geq f(u)+g(u)+\scal{v}{x-u},
\end{equation}
for $(u,v)$ verifying \eqref{e:fplusg1}--\eqref{e:fplusg2}.
Using \cite[Theorem 1]{rock-vietnam} again, we conclude that
$f+g$ is variationally convex at $\barx$ for $\barv+\nabla g(\barx)$.
\end{proof}

The next result on compositions of averaged mappings improves \cite[Proposition 4.44]{BC},
in which they assume $T_{1}, T_{2}:D\rightarrow D$. Since the proof is very similar, we omit it.
\begin{lemma}\label{l:composition:ave}
Let $D_{1}, D_{2}$ be nonempty subsets of $\R^n$, let $T_1:D_{1}\rightarrow D_{2}$
be $\alpha_{1}$-averaged, and let $T_{2}: D_{2}\rightarrow \R^n$ be $\alpha_{2}$-averaged,
where $\alpha_{i}\in ]0,1[$ for $i=1, 2$. Set $$T:=T_{2}T_{1} \text{ and }
\alpha:=\frac{\alpha_{1}+\alpha_{2}-2\alpha_{1}\alpha_{2}}{1-\alpha_{1}\alpha_{2}}.$$
Then $\alpha\in ]0,1[$, and $T$ is $\alpha$-averaged on $D_{1}$.
\end{lemma}

Utilizing Theorems~\ref{thm:coco env grad} and \ref{thm:variational str cvx}, we propose a localized
proximal gradient algorithm.
\begin{theorem}[proximal gradient algorithm under variational convexity]\label{t:f:plusg}
 Suppose that $f:\R^n\rightarrow\OR$ is proper, lsc and prox-bounded with threshold $\lambda_{f}>0$, and that
$g:\R^n\rightarrow\R$ is convex, differentiable, and has a Lipschitz continuous gradient. Let $f$ be
 variationally convex at $\bar{x}$ for $-\nabla g(\bar x)\in\partial f(\bar x)$.
Define
$T:=P_{\lambda}f(\Id-\lambda \nabla g)$. Then
$\barx$ is a local minimizer of $f+g$, and
$\exists {0<\bar{\lambda}}\leq\lambda_{f}$ such that
\begin{enumerate}
\item\label{i:ave1}
 $\forall \lambda\in (0,\bar{\lambda})$ there exists a
neighborhood $U_{\lambda}$ of ${\bar{x}}$ on which $T:U_{\lambda}\rightarrow U_{\lambda}$,
$T$ is averaged on $U_{\lambda}$, and $T(\bar{x})=\bar{x}$.
\item\label{i:ave2}
Let $x_{1}\in U_{\lambda}$. Define the
sequence $(x_{k})_{k\in\N}$ by $x_{k+1}:=Tx_{k}$. Then $(x_{k})_{k\in\N}$ converges
to $x^*\in U_{\lambda}$,  a $\lambda$-level $p$-critical point of $f+g$.
\item\label{i:ave3}
If, in addition, $f$ is variationally strongly convex at $\barx$ for
$-\nabla g(\bar x)\in\partial f(\barx)$,
then, by choosing a smaller neighborhood $U_\lambda$,
the sequence given in \ref{i:ave2} converges linearly to $\barx$.
\end{enumerate}
\end{theorem}

\begin{proof} Since
$f$ is variationally convex at $\bar{x}$ for $-\nabla g(\bar x)\in\partial f(\bar x)$, Lemma~\ref{l:c1sum}
assures that $f+g$ is variationally convex at $\barx$ for $0\in\partial(f+g)(\barx)$.
It follows from Fact~\ref{f:localmin:v} that $\bar{x}$ is a local minimum of
$f+g$.

\ref{i:ave1}:
By Theorem~\ref{thm:coco env grad}, $\exists 0<\bar{\lambda}\leq \lambda_{f}$ such
that for $0<\lambda<\bar{\lambda}$ the mapping
$P_{\lambda}f$ is firmly nonexpansive on
a neighborhood of $\bar{x}-\lambda\nabla g(\bar{x})$, say
$\Ball_{\delta}(\bar{x}-\lambda\nabla g(\bar{x}))$; and that
$T(\bar{x})=P_{\lambda}f(\bar{x}-\lambda\nabla g(\bar{x}))=\bar{x}$.
Because $\nabla g$ is Lipschitz, by choosing a smaller $\bar{\lambda}$, we can make
$\Id-\lambda\nabla g$ firmly nonexpansive on $\R^n$ by the Baillon-Haddad theorem
\cite[Corollary 18.17]{BC}.
%In view of
%\eqref{e:cone}-\eqref{e:e:ctwo}, $T(\bar{x})=\bar{x}$.
For $x\in U_{\lambda} :=\Ball_{\delta}(\bar x)$, we have
$$\|(x-\lambda \nabla g(x))-
(\bar{x}-\lambda\nabla g(\bar{x}))\|\leq \|x-\bar{x}\|\leq \delta,$$
so that
\begin{align}
\|Tx-\bar{x}\| &=\|P_{\lambda}f(x-\lambda\nabla g(x))-
P_{\lambda}(\bar{x}-\lambda\nabla g(\bar{x})\|\\
&\leq \|(x-\lambda\nabla g(x))-(\bar{x}-\lambda\nabla g(\bar{x})\|\leq \delta
\end{align}
because $\Id-\lambda\nabla g$ is nonexpansive on $\R^n$, and $P_{\lambda}f$ is
nonexpansive on $\Ball_{\delta}(\bar{x}-\lambda\nabla g(\bar{x}))$. Thus, $T:U_{\lambda}\rightarrow
U_{\lambda}$.
Being a composition of firmly nonexpansive ($1/2$-averaged) mappings,
$T$ is $2/3$-averaged on $U_{\lambda}$
by Lemma~\ref{l:composition:ave}.

\ref{i:ave2}: Applying \cite[Proposition 5.16]{BC} to $T$ on $U_{\lambda}$ gives that
$x_{k}\rightarrow x^*\in U_{\lambda}$ and $x^*=Tx^*$. It suffices to use Lemma~\ref{l:prox:fixset}.
%The proof is
%complete by using \eqref{e:cone}--\eqref{e:ctwo}.

\ref{i:ave3}: By Theorem~\ref{thm:variational str cvx}\ref{thm:variational str cvx:lip}
$\exists \bar{\lambda}$ such that for $0<\lambda<\bar{\lambda}$ the mapping
$P_{\lambda}f$ is a Banach contraction on
a neighborhood of $\bar{x}-\lambda\nabla g(\bar{x})$, say
$\Ball_{\delta}(\bar{x}-\lambda\nabla g(\bar{x}))$; and that
$T(\bar{x})=P_{\lambda}f(\bar{x}-\lambda\nabla g(\bar{x}))=\bar{x}$.
Because $\nabla g$ is Lipschitz, by choosing a smaller $\bar{\lambda}$, we can make
$\Id-\lambda\nabla g$ firmly nonexpansive on $\R^n$.
Then $T: U_{\lambda}\rightarrow U_\lambda$ is a Banach contraction with
$U_{\lambda}:=\Ball_{\delta}(\bar x)$.
It follows from the Banach fixed point theorem \cite[Theorem 1.50]{BC} or \cite[Theorem 5.1-2]{kreyszig89} that
for every $x_{1}\in U_\lambda$ the sequence
$(x_{k})_{k\in\NN}$ must converge linearly to $\barx$.
\end{proof}
\begin{remark} Theorem~\ref{t:f:plusg} still holds if the global assumption on $g$ is replaced by a local one,
i.e., $g$ being
convex and differentiable
with $\nabla g$ cocoercive on a neighborhood of $\bar{x}$. The proof goes through because
$\Id-\lambda \nabla g$ is firmly nonexpansive on the neighborhood when $\lambda$ is sufficiently small.
Local cocoercivity can be realized by using Lemma~\ref{lem:descent to lip grad} in section~\ref{s:lips}.
\end{remark}

When $f$ is prox-regular on a
neighborhood (e.g., $f$ being strongly amenable at $\bar{x}$
\cite[Proposition 13.32]{rockafellar_variational_1998}),
amazingly the iterates
given by Theorem~\ref{t:f:plusg}\ref{i:ave2} actually converge to a
local minimizer!
\begin{corollary}\label{c:proxreg:neigh}
 Suppose that $f:\R^n\rightarrow\OR$ is proper, lsc and prox-bounded with threshold $\lambda_{f}>0$, and that
$g:\R^n\rightarrow\R$ is convex, differentiable, and has a Lipschitz continuous gradient. Let $f$ be
 variationally convex at $\bar{x}$ for
 $-\nabla g(\bar x)\in\partial f(\bar x)$ and let
 $f$ be prox-regular on a neighborhood
 of $\barx$. Define
$T:=P_{\lambda}f(\Id-\lambda \nabla g)$. Then
$\exists 0<{\bar{\lambda}}\leq\lambda_{f}$ such that
 $\forall \lambda\in (0,\bar{\lambda})$ there exists a
neighborhood $U_{\lambda}$ of ${\bar{x}}$ on which
for every $x_{1}\in U_{\lambda}$, the
sequence $(x_{k})_{k\in\N}$ given by $x_{k+1} :=Tx_{k}$ converges
to $x^*\in U_{\lambda}$,  a local minimizer of $f+g$.
\end{corollary}
\begin{proof} Observe that $f+g$ is prox-bounded because of $f$ being prox-bounded and $g$ being convex.
By assumption, we can choose a neighborhood $U$ of $\barx$ on which $f$ is prox-regular.
Because $g$ is $C^{1+}$, $f+g$ is prox-regular on $U$;
see, e.g., \cite[Exercise 13.35]{rockafellar_variational_1998}.
Theorem~\ref{t:f:plusg} shows that $\exists {\bar{\lambda}}$ such that
 $\forall \lambda\in (0,\bar{\lambda})$ there exists a
neighborhood $U_{\lambda}$ of ${\bar{x}}$ such that $U_{\lambda}\subset U$ and
that for every $x_{1}\in U_{\lambda}$, the
sequence $(x_{k})_{k\in\N}$ given by $x_{k+1}=Tx_{k}$ converges
to $x^*\in U_{\lambda}$ such that $-\nabla g(x^*)\in\partial_{p}^{\lambda}f(x^*)$
by Lemma~\ref{l:prox:fixset}.
Since $f+g$ is variationally convex at $\barx$ for $0\in\partial(f+g)(\barx)$,
by Theorem~\ref{thm:coco env grad}, using a smaller $\bar{\lambda}$ if necessary
we can also make
$e_{\lambda}(f+g)$ convex on $U_{\lambda}$.
Now for $x^*\in U_{\lambda}$, we have that $f+g$ is prox-regular at $x^*$ for
$0\in\partial (f+g)(x^*)$ and that $e_{\lambda}(f+g)$ is convex on $U_{\lambda}$.
Using Theorem~\ref{thm:coco env grad}\ref{thm:coco env grad:: v cvx}$
\Leftrightarrow$\ref{envel:convex spec},
we deduce that $f+g$ is variationally convex
at $x^*$ for $0\in\partial (f+g)(x^*)$. Then $x^*$ is a local minimizer for $f+g$
by Fact~\ref{f:localmin:v}.
\end{proof}

As a concrete example, we consider:

 \begin{example}[regularized least squares via variational convexity]
  Consider the best subset selection problem
            $$
                \min\left\lbrace
                    \Vert x\Vert_0|\ Ax = b
                \right\rbrace,
            $$
            where $\Vert x\Vert_0$ is the counting norm,
            $A \in \R^{m\times n}$, and $b \in \R^m$.
            Choose $\gamma>0$ and reformulate the problem as a regularized least squares:
            \begin{equation}\label{e:bestsub:reg}
                \min\left\lbrace
                   f(x) + g(x)|\ x\in\R^n
                \right\rbrace,
            \end{equation}
            where $f(x):= \gamma \Vert x\Vert_0, g(x) := \frac{1}{2}\Vert Ax - b\Vert^2$.
            Since
            $$
               (\forall u\in\R)\ P_{\lambda}\gamma |\cdot|_0 (u) =
                \begin{cases}
                    u, & \text{ if }  |u| \ge \sqrt{2 \gamma \lambda};
                    \\
                    \{0, u\}, & \text{ if } |u| = \sqrt{2 \gamma \lambda};
                    \\
                    0 & \text{otherwise},
                \end{cases}
            $$
            we have
            $
              (\forall {\bf u} \in\R^n)\  P_{\lambda} \gamma \Vert \cdot\Vert_0({\bf u})
                =
                (P_{\lambda}\gamma |\cdot|_0(u_1), \cdots, P_{\lambda}\gamma|\cdot|_0(u_n)).
            $
           The proximal gradient method reads as follows
           \begin{equation}\label{e:hard:thresh}
               (\forall k\in\NN)\  x_{k + 1}\in Tx_{k}:=P_{\lambda}\gamma \Vert \cdot\Vert_0
                (x_{k}- \lambda (A^{\intercal} Ax_{k}- A^{\intercal} b)),
           \end{equation}
            where $0 < \lambda < 1/\|A\|_F^2$, square of the Frobenius norm of $A$.
            %The choice of $\lambda$
%            makes the mapping $\Id-\lambda(A^{\intercal}A-A^{\intercal}b)$ firmly nonexpansive.
            Because that $f$ is variationally convex, see, e.g.,
            \cite[Example 2.5]{Khanh2023variational}, and
            that $g$ is convex and has a Lipschitz continuous gradient, \eqref{e:bestsub:reg}
            fits the framework of Corollary~\ref{c:proxreg:neigh} perfectly!
            %by Theorem~\ref{thm:coco env grad}
%            we have that
%            for $\lambda$ sufficiently small,
%            $P_{\lambda}\gamma\Vert \cdot\Vert_0$ is locally firmly nonexpansive.
%            Hence for $\lambda$ sufficiently small $T$ is an averaged mapping locally.
            It follows from Corollary~\ref{c:proxreg:neigh} that
            the sequence $(x_k)_{k \in \NN}$ given by \eqref{e:hard:thresh}
             converges to a local minimizer $x^*$
            of \eqref{e:bestsub:reg}
            %$$
%                \gamma\Vert x\Vert_0 + \frac{1}{2} \Vert Ax - b\Vert^2,
%            $$
            if $x_{1}$ is sufficiently nearby the solution.
            %Since $f$ is variationally convex,
%             Lemma~\ref{l:c1sum} assures that $f+g$ is variationally convex
%            at $x^*$ for $0\in \partial(f+g)(x^*)$. Therefore, $x^*$ is in fact a local minimizer of $f+g$.
        \end{example}

\begin{remark} Observe that \eqref{e:bestsub:reg} has been studied by Attouch, Bolte \& Svaiter in
\cite[Example 5.4(a)]{attouch2013convergence}
via a different approach, namely $f+g$ is a KL function. Under the assumption that the sequence
$(x_{k})_{k\in\NN}$ is bounded, in \cite[Example 5.4(a)]{attouch2013convergence}
the authors obtain that $(x_{k})_{k\in\NN}$ converges to
a critical point. Because $f+g$ is variationally convex, using variational sufficiency
we conclude that the critical point is
indeed a local minimizer! However, $(x_{k})_{k\in\NN}$ being bounded is indispensable in their proof.
\end{remark}

Finally, we provide an application of Theorem~\ref{thm:coco env grad}, where the celebrated Krasnosel'ski\u{i}-Mann iteration \cite[Theorem 5.15]{BC} kicks in locally on a variationally convex (not necessarily convex) problem.

\begin{theorem}[Krasnosel'ski\u{i}-Mann algorithm under variational convexity]\label{cor:application}
	Let \(f:\Rn\to\overline{\R}\) be proper, lsc and prox-bounded with threshold \(\lambda_f>0\).
	Suppose that \(f\) is variationally convex at \(\bar x\) for \(0\in\partial f(\bar x)\) and that \(f\) is prox-regular around \(\bar x\). Let \(\lambda>0\) and \(U_\lambda\) be those described in Theorem \ref{thm:coco env grad}.
	Take \(x_1\) sufficiently close to \(\bar x\), that is, \(x_1\) belonging to a closed ball \(C\subset
\inte U_\lambda\) centered at \(\bar x\) on which \(f\) is prox-regular.
	Define
	\[
	(\forall k\in\NN)\ x_{k+1} :=\mu_k P_\lambda f(x_k)+(1-\mu_k)x_k,
	\]
	where
	\((\mu_k)_{k\in\N}\) is a sequence in $[0,1]$
	and
	\(
	\sum_{k=1}^\infty \mu_k(1-\mu_k)=\infty
	\).
	Then \(\Fix P_\lambda f\cap C\neq\varnothing\) and the following hold:
	\begin{enumerate}
	\item\label{cor:application::fejer}
	The sequence \((x_k)_{k\in\N}\) stays in \(C\) and is Fej\'er monotone with respect to \(\Fix P_\lambda f\cap C\);
\item $(P_\lambda f(x_k)-x_k)_{k\in\N}$ converges to $0$;
\item\label{cor:application::fejer1}
 $(x_{k})_{k\in\N}$ converges to a point $x^*\in\Fix P_{\lambda}f\cap C$;
	\item\label{cor:application::convergence}
	\(x^*\) is a local minimizer of \(f\).
	\end{enumerate}
\end{theorem}
\proof
By Theorem~\ref{thm:coco env grad}, the mapping
$P_{\lambda}f$ is nonexpansive on $C$ and $P_{\lambda}f(\bar x)={\bar x}$, so
%The assumption on variational convexity implies that \(f\) is
%prox-regular at \(\bar x\) for \(0\in\partial f(\bar x)\),
%thus \(\bar x\in P_\lambda f(\bar x)\) and moreover
\(\bar x\in\Fix P_\lambda f\cap C\). Then $P_{\lambda}f: C\rightarrow C$
and $(x_{k})_{k\in\N}$ lie in the closed convex set $C$.

\ref{cor:application::fejer}--\ref{cor:application::fejer1}:
Apply \cite[Theorem 5.15]{BC} on $C$.
% They follow from \cite[Theorem 5.15]{BC}.

%%	Invoke Theorem \ref{thm:coco env grad} to see that \(P_\lambda f\) is nonexpansive on \(C\), in which case
%It is routine to verify the claim.
%	Here we provide full detail following \cite[Theorem 5.14]{BC} for the sake of completeness.
%	Pick \(x\in\Fix P_\lambda f\cap C\).
%	Then for every \(k\in\N\),
%	\begin{align*}
%		\norm{x_{k+1}-x}^2
%	&=
%		\norm{\mu_k (P_\lambda f(x_k)-P_\lambda f(x))+(1-\mu_k)(x_k-x)}^2\\
%	&=
%		\mu_k\norm{P_\lambda f(x_k)-P_\lambda f(x)}^2+(1-\mu_k)\norm{x_k-x}^2-\mu_k(1-\mu_k)\norm{P_\lambda f(x_k)-x_k}^2\\
%	&\leq
%		\norm{x_k-x}^2-\mu_k(1-\mu_k)\norm{P_\lambda f(x_k)-x_k}^2,
%	\end{align*}
%	which justifies the Fej\'er monotonicity and entails
%	\(
%		\sum_{k=1}^\infty\mu_k(1-\mu_k)\norm{P_\lambda f(x_k)-x_k}^2
%	\leq
%		\norm{x_1-x}^2
%	<
%		\infty
%	\).
%	So we must have \(\liminf_{k\to\infty}\norm{P_\lambda f(x_k)-x_k}=0\), otherwise the assumption on  \(\mu_k\) would enforce \(\sum_{k=1}^\infty\mu_k(1-\mu_k)\norm{P_\lambda f(x_k)-x_k}^2=\infty\).
%	Nevertheless \(P_\lambda f\) is nonexpansive on \(C\), which contains all \((x_k)_{k\in\N}\), so
%	\begin{align*}
%		\norm{P_\lambda f(x_{k+1})-x_{k+1}}
%	&=
%		\norm{P_\lambda f(x_{k+1})-P_\lambda f(x_k)+(1-\mu_k)(P_\lambda f(x_{k})-x_k)}\\
%	&\leq
%		\norm{x_{k+1}-x_k}+(1-\mu_k)\norm{P_\lambda f(x_{k})-x_k}
%	=
%		\norm{P_\lambda f(x_{k})-x_k},
%	\end{align*}
%	implying that \(\norm{P_\lambda f(x_{k})-x_k}\) is in fact convergent.
%	Thereby \(P_\lambda f(x_{k})-x_k\to0\).

\ref{cor:application::convergence}: \ref{cor:application::fejer1} gives
%tells us that
%%	The behavior of \((x_k)_{k\in\N}\) established in \ref{cor:application::fejer} ensures that all cluster points of \((x_k)_{k\in\N}\) lie in \(\Fix P_\lambda f\cap C\), entailing that
$x^*\in\Fix P_\lambda f\cap C$,  hence \(0\in\partial_{p}^{\lambda}f(x^*)\subseteq \partial f(x^*)\)
by Proposition~\ref{thm:optimality}\ref{i:fixed}.
Since \(f\) is prox-regular on $C$,
 \(f\) is
prox-regular at \(x^*\) for \(0\in\partial f(x^*)\).	
	Then Theorem \ref{thm:coco env grad} is applicable at \(x^*\), which in turn
the convexity of \(e_\lambda f\) on \(U_\lambda\) yields that \(f\) is variationally convex at \(x^*\)
for \(0\in\partial f(x^*)\). We can therefore apply
Fact~\ref{f:localmin:v} to conclude.
%entails \(x^*\) to be a local minimizer of \(f\).
\endproof

\begin{remark} Proximal gradient method and others for nonconvex optimization problems
involving KL analytic features have been thoroughly investigated by Attouch, Bolte, and Svaiter in
their seminal works \cite{bolte09,bolte2014proximal}. They have established a fundamental and rich
 convergence theory,
even though the proximal mapping is set-valued. Our work here concentrates on variationally convex
optimization problems, not necessarily having KL analytic features. As shown in
\cite[Section 5, Corollary 4]{edouard}, a differentiable convex function might not have the
KL analytic feature.
\end{remark}

Of course, our results rely on how rich the class of variationally convex functions is. To emphasize
the potential wide applicability of our algorithms, we end this section with a separable
sum rule on variationally convex functions and some examples.

\begin{proposition}\label{p:sep:rule} Suppose that each lsc function $f_{i}:\R^{n_{i}}\rightarrow\OR$ is
 is variationally convex at $\bar{x}_{i}$
for $\bar{v}_{i}\in\partial f_{i}(\bar{x}_{i})$, where $n_{i}\in\NN$ and $i=1,\ldots, m$. Then the separable sum
$f:=f_{1}\oplus\cdots \oplus f_{m}$ is variationally convex at $\bar{x}:=(\bar{x}_{1},\ldots,\bar{x}_{m})$
for $\bar{v}:=(\bar{v}_{1},\ldots, \bar{v}_{m})\in\partial f(\bar{x}).$
\end{proposition}

\begin{proof} By \cite[Theorem 1]{rock-vietnam}, the assumption implies that there exist $U_{i}$ a neighborhood
of $\bar{x}_{i}$, $V_{i}$ a neighborhood of $\bar{v}_{i}$, and $\varepsilon>0$ such that
\begin{equation}\label{e:ifunction}
\gph\partial f_{i}\cap (U_{i,\varepsilon}\times V_{i}) \text {  is  monotone},
\end{equation}
where $U_{i,\varepsilon}:=\{x_{i}\in U_{i}|\ f_{i}(x_{i})<f_{i}(\bar{x}_{i})+\varepsilon\}$.
Since each $f_{i}$ is lsc at $\bar{x}_{i}$, we can choose $\delta>0$ sufficiently small such that
$(\forall i=1,\ldots, m)\ \Ball_{\delta}(\bar{x}_{i})\subset U_{i}$ and
\begin{equation}\label{e:lscbound}
(\forall x_{i}\in \Ball_{\delta}(\bar{x}_{i}))\ f_{i}(x_{i})>f_{i}(\bar{x}_{i})-\frac{\varepsilon}{2m}.
\end{equation}
Let
\begin{align}
U:&= \Ball_{\delta}(\bar{x}_{1})\times\cdots\times \Ball_{\delta}(\bar{x}_{m}),\\
V:&= V_{1}\times \cdots\times V_{m},\\
U_{\varepsilon}:&=\{x=(x_{1},\ldots, x_{m})\in U|\ f_{1}(x_{1})+\cdots+f_{m}(x_{m})<f_{1}(\bar{x}_{1})+\cdots+f_{m}(\bar{x}_{m})+\varepsilon/2\}.
\end{align}
Then $U\times V$ is a neighborhood of $(\bar{x},\bar{v})$. Observe that
$\partial f(x)=\partial f_{1}(x_{1})\times \cdots\times \partial f_{m}(x_{m})$.
For $(x,v)\in\gph\partial f\cap (U_{\varepsilon}\times V)$,
we have
\begin{equation}\label{e:setsun}
\text{$x_{i}\in \Ball_{\delta}(\bar{x}_{i})\subseteq U_{i}$, $v_{i}\in V_{i}$, $(x_{i},v_{i})\in\gph\partial f_{i}$,}
\end{equation}
 and
\begin{align}
f_{i}(x_{i})& <f_i(\bar{x}_{i})+\sum_{j\neq i}(f_{j}(\bar{x}_{j})-f_{j}(x_{j}))+\frac{\varepsilon}{2}\\
&< f_i(\bar{x}_{i})+(m-1)\frac{\varepsilon}{2m}+\frac{\varepsilon}{2}=f_i(\bar{x}_{i})+\varepsilon\label{e:functionsun}
\end{align}
by \eqref{e:lscbound}. Thus,
\eqref{e:setsun}--\eqref{e:functionsun} together mean
that $x_{i}\in U_{i,\varepsilon}$, so that $(x_{i},v_{i})\in\gph\partial f_{i}\cap
(U_{i,\varepsilon}\times V_{i})$
for $i=1,\ldots, m$.
Using $(x,v)=((x_{1},\ldots, x_{m}), (v_{1},\ldots, v_{m}))$ and \eqref{e:ifunction}, we deduce that
$\gph \partial f\cap (U_{\varepsilon}\times V)$ is monotone.
The proof is complete by using \cite[Theorem 1]{rock-vietnam} again.
\end{proof}

\begin{example}\label{e:thurs1} Let $\varphi_i:\R\rightarrow\R$ be piecewise constants and lower semicontiunous continuous with $i=1,\ldots, n$.
Then $\Phi:\Rn\rightarrow\R: x\mapsto \sum_{i=1}^n\varphi_{i}(x_{i})$ is varitationally convex everywhere.
\end{example}

\begin{example}\label{e:thurs2}
Let $\varphi_{i}:\R\rightarrow\overline{\R}$ be piecewise affine, lower semicontiunous continuous, and
$(\forall t\in\dom\varphi_{i})\ \hat{\partial}\varphi_i(t)\neq\varnothing$ with $i=1,\ldots, n$.
Then $\partial \varphi_{i}=\hat{\partial}\varphi_{i}$, and $\varphi_{i}$ is variationally convex at every $t\in\dom\varphi_i$ for $s\in\hat{\partial}\varphi_{i}(t)$. Consequently,
$\Phi:\Rn\rightarrow\overline{\R}: x\mapsto \sum_{i=1}^n\varphi_{i}(x_{i})$ is
variationally convex at $x\in\dom\Phi$ for $v\in\hat{\partial}\Phi(x)$.
\end{example}

\begin{example}\label{e:thurs3}
 Let $0<p_{i}< 1$ for $i=1,\ldots, n$. The function
$\Phi:\Rn\rightarrow\R: x\mapsto \sum_{i=1}^n |x_{i}|^{p_{i}}$ is variationally convex at $0$
for every $v\in\partial \Phi(0)$. Note that for
$\varphi_{i}:\R\rightarrow\R: t\mapsto |t|^{p_{i}}$ one has
$$\hat{\partial}\varphi_{i}(t)=\partial\varphi_{i}(t)=\begin{cases}
p_{i}t^{p_{i}-1}, & \text{ if $t>0$;}\\
\R, & \text{ if $t=0$;}\\
-p_{i}(-t)^{p_{i}-1}, & \text{ if $t<0$,}
\end{cases}
$$
thus $\hat{\partial}\Phi(x)=\partial\Phi(x)=\partial\varphi_{1}(x_{1})\times\cdots\times\partial\varphi_{n}(x_{n})$.
For applications of $\Phi$, see \cite{candes08,attouch2013convergence}.
\end{example}

\begin{example}\label{e:thurs4} Let $\varepsilon>0$. The function $\Phi:\Rn\rightarrow\R:x\mapsto\sum_{i=1}^{n}\ln(|x_{i}|+\varepsilon)$
is variationally convex at $0$ for $v\in\inte\partial \Phi(0)=(-1/\varepsilon, 1/\varepsilon)^n$.
When $\varepsilon=1$, the result is given in \cite[Example 2.6]{Khanh2023variational}.
For application of $\Phi$, see \cite{candes08}.
\end{example}

The effectiveness of proximal gradient algorithms in Theorem~\ref{t:f:plusg} relies on the computation of proximal mappings of variationally convex functions. While in general the computation of proximal mappings might not be easy,
for functions of separable nature given in Proposition~\ref{p:sep:rule}, this reduces to one dimensional
optimization problems, see, e.g., \cite[Section 6.3]{beck2017}.

\section{Locally single-valued proximal mappings}\label{s:locally}

For a proper, lsc and prox-bounded function $f:\Rn\rightarrow\OR$ with threshold $\lambda_{f}>0$,
when $0<\lambda<\lambda_{f}$ the proximal mapping $P_{\lambda}f$ is always nonempty, compact-valued, and upper
semicontinuous, see \cite[Theorem 1.25]{rockafellar_variational_1998}. If $P_{\lambda}f(x)$ is a singleton,
then $P_{\lambda}f$ is continuous at $x$.
Without variational convexity of $f$, the single-valuedness of proximal mapping
$P_{\lambda}f$ at a point $x$ is
equivalent to the strict differentiability of $e_{\lambda}f$ at $x$, as we show in this section.
On an open set, this result extends \cite[Theorem 3.5]{wang2010Chebyshev}:
	``For a proper, lsc and prox-bounded function $f:\Rn\to\overline{\R}$,
	$P_\lambda f$
is single-valued if and only if $f+\lambda^{-1}j$ is essentially strictly convex."
Our proof hinges on the pointwise duality.
 %we show in this section that $P_\lambda f$ is single-valued at a point $x$ if and only if $e_{\lambda}f$ is strictly differentiable at $x$.
%}
Denote by $\Gamma_0$ the class of functions $\psi:\mathbb{R}_+\to
\overline{\mathbb{R}}_+ :=[0,\infty]$ that are proper, lsc, and convex,
and satisfy $\psi(t)=0$ if and only if $t=0$.
The following is a pointwise variation of \cite[Defnition 2(i)]{Stromberg2011}; see also \cite[Corollary 3.4.4(i)]{zalinescu2002convex}.

\begin{defn}
	Let $g:\mathbb{R}^n\to \overline{\mathbb{R}}$ be proper and lsc.
	Let $u\in\dom g$ and $x\in \mathbb{R}^n$.
	We say $g$ is essentially strongly convex  at $u$ for $x$, if there exists a modulus $\psi\in\Gamma_0$ such that
	$$
	(\forall y\in \mathbb{R}^n)~
	g(y)\geq g(u)+\left\langle y-u,x \right\rangle+\psi \left(\left\|y-u\right\|\right).
	$$
\end{defn}
Two lemmas concerning the differentiability of Fenchel conjugate are convenient
 for the proof of main result in this section.
\begin{lemma}\label{lem:pointwise dual of differentiability}
	\emph{\cite[Corollary 3.4.4(i)\&(vii)]{zalinescu2002convex}}
	Let $g:\Rn\to\bR$ be a proper, lsc and convex function,
and let $(x,u)\in \operatorname{gph}\partial g^*$.
	Then $x\in \operatorname{int}\dom g^*$ and $g^*$ is differentiable at $x$ if and only if $g$ is essentially strongly convex  at $u$ for $x$.
\end{lemma}

%Below is a finite-dimensional version of a result by Fabian-Zizler \cite{fabian1999}, which will be useful soon.
\begin{lemma}\label{lem: differentiable conjugate}
	\emph{\cite[Lemma 3]{fabian1999}}
	Let $f:\Rn\to\overline{\R}$ be proper and lsc with $\inf f>-\infty$, and suppose that the conjugate $f^*$ is differentiable at $x\in\Rn$.
	Let $u :=\nabla f^*(x)$.
	Then
	$f^{**}(u)=f(u)$.
\end{lemma}

We are now ready to present the main result of this section.
%It turns out that single-valued $P_\lambda f(x)$

\begin{theorem}[proximal mapping single-valued at a point] \label{thm:equi of singleval prox}
	Let $f:\Rn\to\overline{\R}$ be proper, lsc and prox-bounded with
threshold $\lambda_f>0$. Let $ 0<\lambda<\lambda_f$ and let $x\in\dom f$.
	%Write $F=\lambda f+j$ for simplicity.
	Then the following are equivalent:
	\begin{enumerate}
		\item  \label{thm:equi of singleval prox:prox}
		$P_\lambda f(x)$ is a singleton; %\red{(TODO: add $u$)}
		\item \label{thm:equi of singleval prox:envelope}
		$e_\lambda f$ is strictly differentiable at $x$;%  \red{(TODO: add $u$)}
		\item \label{thm:equi of singleval prox:conjugate}
		$(\lambda f+j)^*$ is (strictly) differentiable at $x$;
%with gradient $u=\nabla (\lambda f+j)^*(x)$.
		\item  \label{thm:equi of singleval prox:cvs hull} $(\exists u\in P_\lambda f(x))$ $\lambda f+j$ is essentially strongly convex  at $u$ for $x$.
%\item  \label{i:convex:strong} $(\exists u\in P_\lambda f(x))$ $\conv(\lambda f+j)$ is essentially strongly convex  at $u$ for $x$, and $\conv (\lambda f+j)(u)=(\lambda f+j)(u)$.
%\in\lambda \partial_p^\lambda f(u)+u$ with $(\lambda f+j)(u)=\conv (\lambda f+j)(u)$.
		%	 \item \label{thm:equi of singleval prox:level prox}	 $(\exists x^*\in\Rn)$ $\partial_p^\lambda f(x^*)\neq\emptyset$.
	\end{enumerate}
	%Then $\ref{thm:equi of singleval prox:prox}\Leftrightarrow \ref{thm:equi of singleval prox:envelope}\Leftrightarrow \ref{thm:equi of singleval prox:conjugate}\Leftrightarrow\ref{thm:equi of singleval prox:cvs hull}$.
	%Assume in addition that $u\in\inte \dom (\lambda f+j)^*$. Then $\ref{thm:equi of singleval prox:conjugate}\Leftrightarrow\ref{thm:equi of singleval prox:cvs hull}$.
When one of the above holds, $P_\lambda f(x)=\nabla(\lambda f+j)^*(x)$.
\end{theorem}
\proof
We shall prove the equivalence \ref{thm:equi of singleval prox:prox}$\Leftrightarrow$ \ref{thm:equi of singleval prox:envelope}$\Leftrightarrow $\ref{thm:equi of singleval prox:conjugate} first, then justify \ref{thm:equi of singleval prox:conjugate}$\Leftrightarrow$\ref{thm:equi of singleval prox:cvs hull}.

``\ref{thm:equi of singleval prox:prox}$\Rightarrow$\ref{thm:equi of singleval prox:envelope}":
Note that $\partial(e_\lambda f)(x)\subseteq\lambda^{-1}[x-P_\lambda f(x)]$ and $e_\lambda f$ is locally Lipschitz; see, e.g.,~\cite[Example 10.32]{rockafellar_variational_1998}.
The assumption implies that $\partial(e_\lambda f)(x)$ is at most single-valued.
% and $e_\lambda f$ is strictly continuous at $x$.
Invoking \cite[Theorem 9.18]{rockafellar_variational_1998}
 yields that $e_\lambda f$ is strictly differentiable at $x$.

``\ref{thm:equi of singleval prox:envelope}$\Rightarrow$\ref{thm:equi of singleval prox:prox}": By the
strict differentiability of $e_{\lambda}f$ at $x$, $\partial(-e_\lambda f)(x)$ is a singleton, then so is $\conv P_\lambda f(x)=\lambda \partial(-e_\lambda f)(x)+x$ by \cite[Example 10.32]{rockafellar_variational_1998} .

``\ref{thm:equi of singleval prox:envelope}$\Leftrightarrow$\ref{thm:equi of singleval prox:conjugate}":
The desired equivalence follows from the identity
\begin{equation}\label{e:envel:conjug}
\lambda^{-1}j-e_\lambda f=
\left(
f+\lambda^{-1}j
\right)^*\left(\frac{\cdot}{\lambda}\right)
\Leftrightarrow
j-\lambda e_\lambda f=\left(\lambda f+j\right)^*,
\end{equation}
where the first equality on the left hand side holds thanks to \cite[Example 11.26]{rockafellar_variational_1998}, and the displayed equivalence is a consequence of Fenchel conjugate calculus; see, e.g.,~\cite[Proposition 13.23]{BC}.

%Assume now $u\in\inte \dom (\lambda f+j)^*$.
To prove \ref{thm:equi of singleval prox:conjugate}$\Leftrightarrow$\ref{thm:equi of singleval prox:cvs hull},
let us put $g :=\conv (\lambda f+j)=\overline{\conv }(\lambda f+j)$, where the last equality holds by Lemma~\ref{lem: proximal and closed cvx hull}.
Then $g$ is proper, lsc and convex, and $g^*=(\lambda f+j)^*$.
%to which  Lemma~\ref{lem:pointwise dual of differentiability} is applicable.

``\ref{thm:equi of singleval prox:conjugate}$\Rightarrow$\ref{thm:equi of singleval prox:cvs hull}":
By assumption, $g^*$ is differentiable at $x$. With $u :=\nabla g^*(x)$, we have $u\in\partial g^*(x)\Leftrightarrow x\in\partial g(u)$.
Then Lemma~\ref{lem:pointwise dual of differentiability} yields that $g$ is essentially strongly convex  at $u$ for $x$.
In turn,
\begin{equation}\label{eq: prox thm}
	(\forall y\in \mathbb{R}^n)~(\lambda f+j)(y)\geq g(y)\geq g(u)+\left\langle y-u,x \right\rangle+\psi(\left\|y-u\right\|),
	%	&=(\lambda f+j)(u)+\left\langle y-u,x \right\rangle+\psi(\left\|y-u\right\|).
\end{equation}
where $\psi\in\Gamma_{0}$.
We claim that $g(u)=(\lambda f+j)(u)$. If so, Fact~\ref{fact: level and Fenchel sub} and
\eqref{eq: prox thm} imply that
$$
x\in\partial\conv (\lambda f+j)(u)
=\partial_F(\lambda f+j)(u)
=\lambda\partial_F(f+\lambda ^{-1}j)(u)
=\lambda\partial_p^\lambda f(u)+u,
$$
i.e., $u\in P_{\lambda}f(x)$;
and by (\ref{eq: prox thm}) again
\(
(\forall y\in \mathbb{R}^n)~(\lambda f+j)(y)\geq(\lambda f+j)(u)+\left\langle y-u,x \right\rangle+\psi(\left\|y-u\right\|),
\)
completing the proof.
It remains to justify  $g(u)=(\lambda f+j)(u)$.
To see this, invoking the identity $(\lambda f+j)^*=g^*$, Lemma~\ref{lem: differentiable conjugate} entails that
$
( \lambda f+j)(u)=(\lambda f+j)^{**}(u)=\conv (\lambda f+j)(u)=g(u).
$

``\ref{thm:equi of singleval prox:cvs hull}$\Rightarrow$ \ref{thm:equi of singleval prox:conjugate}":
By assumption, $\exists \psi\in\Gamma_{0}$ such that
$
(\forall y\in \mathbb{R}^n)~
(\lambda f+j)(y)\geq (\lambda f+j)(u)+\left\langle y-u,x\right\rangle+\psi(\left\|y-u\right\|),
$
which gives
\begin{equation}\label{e:convex:ineq}
(\forall y\in \mathbb{R}^n)~
g(y)=\conv (\lambda f+j)(y)\geq (\lambda f+j)(u)+\left\langle y-u,x\right\rangle+\psi(\left\|y-u\right\|),
\end{equation}
from which $g(u)\geq (\lambda f+j)(u)$. Since $(\lambda f+j)(u)\geq g(u)$ always holds, we have
$g(u)=(\lambda f+j)(u)$.
Thus, \eqref{e:convex:ineq} can be written as
$$(\forall y\in \mathbb{R}^n)~
g(y)\geq g(u)+\left\langle y-u,x\right\rangle+\psi(\left\|y-u\right\|),$$
in particular, $x\in\partial g(u)$.
%Now, note that $(\lambda f+j)(u)=g(u)$. Then $\partial_F(\lambda f+j)(u)=\partial g(u)$ by Lemma \ref{lem: Fenchel of conv} and appealing to Fact~\ref{fact: level and Fenchel sub} yields
%$$
%x\in\lambda\partial_p^\lambda f(u)+u=\lambda\left[\partial_F(f+\lambda^{-1}j)(u)-\lambda^{-1}(u) \right]+u=\partial_F (\lambda f+j)(u)=\partial g(u).
%$$
%In particular $u\in P_\lambda f(x)$ by Fact~\ref{fact: prox identity}.
Then $g$ is essentially strongly convex  at $u$ for $x\in\partial g(u)$, so that $g^*$ is differentiable at $x$ by Lemma~\ref{lem:pointwise dual of differentiability}. Since $g^*=(\lambda f+j)^*$,
\ref{thm:equi of singleval prox:conjugate} is established.

Finally, when one of \ref{thm:equi of singleval prox:prox}--\ref{thm:equi of singleval prox:cvs hull} holds,
by the identity \eqref{e:envel:conjug} and
$\partial [-e_{\lambda}f](x)=\lambda^{-1}[\conv P_{\lambda}f(x)-x],$ we obtain
$\conv P_{\lambda}f(x)=\nabla (\lambda f+j)^*(x)$, i.e.,
$P_{\lambda}f(x)=\nabla (\lambda f+j)^*(x)$.
\endproof

Immediately, one has the following result.
\begin{corollary}\label{c:prox-sing}
	Let $f:\Rn\to\overline{\R}$ be proper, lsc and prox-bounded with threshold $\lambda_f>0$. Let $U\subseteq \Rn$ and $ 0<\lambda<\lambda_f$.
	%Write $F=\lambda f+j$ for simplicity.
	Then the following are equivalent:
	\begin{enumerate}
		\item
		$P_\lambda f$ is single-valued (so continuous) on $U$;
		\item
		$e_\lambda f$ is strictly differentiable on $U$;
		\item
		$(\lambda f+j)^*$ is (strictly) differentiable on $U$;
		\item $(\forall x\in U)(\exists u\in P_\lambda f(x))$ $\lambda f+j$ is essentially strongly convex  at $u$ for $x$.
	\end{enumerate}
When one of the above holds, $\ P_\lambda f=\nabla(\lambda f+j)^*$ on $U$.
\end{corollary}
%Equipped with Theorem \ref{thm:equi of singleval prox},

\begin{remark}
It is worthwhile to compare Corollary~\ref{c:prox-sing} (single-valued $P_{\lambda}f$
related to the differentiability of $(f+\lambda^{-1} j)^*$) to
Corollary~\ref{c:levelsub-sing} (single-valued $\partial_{p}^{\lambda}f$
related to the differentibility of $\conv(f+\lambda^{-1} j)$).
\end{remark}
\begin{remark} Readers might ask what is the connection between results in section~\ref{s:locally}
and results in section~\ref{s:variat}? The picture is clear.
Corollary~\ref{c:prox-sing} says that $P_{\lambda}f$
being single-valued (of course continuous) on a neighborhood corresponds to $e_{\lambda}f$ being
strictly differentiable; while Theorems~\ref{thm:coco env grad} and \ref{thm:variational str cvx}
say that for a prox-regular function $P_{\lambda}f$ being
Lipschitz continuous with various moduli on a neighborhood corresponds to $e_{\lambda}f$ being
convex with appropriate moduli. These intrinsic and novel properties are certainly worth of pointing out.
\end{remark}

When $U=\Rn$, we obtain the following result which improves \cite[Theorem 3.5]{wang2010Chebyshev}.
\begin{corollary}\label{cor:prox singlevalued}
Let $f:\mathbb{R}^n\to\overline{\mathbb R}$ be proper, lsc and prox-bounded with threshold $\lambda_f>0$.
Let $0<\lambda<\lambda_f$.
Then the following are equivalent:
\begin{enumerate}
	\item\label{cor:prox singlevalued:a} $P_\lambda f$ is single-valued everywhere;
	\item\label{cor:prox singlevalued:b} $(\forall x\in \mathbb{R}^n)(\exists u\in P_\lambda f(x))$ $\lambda f+j$ is essentially strongly convex at $u$ for $x$;
%\in\lambda \partial_p^\lambda f(u)+u$ and $(\lambda f+j)(u)=\operatorname{conv}(\lambda f+j)(u)$;
	%In particular $u\in P_\lambda f(x)$ by Theorem~\ref{fact: prox identity};
	\item\label{cor:prox singlevalued:c} $\lambda f+j$ is essentially strictly convex.
\end{enumerate}
When one of the above holds, we have $P_{\lambda}f=\nabla (\lambda f+j)^*$.
\end{corollary}
\proof
``\ref{cor:prox singlevalued:a}$\Rightarrow$\ref{cor:prox singlevalued:b}": Apply Theorem~\ref{thm:equi of singleval prox}.

``\ref{cor:prox singlevalued:b}$\Rightarrow$\ref{cor:prox singlevalued:c}": For convenience,
set $g :=\lambda f+j$.
The assumption asserts that to each $x$ corresponds some $u$
and modulus $\psi_u\in\Gamma_{0}$ such that
$
(\forall y\in\mathbb{R}^n)\
g (y)\geq g (u)+\left\langle y-u,x \right\rangle+\psi_u\left(\left\|y-u\right\|\right),
$
further implying that
$
(\forall y\in\mathbb{R}^n)\
(\operatorname{conv}g) (y)\geq
g(u)+\left\langle y-u,x \right\rangle+\psi_u\left(\left\|y-u\right\|\right).
$
Since $g(u)\geq (\conv g)(u)$, we have
$$
(\forall y\in\mathbb{R}^n)\
(\operatorname{conv}g) (y)\geq
(\conv g)(u)+\left\langle y-u,x \right\rangle+\psi_u\left(\left\|y-u\right\|\right).
$$
By Lemma~\ref{lem:pointwise dual of differentiability}, we conclude that  $(\conv g)^*$ is differentiable on $\inte\dom (\conv g)^*=\Rn$. Since $g^*=(\conv g)^*$ and $g^*$ is differentiable on $\Rn$,  we derive that $g$ is convex by \cite[Corollary 2.3]{soloviov}.
Because $g^*$ is differentiable, $g$ is essentially strictly convex.

%Note that $\clconv g=\conv g$ by Lemma \ref{lem: proximal and closed cvx hull}.
%Then the function $\conv g=(\clconv g)^{**}=(\conv g)^{**}=\left[(\conv g)^* \right]^{*}$ is essentially strictly convex; see, e.g., \cite[Theorem 5.4]{Bauschke2001}.

%To see that $g$ is essentially strictly convex, it suffices to show $\conv g=g$ on $\dom\partial_F g$ because $\dom\partial_F g\subseteq \dom(\partial\conv g)$.
%Invoking Lemma \ref{lem: differentiable conjugate} entails that for every $x\in\mathbb{R}^n$ and $u=\nabla g ^*(x)$
%$$
%g (u)
%=
%g ^{**}(u)
%=(\clconv g )(u)
%=(\conv g) (u),
%$$
%where the last equality follows from Lemma \ref{lem: proximal and closed cvx hull}.
%Hence $\lambda f+j=\conv g $ on $\ran\nabla g ^*=\operatorname{dom}\partial g $, completing the proof.

``\ref{cor:prox singlevalued:c}$\Rightarrow$\ref{cor:prox singlevalued:a}":
Observe that when $0<\lambda<\lambda_{f}$,
\begin{equation}\label{e:prox:s}
(\forall x\in\Rn)\ \varnothing\neq P_{\lambda}f(x)=\argmin_{y\in\Rn}\bigg(f(y)+\frac{1}{\lambda}j(y)-\frac{1}{\lambda}\langle y, x\rangle \bigg).
\end{equation}
The assumption implies that $y\mapsto f(y)+\frac{1}{\lambda}j(y)-\frac{1}{\lambda}\langle y, x\rangle$ is essentially strictly
convex, so it has at most one minimizer. Hence $P_{\lambda}f(x)$ is single-valued.
%Finally, to find $P_{\lambda}f$, we use \eqref{e:prox:s} to obtain
%$y=P_{\lambda}f(x)$ if and only if
%$0\in \partial(\lambda f+j)(y)-x$, equivalently, $y\in \partial (\lambda f+j)^*(x)=\nabla(\lambda f+j)^*(x).$
\endproof

We refer readers to \cite{jourani2014,thibault} for more details on the differentiability of Moreau envelopes.

\section{Integration of level proximal subdifferentials}\label{s:integ}
Given the fact that the domain of level proximal subdifferential of a function might be small,  can it
determine the function uniquely up to a constant?
This section concerns consequences of proper, lsc and prox-bounded functions $f_1$ and $f_2$ having the same level proximal subdifferential. We show that level proximal subdifferential can determine
the proximal hull of the function uniquely up to a constant and
that Rockafellar's technique
of integrating cyclically monotone operators, see,
e.g., \cite[Theorem 12.25]{rockafellar_variational_1998},
can be applied to find the proximal hull.

%\begin{lemma}\label{lem:prox coincidence}
%Let $f_i:\Rn\to\overline{\R}$ be proper, lsc, and prox-bounded with threshold $\lambda_{f_i}>0$
%for $i=1, 2$, and let $0<\lambda<\min\{\lambda_{f_1},\lambda_{f_2}\}$.
%\begin{enumerate}
%\item
%If $P_\lambda f_1=P_\lambda f_2$, then
%$
%(\exists c\in\mathbb R)\ \clconv f_1=\clconv f_2+c,
%$
%%for some constant $c\in\R$,
%provided that $\conv f_i$ is proper for $i=1,2$.
%\item If $e_{\lambda}f_{1}=e_{\lambda}f_{2}$, then
%$\clconv f_1=\clconv f_2,
%$
%%for some constant $c\in\R$,
%provided that $\conv f_i$ is proper for $i=1,2$.
%\end{enumerate}
%\end{lemma}
%\proof
%\endproof

\begin{theorem}\label{thm:coincidence}
Let $f_i:\Rn\to\overline{\mathbb R}$ be proper, lsc and prox-bounded with threshold
$\lambda_{f_i}>0$ for $i=1, 2$, and let
 $0<\lambda<\min\{\lambda_{f_1},\lambda_{f_2}\}$.
Consider the following statements:
\begin{enumerate}
\item\label{thm:coincidence:level sub} $\partial_p^\lambda f_1=\partial_p^\lambda f_2$;
\item\label{thm:coincidence:prox} $P_\lambda f_1=P_\lambda f_2$;
\item\label{thm:coincidence:envlope} $(\exists c\in\mathbb R)$ $e_\lambda f_1=e_\lambda f_2+c$;
\item\label{thm:coincidence:proximal hull} $(\exists c\in\mathbb R)$ $h_\lambda f_1=h_\lambda f_2+c$;
\item\label{thm:coincidence:cvx hull} $(\exists c\in\mathbb R)$ $\overline{\conv }f_1=\overline{\conv }f_2+c$, provided that $\conv f_i$ is proper for $i=1, 2$.
\end{enumerate}
Then $\ref{thm:coincidence:level sub}\Leftrightarrow\ref{thm:coincidence:prox}\Rightarrow\ref{thm:coincidence:envlope}\Leftrightarrow
\ref{thm:coincidence:proximal hull}$,
 and $\ref{thm:coincidence:level sub}\Rightarrow\ref{thm:coincidence:cvx hull}$.
\end{theorem}
\proof
``\ref{thm:coincidence:level sub}$\Leftrightarrow$\ref{thm:coincidence:prox}":
$\partial_p^\lambda f_1=\partial_p^\lambda f_2\Leftrightarrow\text{Id}+\lambda\partial_p^\lambda f_1=\text{Id}+\lambda\partial_p^\lambda f_2\Leftrightarrow P_\lambda f_1=P_\lambda f_2$.

``\ref{thm:coincidence:prox}$\Rightarrow$\ref{thm:coincidence:envlope}":
%Apply Lemma~\ref{lem:prox coincidence}.
This implication can be found in the proof of \cite[Proposition 3.16(ii)]{wang2010Chebyshev}.
Here we include a similar argument for the sake of completeness.
Recall from \cite[Example 10.32]{rockafellar_variational_1998} that $-e_\lambda f_i$ is
lower-$C^2$ for each $i=1,2$. Therefore, $-e_\lambda f_i$ is primal lower-nice in the sense of \cite[Definition 3.1]{poliquin1991} and
$\partial_p(-e_\lambda f_i)=\partial(-e_\lambda f_i)$ by \cite[Proposition 3.5]{poliquin1991} for $i=1,2$.
Furthermore,
\[
(\forall x\in\Rn)~
\partial_p(-e_\lambda f_1)(x)=
%\partial(-e_\lambda f_1)(x)=
\lambda^{-1}\left[\conv P_\lambda f_1(x)-x\right]
=
\lambda^{-1}\left[\conv P_\lambda f_2(x)-x\right]
%=\partial(-e_\lambda f_2)(x)
=
\partial_p(-e_\lambda f_2)(x),
\]
where the first equality follows from \cite[Example 10.32]{rockafellar_variational_1998}, and the second from our assumption.
In turn~\cite[Theorem 4.1]{poliquin1991} implies that there exists
a neighborhood of $x$ on which $e_\lambda f_1-e_\lambda f_2$ is a constant,
and this holds for every $x\in\Rn$. Since both
$e_{\lambda}f_{1}, e_{\lambda}f_{2}$ are locally Lipschitz on $\Rn$, we deduce that
$(\exists c\in\R)\ e_{\lambda}f_{1}=e_{\lambda}f_{2}+c$.

``\ref{thm:coincidence:envlope}$\Rightarrow$\ref{thm:coincidence:proximal hull}":
Applying the identity $h_\lambda f_i=-e_\lambda(-e_\lambda f_i)$, we have
$
h_\lambda f_1=-e_\lambda(-e_\lambda f_1)
=
-e_\lambda(-e_\lambda f_2-c)
=
-[e_\lambda(-e_\lambda f_2)-c]
=
h_\lambda f_2+c.$
``\ref{thm:coincidence:proximal hull}$\Rightarrow$\ref{thm:coincidence:envlope}":
Apply Lemma~\ref{lem:prox of proximal hull}\ref{i:rock:prox}.
%the identity $e_{\lambda}[h_{\lambda} f]=e_{\lambda}f$, see \cite[Example 1.44]{rockafellar_variational_1998}.

``\ref{thm:coincidence:level sub}$\Rightarrow$\ref{thm:coincidence:cvx hull}":
%Apply Lemma~\ref{lem:prox coincidence}.
Apply \cite[Proposition 3.16(ii)]{wang2010Chebyshev}.
\endproof

The example below shows that the coincidence of closed convex hulls up to a constant does not imply that of level proximal subdifferentials.
\begin{example}
	\label{example:coincidience}
Consider the zero ``norm"
\begin{align*}
f:\R\rightarrow\R:
x\mapsto
	\begin{cases}
	1,&\text{ if }x\neq0;\\
	0,&\text{ if }x=0,
	\end{cases}
\end{align*}
and $g :=0$.
Then the following hold:
\begin{enumerate}
	\item\label{example:coincidience:hull} $\clconv f=\clconv g=g=0$;
	\item\label{example:coincidience:subd}$(\forall \lambda>0)$ $\partial_p^\lambda f\neq\partial_p^\lambda g$
because for $x\in\R$,
$$\partial_{p}^{\lambda}f(x)=\begin{cases}
\varnothing, &\text{ if $0<|x|<\sqrt{2\lambda}$;}\\
[-\sqrt{2/\lambda}, \sqrt{2/\lambda}], & \text{ if $x=0$;}\\
0, &\text{ if $|x|\geq \sqrt{2\lambda}$},
\end{cases}
$$
and $\partial_{p}^\lambda g(x)=0.$
\end{enumerate}
\end{example}
\proof
\ref{example:coincidience:hull}:
We claim that $\clconv f=0$.
%To justify the claim, since $\clconv
%f(0)\leq f(0)=0$,
Because $\clconv f\geq 0$, it suffices to show $\clconv f\leq 0$.
Indeed,
\begin{align*}
	(\forall x\in\mathbb{R})~~
	(\clconv f)(x)&\leq
(\conv f)(x)
	\leq \inf_{\substack{0<\lambda<1\\ u,v\in\mathbb{R}}}
	\{
	\lambda f(u)+(1-\lambda )f(v): x=\lambda u+(1-\lambda)v
	\}\\
	&\leq
	\inf_{0<\lambda<1}
	\left(
	\lambda f(\lambda^{-1}x)+(1-\lambda)f(0)
	\right)
	\leq
	\inf_{0<\lambda<1}
	\lambda
	=0.
\end{align*}

\ref{example:coincidience:subd}: Apply \cite[Example 4.2]{wang23level}.
\endproof

The next example shows that Theorem \ref{thm:coincidence}\ref{thm:coincidence:envlope} does not imply Theorem \ref{thm:coincidence}\ref{thm:coincidence:level sub} either.
\begin{example}\label{example:fail coin}
	For $\lambda>0$, define $f,g:\R\rightarrow\R$ by
	\begin{align*}
		f(x)
		:=
		\begin{cases}
			1,&\text{ if }x\neq0;\\
			0,&\text{ if }x=0,
		\end{cases}
		\text{ and }
		g(x)
		:=
		\begin{cases}
			\sqrt{2\lambda^{-1}}|x|-\lambda^{-1}j(x), &\text{ if }|x|\leq\sqrt{2\lambda};\\
			1, &\text{ otherwise }.
		\end{cases}
	\end{align*}
	Then both $f$ and $g$ are proper, lsc and prox-bounded with thresholds $\lambda_f=\lambda_g=\infty$.
	 Moreover, the following hold:
	\begin{enumerate}
	%\item \label{example:fail coin::prox-bd}	Both $f$ and $g$ are prox-bounded with threshold $\lambda_f=\lambda_g=\infty$.
		\item\label{example:fail coin::hull} $e_\lambda f=e_\lambda g$, because of $h_\lambda f=g$;
		\item\label{example:fail coin::prox} $P_\lambda f\neq P_\lambda g$.
		To be specific,
		\begin{align*}
			P_\lambda f(x)=
			\begin{cases}
				\{0\}, &\text{ if }|x|<\sqrt{2\lambda};\\
				\{0,x\}, &\text{ if }|x|=\sqrt{2\lambda};\\
				x,&\text{otherwise},
			\end{cases}
			\text{ and }
			P_\lambda g(x)=
			\begin{cases}
				\{0\}, &\text{ if }|x|<\sqrt{2\lambda};\\
				[0,\sqrt{2\lambda}], &\text{ if }x=\sqrt{2\lambda};\\
				[-\sqrt{2\lambda},0], &\text{ if }x=-\sqrt{2\lambda};\\
				x,&\text{otherwise},
			\end{cases}
		\end{align*}
		whereas
		\begin{align*}
			e_\lambda f(x)=e_\lambda g(x)=
			\begin{cases}
				1, &\text{ if }|x|\geq \sqrt{2\lambda};\\
				\dfrac{x^2}{2\lambda}, &\text{ if }|x|<\sqrt{2\lambda}.
			\end{cases}
		\end{align*}
	\end{enumerate}
\end{example}
\proof
\ref{example:fail coin::hull}: To see $e_\lambda f=e_\lambda g$, by
\cite[Example 1.44]{rockafellar_variational_1998} or
Lemma~\ref{lem:prox of proximal hull}, it suffices to show $h_\lambda f=g$.
% it is easy  in view of \cite[Example 1.44]{rockafellar_variational_1998}.
%We justify $h_\lambda f=g$.
To this end, observe that for every $x\in\R$
\begin{align*}
\clconv(f+\lambda^{-1}j)(x)
	&\leq
	\inf_{0<\mu<1}
	\{\mu(f+\lambda^{-1}j)(\mu^{-1}x)\}
	=
	\inf_{0<\mu<1}
	\{
	\mu
	[1+(2\lambda)^{-1}\mu^{-2}x^2]
	\} \\
	&=
	\inf_{0<\mu<1}
	\left\{
	\mu+\frac{x^2}{2\lambda}\mu^{-1}
	\right\},
\end{align*}
where the first inequality holds because $x=(1-\mu)0+\mu(\mu^{-1}x)$.
We claim that
\begin{equation}\label{claim}
	\inf_{0<\mu<1}
\left\{
\mu+\frac{x^2}{2\lambda}\mu^{-1}
\right\}
	=(g+\lambda^{-1}j)(x).
\end{equation}
For simplicity, let $a :=x^2/2\lambda$.
Then easy calculus shows that $|a|\leq 1\Leftrightarrow |x|\leq \sqrt{2\lambda}$ and
\begin{align*}
	\inf_{0<\mu<1}
	\{
	\mu+a\mu^{-1}
	\}
	=
	\begin{cases}
		2\sqrt{a},&\text{ if }|a|\leq1;\\
		1+a,&\text{ if }|a|>1,
	\end{cases}
\end{align*}
from which (\ref{claim}) readily follows.
Moreover, the convex function $g+\lambda^{-1}j$ is dominated by $f+\lambda^{-1}j$.
To see this, note that $(g+\lambda^{-1}j)(x)=1+\lambda^{-1}j(x)=(f+\lambda^{-1}j)(x)$ for $|x|>\sqrt{2\lambda}$ and
$$
(\forall |x|\leq\sqrt{2\lambda})~
(g+\lambda^{-1}j)(x)
=\sqrt{2\lambda^{-1}}|x|
\leq
\frac{1}{2}
\left(
\lambda (2\lambda^{-1})+\lambda^{-1}x^2
\right)
=
(f+\lambda^{-1}j)(x),
$$
where the inequality follows from the Fenchel-Young inequality.
Altogether $\clconv(f+\lambda^{-1}j)=g+\lambda^{-1}j$, and therefore \cite[Example 11.26(c)]{rockafellar_variational_1998} entails $h_\lambda f=g$.

\ref{example:fail coin::prox}:
The formula of $P_\lambda f$ is well-known; see, e.g., \cite[Example 4.2]{wang23level}.
Fact \ref{fact: level and Fenchel sub} yields
\begin{align}\label{lp sub}
\partial_p^\lambda g(x)
=
\partial_F(g+\lambda^{-1}j)(x)-\lambda^{-1}x
=
\begin{cases}
[-\sqrt{2\lambda^{-1}},\sqrt{2\lambda^{-1}}],&\text{ if }x=0;\\
\sqrt{2\lambda^{-1}}\sgn(x)-\lambda^{-1}x, &\text{ if }0<|x|\leq\sqrt{2\lambda};\\
0, &\text{ otherwise},
\end{cases}
\end{align}
and moreover
\begin{align*}
\left(
\Id+\lambda\partial_p^\lambda g
\right)(x)
=
\begin{cases}
	[-\sqrt{2\lambda},\sqrt{2\lambda}],&\text{ if }x=0;\\
	\sqrt{2\lambda}\sgn(x), &\text{ if }0<|x|\leq\sqrt{2\lambda};\\
	x, &\text{ otherwise}.
\end{cases}
\end{align*}
The claimed formula of $P_\lambda g$ follows from inverting the above and Fact \ref{fact: prox identity}, after which simplifying the identity $e_\lambda g(x)=g(P_{\lambda }g(x))+( P_{\lambda }g(x)-x)^2/2\lambda$ completes the proof.
\endproof

\begin{figure}[h!]\label{fig:prox}
	\begin{subfigure}[b]{0.3\textwidth}
		\centering
		\resizebox{\linewidth}{!}
		{
			\begin{tikzpicture}
				\begin{axis}[xmin = -2, xmax = 2,ymin = -2, ymax = 2,
					legend style=
					{
						legend columns=2,
						anchor=south,
						at={(0.5,1)}
					},
					xtick={-2,-1,0,1,2},
					ytick={-2,-1,0,1,2},
				xticklabels={, $-\sqrt{2}$,0,$\sqrt{2}$, $2\sqrt{2}$},
				yticklabels={$-2\sqrt{2}$, $-\sqrt{2}$,0,$\sqrt{2}$, $2\sqrt{2}$}
				]
					\addplot[domain =-2:-1,samples = 10,smooth,ultra thick,blue]{x};
					\addplot[domain =1:2,samples = 10,smooth,ultra thick,blue] {x};
					\addplot[domain =-1:1,samples = 10,smooth,ultra thick,blue]{0};
				
				\end{axis}
			\end{tikzpicture}
		}
		\caption{$P_\lambda f$}
	\end{subfigure}
	\begin{subfigure}[b]{0.3\textwidth}
		\centering
		\resizebox{\linewidth}{!}{
			\begin{tikzpicture}
				\begin{axis}[xmin = -2, xmax = 2,ymin = -2, ymax = 2,
					legend style={
						legend columns=2,
						anchor=south,
						at={(0.5,1)}
					},
					xtick={-2,-1,0,1,2},
					ytick={-2,-1,0,1,2},
					xticklabels={, $-\sqrt{2}$,0,$\sqrt{2}$, $2\sqrt{2}$},
					yticklabels={$-2\sqrt{2}$, $-\sqrt{2}$,0,$\sqrt{2}$, $2\sqrt{2}$}
					]
					
					\addplot[domain =-2:-1,samples = 10,smooth,ultra thick,blue]{x};
					\addplot[domain =1:2,samples = 10,smooth,ultra thick,blue] {x};
					\addplot[domain =-1:1,samples = 10,smooth,ultra thick,blue]{0};
					\addplot[smooth,ultra thick,blue] coordinates {(1,0) (1, 1)};
					\addplot[smooth,ultra thick,blue] coordinates {(-1,0) (-1, -1)};
				\end{axis}
				
			\end{tikzpicture}
		}
		\caption{$P_\lambda g$}
	\end{subfigure}
	\begin{subfigure}[b]{0.3\textwidth}
	\centering
	\resizebox{\linewidth}{!}{
		\begin{tikzpicture}
			\begin{axis}[xmin = -2, xmax = 2,ymin = 0, ymax = 3,
				legend style={
					legend columns=2,
					anchor=south,
					at={(0.5,1)}
				},
				xtick={-2,-sqrt(2),0,sqrt(2),2},
				ytick={0,1,2,3},
				xticklabels={,$-\sqrt{2}$,0,$\sqrt{2}$, $2\sqrt{2}$},
				yticklabels={$0$,$1$,$2$,$2\sqrt{2}$}
				]
				
				\addplot[domain =-2:-sqrt(2),samples = 10,smooth,ultra thick,blue]{1};
				\addplot[domain =sqrt(2):2,samples = 10,smooth,ultra thick,blue] {1};
				\addplot[domain =-sqrt(2):sqrt(2),samples = 10,smooth,ultra thick,blue]{x^2/2};
	
			\end{axis}
			
		\end{tikzpicture}
	}
	\caption{$e_\lambda g$}
\end{subfigure}
	\caption{Proximal operators and Moreau envelope in Example \ref{example:fail coin} with $\lambda=1$.}
\end{figure}

Nonetheless, we have an integration result of level proximal subdifferentials. At the core of its proof
is the integration technique for the Fenchel subdifferential by Benoist and Daniilidis \cite{benoist2002}
or Bachir, Daniilidis and Penot \cite{penot2002}.

\begin{theorem}[integration of level proximal subdifferential]\label{convexhull}
Let $f:\Rn\to\overline{\mathbb R}$ be proper, lsc and prox-bounded with thresholds $\lambda_{f}>0$.
Let $0<\lambda<\lambda_{f}$ and let $x_{0}\in\dom \partial_{p}^{\lambda}f$.
Define $\widehat{f}:\Rn\to\overline{\mathbb R}$ by
\begin{equation}\label{e:antiderivative}
x\mapsto f(x_{0})+\lambda^{-1}j(x_{0})+
\sup\left\{\sum_{i=0}^{k-1}\ip{x_{i}^*+\lambda^{-1}x_{i}}{x_{i+1}-x_{i}}+\ip{x_{k}^*+\lambda^{-1}
x_{k}}{x-x_{k}}
\right\}-\lambda^{-1}j(x),
\end{equation}
where the supremum is taken for $k\in\NN$,
all $x_{1},\ldots, x_{k}\in\dom \partial_p^{\lambda}f$ and
all $x_0^*\in\partial_{p}^{\lambda}f(x_{0}), \ldots, x_k^*\in\partial_{p}^{\lambda}
f(x_{k})$.
Then the following hold:
\begin{enumerate}
\item\label{i:integ1}
 $\widehat{f}=h_{\lambda}f$ the proximal hull of $f$, and
 $e_{\lambda}\widehat{f}=e_{\lambda}f;$
\item\label{i:integ2}
 If, in addition, $f$ is $\lambda$-proximal, i.e., $f+\lambda^{-1}j$ is convex,
then $\widehat{f}=f$.
\end{enumerate}
\end{theorem}
\proof
\ref{i:integ1}:
By Fact~\ref{fact: level and Fenchel sub},
$\partial_{F}(f+\lambda^{-1}j)(x)=\partial_{p}^{\lambda}f(x)+\lambda^{-1}x$.
Since $\dom \partial_{F}(f+\lambda^{-1}j)=\dom \partial_{p}^{\lambda}f\neq\varnothing$ by Proposition~\ref{thm: always nonempty}, it is possible to find $x_{i}\in\dom \partial_{p}^{\lambda}f$
for $i=0, \ldots, k$ so that \eqref{e:antiderivative} is well defined.
Furthermore,
$x_{i}\in\dom\partial_{F}(f+\lambda^{-1}j) \Leftrightarrow  x_{i}\in\dom \partial_p^{\lambda}f$,
and $\partial_{F}(f+\lambda^{-1}j)(x_{i})=\partial_{p}^{\lambda}f(x_{i})+\lambda^{-1}x_{i}$ for $i=0,\ldots, k$.
The function $f+\lambda^{-1}j$ is $1$-coercive (so epi-pointed) and
$\conv(f+\lambda^{-1}j)=\clconv(f+\lambda^{-1}j)$ by
Lemma~\ref{lem: proximal and closed cvx hull}\ref{lem: closed cvx hull}.
Employing \cite[Proposition 2.7]{penot2002}
or \cite[Theorem 3.5]{benoist2002}, we have
\begin{align}
(\forall x\in\Rn)\ & \conv (f+\lambda^{-1}j)(x) \nonumber\\
 & =f(x_{0})+\lambda^{-1}j(x_{0})+
\sup\left\{\sum_{i=0}^{k-1}\ip{x_{i}^*+\lambda^{-1}x_{i}}{x_{i+1}-x_{i}}+\ip{x_{k}^*
+\lambda^{-1}x_{k}}{x-x_{k}}
\right\},
\end{align}
where the supremum is taken for $k\in\NN$,
all $x_{1},\ldots, x_{k}\in\dom \partial_{F}(f+\lambda^{-1}j)$ and
all $x_0^*+\lambda^{-1}x_{0}\in\partial_{F}(f+\lambda^{-1}j)(x_{0}), \ldots, x_k^*
+\lambda^{-1}x_{k}\in\partial_{F}
(f+\lambda^{-1}j)(x_{k})$.
It follows that
 $\widehat{f}(x)=\conv(f+\lambda^{-1}j)(x)-\lambda^{-1}j(x)=h_{\lambda}(x).$
 Moreover, $e_{\lambda}\widehat{f}=e_{\lambda}[h_{\lambda}f]=e_{\lambda}f$
 by Lemma~\ref{lem:prox of proximal hull}\ref{i:rock:prox}.

\ref{i:integ2}: Immediate from \ref{i:integ1}.
\endproof

\begin{remark}
To put Theorem~\ref{convexhull} into perspective, for possibly nonconvex functions,
under mild conditions on $f$ while $\partial_{F}f$ determines the closed convex hull of $f$ uniquely up to a constant
\cite{benoist2002, penot2002},
$\partial_{p}^{\lambda}f$ determines the proximal hull of $f$ uniquely up to a constant.
See \cite{clarke93,zili2000,thibault1995,poliquin1991} for integration results
of other subdifferentials.
\end{remark}

%Besides  the coincidence results in Theorem~\ref{thm:coincidence}, another natural question to ask is
%\begin{center}
%\emph{Does $\partial_p^{\lambda_1}f(x)=\partial_p^{\lambda_2}f(x)\Leftrightarrow\lambda_1=\lambda_2$?}
%\end{center}
%One direction is certainly trivial.
%However, the opposite may not hold true, as the example below demonstrates.
%\begin{example}
%Let $f(x)=-\|x\|$ and $\lambda_1>\lambda_2>0$.
%Then $\partial_p^{\lambda_1}f(0)=\emptyset=\partial_p^{\lambda_2}f(0)$.
%\end{example}
%\proof
%See \cite[Example 4.4]{wang23level} for the level proximal subdifferential formula of $f$.
%\endproof

\section{Pointwise quadratic approximation (or Lipschitz smoothness)}\label{s:lips}

Functions with Lipschitz continuous gradients ($L$-smooth functions) are usually studied for convex functions, and they are used widely in optimization algorithms; see, e.g.,
\cite{BC, bauschke2017}, \cite[pages 110--112]{beck2017}, \cite{nesterov2018} and \cite{wang2022mirror}.
In this section, we propose the pointwise quadratic approximation (or Lipschitz smoothness)
for possibly nonconvex functions, which
opens the door to study Lipschitz smoothness on open sets. As shown below, the bridge to investigate
the pointwise quadratic approximation is
a two-sided condition utilizing the level proximal
subdifferentials $\partial_p^\lambda(\pm f)$. If the pointwise quadratic approximation property holds on an open set,
 it forces the function to have a gradient
which is Lipschitz on the open set.

Let $f :\Rn\rightarrow\bR$ be proper lsc and differentiable at $x\in\dom f$. We say
that $f$ \emph{admits a quadratic approximation} (or is \emph{Lipschitz smooth}) at $x$ if
\begin{equation}\label{e:lipsmooth}
(\exists L>0)(\forall y\in\Rn)\ |f(y)-f(x)-\ip{\nabla f(x)}{y-x}|\leq L j(y-x),
\end{equation}
and that $f$ is Lipschitz smooth on $U\subseteq\Rn$ if $f$ is Lipschitz smooth at
every $x\in U$.  If $f$ is an $L$-smooth function \cite[Chapter 5]{beck2017}, namely,
$(\forall x,y\in\Rn)\ \|\nabla f(x)-\nabla f(y)\|\leq L\|x-y\|$, then \eqref{e:lipsmooth} holds by the descent lemma
\cite[Lemma 5.7]{beck2017} applied to both $f, -f$. Amazingly, the converse is true when $f$ is convex and
\eqref{e:lipsmooth} holds over the entire space for both $x, y\in\Rn$, see, e.g., \cite[Theorem 5.8]{beck2017}.
It is natural to ask what happens if $f$ is possibly nonconvex and
\eqref{e:lipsmooth} holds only for $x\in U$, an open subset of $\Rn$.
This is the motivation for us to define the Lipschitz smoothness at a point.
We note that in \cite{fabian} Fabi\'an studied the Lipschitz smooth point of a convex function when
\eqref{e:lipsmooth} holds only locally for $y\in\Rn$.

\begin{lemma}\label{lem:descent to lip grad}
Let $U$ be a nonempty open subset in $\Rn$ and let $L>0$.
	Suppose that $g:\Rn\to {\mathbb{R}}$ satisfies
\begin{equation}\label{e:global:opt}
(\forall y\in\mathbb{R}^n)(\forall x\in U)\
	0\leq g(y)-g(x)-\left\langle \nabla g(x),y-x \right\rangle\leq L\left\|y-x\right\|^2.
\end{equation}
	Then $(2L)^{-1} \nabla g$ is firmly nonexpansive on $U$, i.e.,
	\begin{equation}\label{eq: proto descent lemma}
(\forall x,y\in U)\ (2L)^{-1}\left\|\nabla g(y)-\nabla g(x)\right\|^2\leq
		\left\langle \nabla g(x)-\nabla g(y),x-y \right\rangle.
	\end{equation}
\end{lemma}
\proof
Fix $x\in U$, and define $\phi:\Rn\rightarrow\bR: y\mapsto g(y)-\left\langle \nabla g(x),y \right\rangle$.
Then \eqref{e:global:opt} implies that $\phi$ attains global minimum at $x$, i.e.,
$
\phi(x)=\min_{z\in\mathbb{R}^n}\phi(z).
$
We claim that for every $y\in U$
\begin{equation}\label{e:otherone}
(\forall z\in\mathbb{R}^n)~
\phi(z)-\phi(y)-\left\langle \nabla\phi(y),z-y \right\rangle\leq L\left\|y-z\right\|^2,
\end{equation}
implying that
\begin{align}
	\phi(x)
	&=\min_{z\in\mathbb{R}^n}\phi(z)\leq
	\min_{z\in\mathbb{R}^n}
	\left[
	\phi(y)+\left\langle \nabla\phi(y),z-y \right\rangle+L\left\|z-y\right\|^2
	\right]\nonumber\\
	&\leq\min_{r\geq0}
	\left[
	\phi(y)-\left\|\nabla \phi(y)\right\|r+Lr^2
	\right]
	=\phi(y)-\frac{\left\|\nabla \phi(y)\right\|^2}{4L},\label{eq: bound}
\end{align}
where the last equality holds by minimizing the quadratic function $r\mapsto -\left\|\nabla \phi(y)\right\|^2r+Lr^2$.
To see \eqref{e:otherone}, by \eqref{e:global:opt}, for every $y\in U$ we have
\begin{align*}
(\forall z\in\mathbb{R}^n)~
	&\phi(z)-\phi(y)-\left\langle \nabla\phi(y),z-y  \right\rangle\\
	& =
	\left(
	g(z)-\left\langle \nabla g(x),z \right\rangle
	\right)
	-
	\left(
	g(y)-\left\langle \nabla g(x),y \right\rangle
	\right)
	-
	\left\langle \nabla g(y)-\nabla g(x),z-y \right\rangle\\
	& =
	g(z)-g(y)-\left\langle \nabla g(y),z-y \right\rangle\leq L\left\|z-y\right\|^2,
\end{align*}
as claimed.
Therefore (\ref{eq: bound}) implies
$$
(\forall y\in U)~
g(x)-g(y)-\left\langle \nabla g(x),x-y  \right\rangle=\phi(x)-\phi(y)\leq-\frac{\left\|\nabla g(y)-\nabla g(x)\right\|^2}{4L},
$$
owing to the identity $\nabla\phi(y)=\nabla g(y)-\nabla g(x)$ and the definition of $\phi$.
Recall that $x\in U$ is arbitrary. Then
\begin{equation}\label{eq1}
	(\forall x,y\in U)~
	g(x)+\left\langle \nabla g(x),y-x \right\rangle+\frac{\left\|\nabla g(y)-\nabla g(x)\right\|^2}{4L}
	\leq
	g(y).
\end{equation}
By symmetry
\begin{equation}\label{eq2}
	(\forall x,y\in U)~
	g(y)+\left\langle \nabla g(y),x-y \right\rangle+\frac{\left\|\nabla g(y)-\nabla g(x)\right\|^2}{4L}\leq g(x).
\end{equation}
Adding (\ref{eq1}) and (\ref{eq2}) yields
$$
(\forall x,y\in U)~
\left\langle \nabla g(x)-\nabla g(y),y-x \right\rangle+\dfrac{\left\|\nabla g(y)-\nabla g(x)\right\|^2}{2L}\leq0,
$$
which is (\ref{eq: proto descent lemma}).
\endproof

\begin{remark} The proof of Lemma~\ref{lem:descent to lip grad} was inspired by that of \cite[Theorem 2.15]{nesterov2018}, in which Nesterov assumed that
$g:\Rn\rightarrow\R$ is convex and that \eqref{e:global:opt} holds for every $x,y\in\Rn$.
\end{remark}

%\begin{lemma}\label{lem:coincide cvxhull}
%Let $g:\Rn\to\overline{\R}$ be proper, and let $O\subseteq\dom g$ be a nonempty open convex set. Then $g=\clconv g$ on $O$ if and only if $g$ is convex and $\partial_Fg\neq\emptyset$  on $O$.
%\end{lemma}
%\proof Suppose that $g=\clconv g$ on $O$. If $g$ was not convex, then we would have $x_0,x_1\in O$ and $\lambda\in(0,1)$ such that
%\(
%(\clconv g)(x_\lambda)=g(x_\lambda) >(1-\lambda)g(x_0)+\lambda g(x_1)=(1-\lambda)(\clconv g)(x_0)+\lambda (\clconv g)(x_1),
%\)
%where $x_\lambda=(1-\lambda)x_0+\lambda x_1\in O$, which is absurd.
%Furthermore, note that $O\subseteq\inte\dom(\clconv g)$.
%Then Remark \ref{rem:clconv hull subd} and \cite[Proposition 16.27]{BC} enforce $\partial_F g=\partial(\clconv g)\neq\emptyset$ on $O$.
%
%Conversely, it suffices to show that $\clconv g\geq g$ on $O$.
%By an appeal to the assumption $\partial_Fg\neq\emptyset$ on $O$, to each $x\in O$ corresponds affine minorants of $g$ taking the form $y\mapsto f(x)+\ip{v}{y-x}$ for some $v\in\partial_F g(x)$, which means $(\clconv g)(x)\geq g(x)$ owing to the fact that $\clconv g$ is the supremum of all affine minorants of $g$; see, e.g., \cite[Proposition 2.5.2]{Hiriart-Urruty1993ConvexAnalysis}. \endproof
Armed with Lemma~\ref{lem:descent to lip grad}, we are able to study Lipschitz
smoothness locally.
\begin{theorem}[Lipschitz smooth on an open set]\label{thm:Lip}
Let $f:\mathbb{R}^n\to \mathbb{R}$ and let $U\subseteq \mathbb{R}^n$ be a nonempty open convex set.
Suppose that both $f$ and $-f$ are prox-bounded with thresholds $\lambda_{f}, \lambda_{-f}>0$ respectively, and
that $0<1/L<\min\{\lambda_{f},\lambda_{-f}\}$.
%\new{be such that $L>1/\lambda_f$}.
Then the following are equivalent:
	\begin{enumerate}
		\item\label{thm:Lip::two-sided sub} $(\forall x\in U)\ \partial_p^{1/L}\left(\pm f\right)(x)\neq\varnothing$;
		\item\label{thm:Lip::two-sided cvx}
%$Lj\pm f$ are convex on $O$ and
$(\forall x\in U)\ \partial_F(Lj\pm f)(x)\neq\varnothing$;
		\item\label{thm:Lip::descent}  The function $f$ is differentiable on $U$, and
\begin{equation}\label{eq:local descent lemma}
	(\forall y\in\mathbb{R}^n)(\forall x\in U)\ |f(y)-f(x)-\left\langle \nabla f(x),y-x \right\rangle|\leq\frac{L}{2}\left\|y-x\right\|^2.
\end{equation}
	\end{enumerate}
When one of the above holds, $\nabla f$ is $L$-Lipschitz on $U$.
\end{theorem}

\proof
Since $0<1/L<\min\{\lambda_{f},\lambda_{-f}\}$, by Proposition~\ref{thm: always nonempty} we have $\dom\partial_p^{1/L}\left(\pm f\right)\neq\varnothing$.

``\ref{thm:Lip::two-sided sub}$\Rightarrow$\ref{thm:Lip::two-sided cvx}'':
Apply Fact~\ref{fact: level and Fenchel sub}.

%By Theorem \ref{thm:domain} and Lemma \ref{lem: proximal and closed cvx hull}, $\clconv(Lj\pm f)(x)=\conv(Lj\pm f)(x)=(Lj\pm f)(x)$ for every $x\in O$.
%In turn invoking Lemma \ref{lem:coincide cvxhull} with $g=Lj\pm f$ justifies \ref{thm:Lip::two-sided cvx}.

``\ref{thm:Lip::two-sided cvx}$\Rightarrow$\ref{thm:Lip::descent}'':
For every $x\in U$, the assumption $\partial_F(Lj\pm f)(x)\neq\varnothing$ and subdifferential calculus yield
\begin{equation}\label{eq:temp sub inclusion}
	\hat \partial (\pm f)(x)=\hat\partial\left(Lj\pm f-Lj\right)=\hat \partial\left(Lj\pm f\right)(x)-Lx\supseteq\partial_F\left(Lj\pm f\right)(x)-Lx\neq\varnothing,
\end{equation}
entailing $f$ to be differentiable at $x$. Consequently (\ref{eq:temp sub inclusion}) gives
$$\varnothing\neq\partial_F\left(Lj\pm f\right)(x)\subseteq Lx\pm\nabla f(x)\Rightarrow \partial_F\left(Lj\pm f\right)(x)=Lx\pm\nabla f(x),$$
from which we conclude that  $\partial_p^{1/L}\left(\pm f\right)(x)=\partial_F\left(Lj\pm f\right)-Lx=\pm\nabla f(x)$ by Fact \ref{fact: level and Fenchel sub}. Then
$$
(\forall y\in\mathbb{R}^n)~
f(y)\geq f(x)+\left\langle \nabla f(x),y-x \right\rangle-\frac{L}{2}\left\|y-x\right\|^2, \text{ and }
$$
$$
(\forall y\in\mathbb{R}^n)~
-f(y)\geq -f(x)+\left\langle -\nabla f(x),y-x \right\rangle-\frac{L}{2}\left\|y-x\right\|^2,
$$
together furnish \eqref{eq:local descent lemma}.

``\ref{thm:Lip::descent}$\Rightarrow$\ref{thm:Lip::two-sided sub}'':
Inequality (\ref{eq:local descent lemma}) yields $-\nabla f(x)\in\partial_p^{1/L}(-f)(x)$ and $\nabla f(x)\in\partial_p^{1/L}f(x)$ for $x\in U$,
thus justifies \ref{thm:Lip::two-sided sub}.

Finally we show that $\nabla f$ is $L$-Lipschitz on $U$ provided \ref{thm:Lip::descent} holds.
Define $g :=f+Lj$.
Then \ref{thm:Lip::descent} amounts to
$
(\forall y\in\mathbb{R}^n)(\forall x\in U)\
0\leq g(y)-g(x)-\left\langle \nabla g(x),y-x \right\rangle\leq L\left\|x-y\right\|^2,
$
which, by Lemma~\ref{lem:descent to lip grad}, gives that
for $x,y\in U$,
\begin{align}
	\frac
	{
		\norm{\nabla f(x)-\nabla f(y)+L(x-y)}^2
	}
	{2L}
	&=
	\frac{\left\|\nabla g(y)-\nabla g(x)\right\|^2}{2L}\leq
	\left\langle \nabla g(x)-\nabla g(y),x-y \right\rangle\nonumber\\
	&=
	\ip{\nabla f(x)-\nabla f(y)+L(x-y)}{x-y},\label{eq: temp}
\end{align}
where the inequality holds by (\ref{eq: proto descent lemma}).
From \eqref{eq: temp} we derive
\begin{align*}
&\norm{\nabla f(x)-\nabla f(y)}^2 +2L\ip{\nabla f(x)- \nabla f(y)}{x-y}+L^2\norm{x-y}^2 \\
& \leq 2L \ip{\nabla f(x)-\nabla f(y)}{x-y}+2L^2\norm{x-y}^2,
\end{align*}
which reduces to $\norm{\nabla f(x)-\nabla f(y)}^2 \leq L^2\norm{x-y}^2$;
thus $(\forall x, y\in U)\ \norm{\nabla f(x)-\nabla f(y)}\leq L\norm{x-y}$, as required.
\endproof

\begin{example} Having a Lipschitz gradient on an open convex set $U$ cannot guarantee the existence of level-proximal
subdifferential on $U$, as shown by the following two simple functions.
\begin{enumerate}
\item Let $f(x)=1$ for $x\neq0$ and $f(0)=0$. Then $\nabla f=0$ so trivially $1$-Lipschitz on $(0,\sqrt{2})$, however $\partial_p^1 f(x)=\varnothing$ on the same interval \cite[Example 4.2]{wang23level}.
\item Let $f(x)=-|x|$. Then $\nabla f=-1$ so trivially $1$-Lipschitz on $(0,1)$, however
$\partial_{p}^1f(x)=\varnothing$ on the same interval \cite[Example 4.4]{wang23level}.
\end{enumerate}
\end{example}

{
	\begin{corollary}
Let $f:\mathbb{R}^n\to \mathbb{R}$.
Suppose that both $f$ and $-f$ are prox-bounded with thresholds
$\lambda_{f}, \lambda_{-f}>0$ respectively, and
that $0<1/L<\min\{\lambda_{f},\lambda_{-f}\}$.
%	Let $f:\mathbb{R}^n\to\overline{\mathbb{R}}$ be a proper lsc function,
%and prox-bounded with threshold $\lambda_f$,
%and let $L>0$.
% be such that $L>1/\lambda_f$.
Then the following are equivelent:
\begin{enumerate}
\item
$\partial_p^{1/L}\left(\pm f\right)\neq\varnothing$ on $\Rn$;
\item $Lj\pm f$ are convex functions with domain $\Rn$;
\item $f$ is differentiable on $\Rn$, and has a Lipschitz gradient with constant $L>0$.
\end{enumerate}
\end{corollary}
}

\section*{Acknowledgements}
The authors sincerely thank both anonymous referees for their insightful comments and constructive suggestions that
improved the quality of this paper.
The research of HL and XY was partially supported by
the NSF Grants of
China and Chongqing (11991024, 11771064, cstc2021jcyj-msx300,
CQYC20210309536, 20A110029, 12261160365).
The research of XW and ZW was partially supported by
the NSERC Discovery Grant of Canada.

\bibliographystyle{siam}
\bibliography{ref}

\end{document}